\documentclass[10pt,a4paper]{article}
\usepackage[latin,english]{babel}
\usepackage[utf8]{inputenc}
\usepackage[T1]{fontenc}
\usepackage{eufrak}
\usepackage{geometry}
\usepackage{amsmath}
\usepackage{amssymb}
\usepackage{amsbsy}
\usepackage{amsthm}
\usepackage{makeidx}
\usepackage{mathtools}
\usepackage{mathabx, mathrsfs, dsfont}
\usepackage{cases}
\usepackage{braket}
\usepackage[toc,page]{appendix}
\usepackage{dirtytalk}
\usepackage{pdfsync}
\usepackage{graphicx}
\graphicspath{ {./Figures/} } 
\usepackage{epstopdf}
\usepackage{todonotes}
\usepackage{fancyhdr}
\usepackage{math}
\usepackage{cite}
\usepackage{color}
\usepackage[dvipsnames]{xcolor}
\usepackage{tikz}
\usepackage{subcaption}
\usepackage{bbm}
\usepackage{authblk}
\usepackage[colorlinks=true, pdfstartview=FitV, urlcolor=blue, citecolor=red, linkcolor=blue,pdfencoding=unicode,unicode=true,psdextra]{hyperref}
\usepackage{marginnote}

\usepackage{dsfont}
\mathtoolsset{showonlyrefs}

\usetikzlibrary{matrix,decorations.markings}
\usetikzlibrary{decorations.markings}
\usetikzlibrary{quotes,angles}

\usepackage{esint}

\newcommand{\bbN}{{\mathbb{N}}}
\newcommand{\bbC}{{\mathbb{C}}}

\newcommand{\bbP}{{\mathbb{P}}}

\newcommand{\Tr}{\mathrm{Tr}}

\newcommand{\Airy}{\textnormal{Airy}}

\def\pmtwo#1#2#3#4{\left( \begin{array}{cc}#1&#2\\#3&#4\end{array}\right)}

\newtheorem{RHP}[theorem]{RHP}

\numberwithin{equation}{section}

\title{On the asymptotic distribution of the largest eigenvalue of unbalanced complex Jacobi ensembles}
\author[1]{G. Didier}
\author[2]{A. Gkogkou}
\author[3]{J. Lopez Santander}
\author[1]{K. D. T-R McLaughlin}
\affil[1]{Department of Mathematics, Tulane University, New Orleans, LA}
\affil[2]{Department of Applied Mathematics, University of Washington, Seattle, WA}
\affil[3]{Department of Mathematics, University of Manitoba, Winnipeg, MB}

\date{\today}

\begin{document}
	\maketitle
 \begin{abstract}

Let ${\mathbf A} \sim \bbC {\mathcal W}_p(I,m)$ and ${\mathbf B} \sim \bbC {\mathcal W}_p(I,n)$ be two independent complex Wishart random matrices. In this work, we investigate the asymptotic distribution of the largest eigenvalue $\theta_{\max}$ of the random matrix ${\mathbf J} = ({\mathbf A}+{\mathbf B})^{-1} {\mathbf B}$. The high-dimensional regime is potentially unbalanced in the sense that, as $p \to \infty$, $m$ may grow faster than $p$, whereas $n$ grows proportionally to $p$. The asymptotic behavior of $\theta_{\max}$ is characterized by means of a Riemann-Hilbert (RH) problem for Jacobi polynomials with modified Jacobi weights. We establish the type 2 Tracy-Widom asymptotic distribution of $\theta_{\max}$. This provides an alternative, RH-based, proof of the same claim in Johnstone \cite{johnstone:2008} for the balanced regime, and establishes a conjecture in Jiang \cite{jiang:2009} for the unbalanced regime. In addition, we obtain an explicit correction of order \(p^{-2/3}\) to the Tracy–Widom approximation, with a remainder of order \(\mathcal O(p^{-1})\)
, and show the exponential decay, as a function of $p$, of the upper tail of the distribution of $\theta_{\max}$.
 \end{abstract}

\section{Introduction}\label{section Introduction}

The \textit{$\beta$-Jacobi ensemble} is the family of probability density functions (p.d.f.s)
\begin{equation}\label{e:joint_eigenvalue_pdf_beta}
f_{p}(\boldsymbol{\theta})={\mathcal Z}_{p,\beta}\cdot \prod_{i=1}^p (1-\theta_i)^{c_1\beta/2}\theta_i^{c_2 \beta/2}\prod_{i<j}|\theta_i-\theta_j|^{\beta}= {\mathcal Z}_{p,\beta}\cdot \prod_{i=1}^p w_p(\theta_i) \prod_{i<j}|\theta_i-\theta_j|^{\beta}.
\end{equation}
In \eqref{e:joint_eigenvalue_pdf_beta}, $\beta > 0$, $c_1,c_2 > -2/\beta$, and ${\mathcal Z}_{p,\beta}$ is a normalization constant. Also, the so-named \textit{weight function} is defined as
 \begin{equation}\label{eq:eqweight1}
    w_p(x) = (1-x)^{c_1 \beta/ 2}x^{c_2\beta/2}, \quad x \in [0,1].
 \end{equation} 
As originally put forth in Dyson \cite{dyson:1962-I,dyson:1962-II,dyson:1962-III},  the Jacobi ensemble arises in Statistical Mechanics as a model for a so-called \textit{log (Coulomb) gas}  whose potential establishes $\beta$-logarithmic interactions among the charged particles and further confines them to the interval $[0,1]$ (Forrester \cite{forrester:2010}). Now let
\begin{equation}\label{e:X,Y}
{\mathbf X} \textnormal{ and } {\mathbf Y}
\end{equation}
be, respectively, $m \times p$ and $n \times p$ random matrices with i.i.d.\ standard Gaussian entries, where $m,n,p \in \bbN$, $m \geq p$. Also, let ${\mathbf A} = {\mathbf X}^*{\mathbf X}$ and ${\mathbf B} = {\mathbf Y}^*{\mathbf Y}$, where $*$ denotes the adjoint matrix. In \eqref{e:joint_eigenvalue_pdf_beta}, the particular choice of parameter values
\begin{equation*}
    c_1 =m-p+1-\frac{2}{\beta}, \quad c_2 =n-p+1-\frac{2}{\beta}, \quad \beta = 1,2,4,
\end{equation*}
describes the p.d.f.\ of the eigenvalues of the random matrix
\begin{equation}\label{e:(A+B)^(-1)B}
{\mathbf J} = ({\mathbf A}+{\mathbf B})^{-1}{\mathbf B}
\end{equation}
(e.g., James \cite{james1964distributions}, Forrester et al.\ \cite{forrester:rahman:witte:2017}). The Jacobi ensemble is further called \textit{real} (\textit{orthogonal}, $\beta =1$), \textit{complex} (\textit{unitary}, $\beta=2$) or \textit{quaternionic} (\textit{symplectic}, $\beta=4$) depending on all entries of both ${\mathbf X}$ and ${\mathbf Y}$ in \eqref{e:X,Y} being real, complex or quaternionic (Gaussian) random variables, respectively. Random matrix models of the Jacobi ensemble for general $\beta > 0$ were provided in Killip and Nenciu \cite{killip:nenciu:2004} and Edelman and Sutton \cite{edelman:sutton:2008}.

Aside from the prolific role of the general $\beta$-Jacobi ensemble as a model in Mathematical Physics, the real Jacobi ensemble is pervasive in Multivariate Analysis. It emerges in several techniques such as MANOVA, canonical correlation analysis, multivariate regression, covariance testing, and signal detection (Anderson \cite{anderson:2003}, Muirhead  \cite{muirhead:2009}). In turn, the complex Jacobi ensemble naturally appears in Quantum Mechanics, e.g., in chaotic scattering, conductance in mesoscopic systems, and statistical properties of periodically driven quantum systems (Bl\"{u}mel and Smilansky \cite{blumel:smilansky:1988}, Hakke \cite{haake:1991}, Beenakker \cite{beenakker:1997}); see also Example \ref{ex:mesoscopic}. This can be viewed as stemming from the fact that \eqref{e:(A+B)^(-1)B} is equal in distribution to ${\mathbf U}^*{\mathbf U}$, where ${\mathbf U}$ is a truncated Haar random matrix (Collins \cite{collins:2003}, Forrester \cite{forrester:2006}; see also {\.Z}yczkowski and Sommers\cite{zyczkowski:sommers:2000}, Petz and R{\'e}ffy \cite{petz:reffy:2005}, Khoruzhenko et al.\ \cite{khoruzhenko:sommers:zyczkowski:2010}, Dong et al.\ \cite{dong:jiang:li:2012}, Meckes and Stewart \cite{meckes:stewart:2019}). Likewise, the complex Jacobi ensemble is naturally relevant in the Fourier analysis of time series, which hinges upon dealing with complex-valued random vectors (see Example \ref{ex:time_series}).

Due to its great importance, the $\beta$-Jacobi ensemble has motivated the development of a rich literature covering properties of various instances (values of $\beta$) of the ensemble. For example, universality for eigenvalue correlations was established in Kuijlaars and Vanlessen \cite{kuijlaars:vanlessen:2002}. In Dette and Nagel \cite{dette:nagel:2009} and Huang and Harrow \cite{huang:harrow:2024}, the eigenvalues were approximated by roots of Jacobi polynomials. In Dumitriu and Paquette \cite{dumitriu:paquette:2012}, the global fluctuations of linear statistics were characterized. In Holcomb and Moreno Flores \cite{holcomb:morenoflores:2012}, it was shown that the limiting processes of the $\beta$-Jacobi ensemble at the soft and hard edges, respectively, are the stochastic Airy and Bessel point processes. In Nagel \cite{nagel:2014}, the bulk behavior was characterized when one or both parameters $m, n$ grow
faster than the dimension $p$. A central limit theorem for the characteristic polynomial was established in Augeri et al.\ \cite{augeri:butez:zeitouni:2023}.

In many cases -- e.g., in hypothesis tests for the equality of covariance matrices --, a key role is played by the \textit{largest} eigenvalue of \eqref{e:(A+B)^(-1)B}. Nevertheless, its exact distribution is rather difficult to calculate. Over the last many decades, this led to the development of extensive tables or software (e.g., Nanda \cite{nanda:1948}, Foster and Rees \cite{foster:rees:1957}, Chen \cite{chen:2002,chen:2003}), as well as to several fundamental theoretical results. For example, in Dumitriu and Koev \cite{dumitriu:koev:2008}, the exact distributions of the largest eigenvalues of $\beta$-Jacobi ensembles were developed. In Ramli et al.\ \cite{ramli:katzav:castillo:2012}, the Coulomb gas approach was applied in deriving the rate functions for large fluctuations of the two extreme eigenvalues. In Moreno-Pozas et al.\ \cite{moreno:morales-jimenez:mckay:2019},
exact extreme eigenvalue distributions were obtained. Asymptotic formulas for the extreme eigenvalues in terms of hypergeometric functions were constructed in Winn \cite{winn:2024}. Limiting covariances for extreme eigenvalues at zero temperature ($\beta \rightarrow \infty$) were established in Hermann and Voit \cite{hermann:voit:2025}.

In the seminal papers Johnstone \cite{johnstone:2008,johnstone:2009} and Jiang \cite{jiang:2009,jiang:2013}, characterizations of the distribution of the largest eigenvalue(s) of $\beta$-Jacobi ensembles were provided for multiple values of $\beta > 0$ in terms of either random operators or of the so-named \textit{Tracy-Widom} family of distributions. Hereinafter we focus on the complex case ($\beta =2$). More precisely, recall that the celebrated \textit{type 2 Tracy–Widom law} ($\mathrm{TW}_2$) is defined by the c.d.f. 
\begin{equation}\label{e:F2=exp(integ)}
    F_2(x)=\exp\left(-\int_{x}^{\infty}(y-x)q^2(y)dy\right)
\end{equation}
(Tracy and Widom \cite{tracy1994level, tracy2009distributions}). In \eqref{e:F2=exp(integ)}, the function $q$ is the unique solution to the Painlevé II equation
\begin{equation}\label{eq. PII}
    q''(x)=x\,q(x)+2q(x)^3,
\end{equation}
with boundary condition
\begin{equation*}
    q(x)\sim \mathrm{Ai}(x)\quad \text{as}\quad x\rightarrow\infty.
\end{equation*} 
In Johnstone \cite{johnstone:2008}, it was established that -- after centering and rescaling -- the distribution of the largest eigenvalue of the complex Jacobi ensemble asymptotically approaches $\mathrm{TW}_2$, where the error is of order $\mathcal{O}(p^{-2/3})$. In that paper, the asymptotic regime is such that \textit{both} parameters $m$ and $n$ grow proportionally with $p$, namely,
\begin{equation}\label{eq:Johnstone_conditions}
    \lim_{p\rightarrow\infty}\frac{\min(p,n)}{m+n}>0, 
    \hspace{0.5cm} 
    \lim_{p\rightarrow\infty}\frac{p}{m}<1.
\end{equation}
Borrowing terminology from the area of Design of Experiments (e.g., Box et al.\ \cite{box:hunter:hunter:2005}), we refer to \eqref{eq:Johnstone_conditions} as the \textit{balanced} regime. 

Nevertheless, real data is often \textit{unbalanced} (Christensen \cite{christensen:2018}). For example, in multiple fields of application individually recorded time series are inevitably asynchronous or display unequal lengths. By comparison to classical time series assumptions (e.g., Brockwell and Davis \cite{brockwell:davis:1991}), this poses extra challenges for both modeling and classification. Accordingly, there is a vast literature on the subject in the realms of Mathematical Statistics and Machine Learning (e.g., Hayashi and Yoshida \cite{hayashi:yoshida:2005}, Bokil et al.\ \cite{bokil:pupura:schoffelen:thomson:mitra:2007}, T{\'o}th and Kert{\'e}sz \cite{toth:kertesz:2009}, Cuturi \cite{cuturi:2011}, Ruiz et al.\ \cite{ruiz:flynn:large:middlehurst:bagnall:2021}). In addition, the balanced regime \eqref{eq:Johnstone_conditions} only holds \textit{asymptotically} -- hence, knowledge about unbalanced regimes naturally sheds new light on the accuracy of the approximation established for the balanced regime itself.

For these reasons, the present work is dedicated to the distribution of the largest eigenvalue of \eqref{e:(A+B)^(-1)B} in the (generally) \textit{unbalanced} regime where $m$ grows faster than linearly with $p$, while $n$ grows still proportionally to $p$. Specifically, we consider parameters satisfying
\begin{equation}\label{eq: parameters us}
    n-p= b_{\infty}\,p \quad \text{and} \quad m-p=a_{\infty}\,p^{1+\nu}, \quad \nu \geq 0,
\end{equation}
with $a_{\infty},b_{\infty}>0$. Under \eqref{eq: parameters us}, the balanced regime \eqref{eq:Johnstone_conditions} corresponds to the instance $\nu = 0$. In contrast with \eqref{eq:Johnstone_conditions}, relation \eqref{eq: parameters us} implies that
\begin{equation}\label{e:unbalanced_regime}
    \lim_{p \rightarrow\infty}\frac{\min(p,n)}{m+n}=\lim_{p\rightarrow\infty}\frac{p}{m}=0
\end{equation}
for the (unbalanced) instances where $\nu > 0$. Visually, under \eqref{e:unbalanced_regime} the data matrix $\mathbf{X}$ as in \eqref{e:X,Y} becomes increasingly tall, while the data matrix $\mathbf{Y}$ maintains consistent asymptotic aspect ratio.  

In order to study the asymptotic behavior of the largest eigenvalue of \eqref{e:(A+B)^(-1)B} in the scaling regime \eqref{eq: parameters us}, in this paper we combine two interrelated sets of tools: $(i)$ \textit{Fredholm determinants}; and $(ii)$ the \textit{Riemann-Hilbert approach} for the asymptotic analysis of orthogonal polynomials. 

In fact, it is well known that the statistics of eigenvalue spacings for several unitary random matrix ensembles can be described in terms of a Fredholm determinant of an associated integral operator (see Tracy and Widom \cite{tracy1994level, tracy1994fredholm}). To understand this connection, first note that, under \eqref{eq: parameters us}, the weight function \eqref{eq:eqweight1} takes the form
\begin{equation}\label{eq:eqweight}
    w_p(x) = x^{b_{\infty} \, p} (1-x)^{ a_{\infty} \, p^{1+\nu}}.
 \end{equation}
Now, for $j \in \bbN \cup \{0\} $ and $\kappa_j>0$, let $\mathbf{p}_j(x)=\kappa_jx^j+\ell .o.t.$ be the sequence of orthonormal polynomials with respect to $w_p(x)$, where $\ell.o.t.$ stands for ``lower-order terms". Also, let $K_p$ be the so-named \textit{reproducing kernel}
\begin{gather}\label{eq:def Kernel z}
    K_p(z_1,z_2) = \sqrt{w_p(z_1)w_p(z_2)}\,\,\sum_{\ell=0}^{p-1} \mathbf{p}_\ell(z_1)\mathbf{p}_\ell (z_2).
\end{gather}
Further define the associated \textit{integral operator} acting on $L^2[a,1]$ by means of the relation
\begin{gather}\label{eq:eqdefnoper}
\mathcal{K}_p f(z_1)=\int_a^{1} K_p(z_1, z_2) f(z_2) \, d z_2.
\end{gather}
For $\beta = 2$ and for any $p \in \bbN$, let $\theta_{\max} := \theta_{\max,p}$ be the largest eigenvalue of ${\mathbf J}$ as in \eqref{e:(A+B)^(-1)B}. Then, we can express
\begin{equation}\label{eq: prob lambda max}
\mathbb{P}(\theta_{\max}<a)=\det( \mathbb{I}-\mathcal{K}_p)_{L^2[a,1]}, \quad 0 < a < 1
\end{equation}
(e.g., Anderson et al.\ \cite{anderson:guionnet:zeitouni:2010}, Forrester \cite{forrester:2010}). The expression appearing on the right-hand side of \eqref{eq: prob lambda max} is the \textit{Fredholm determinant associated with $\mathcal{K}_p$}, where the determinant is itself defined by means of the series
\begin{equation}\label{def:Fredholm_Kp}
\det( \mathbb{I}-\mathcal{K}_p)_{L^2([a,1])} = 1 + \sum_{j=1}^{p}  \frac{(-1)^{j}}{j!} \int_{a}^{1} \cdots \int_{a}^{1} \det \begin{pmatrix}
K_p(z_1,z_1) & \cdots & K_p(z_1,z_j)\\
\vdots & \cdots & \vdots \\
K_p(z_j,z_1) & \cdots & K_p(z_j,z_j)
\end{pmatrix} \, d^j {\mathbf z}.
\end{equation}

Now note that, under \eqref{eq: parameters us}, the function \eqref{eq:eqweight1} approaches \eqref{eq:eqweight}. Then, the relations \eqref{eq:def Kernel z}--\eqref{eq: prob lambda max} indicate that the analysis of the distribution of the largest eigenvalue of \eqref{e:(A+B)^(-1)B} relies on \textit{establishing the asymptotics (as $p \rightarrow \infty$)} of the corresponding Jacobi polynomials (cf.\ Bleher and Its \cite{bleher:its:1999}, Deift et al.\ \cite{deift:kriecherbauer:mclaughlin:venakides:zhou:1999}, Bothner and Its \cite{bothner:its:2014}). Indeed, for this purpose we follow the celebrated methodology put forth in Fokas et al.\ \cite{fokas:its:kitaev:1992}, namely, the \textit{Riemann-Hilbert approach}. Notably, in Deift and Gioev \cite{deift2007universality} this approach was successfully used in the study of edge universality for the unitary, orthogonal and simplectic random matrix ensembles for a class of weights $w(x)=e^{-V(x)}$, where $V$ is a polynomial with even dominant power. The Riemann-Hilbert approach has now become broadly acknowledged as a prolific toolbox for establishing key results on the statistics of large random matrices (see, for instance, Vanlessen \cite{vanlessen:2007}, Claeys et al.\ \cite{claeys:kuijlaars:vanlessen:2008}, Bertola et al.\ \cite{bertola:lee:mo:2009}, Harnad \cite{harnad:2011}, Bertola and Bothner \cite{bertola:bothner:2015}, Claeys and Fahs \cite{claeys:fahs:2016}, Berestycki et al.\ \cite{berestycki:webb:wong:2018}, Piorkowski \cite{piorkowski:2021}, Little \cite{little:2024}).

In this paper, we establish three main results assuming the scaling regime \eqref{eq: parameters us}. First, we show that, after appropriate centering and scaling, the largest eigenvalue converges to the type 2 Tracy--Widom distribution. That is, there exist sequences $\mu_p$ and $\sigma_p$ such that  
\begin{equation}\label{eq:TW}
     \frac{\theta_{\max} - \mu_p}{\sigma_p} \xrightarrow{\;d\;} \mathrm{TW}_2.
\end{equation}
As a function of $p$, the sequences $\mu_p$ and $\sigma_p$ scale as $p^{-\nu}$ and $p^{-\nu-2/3}$, respectively. Explicit formulas for these quantities, together with related results, are provided in Section~\ref{subsec statement results}. Second, we sharpen the limit \eqref{eq:TW} by showing that
\begin{equation}\label{e:K-hat_p_expansion}
\widehat{K}_{p}(\xi_1,\xi_2) = K_{\textnormal{Airy}}(\xi_1,\xi_2) + p^{-2/3}g_p(\xi_1,\xi_2)+ {\mathcal O}(p^{-1}).
\end{equation}
In \eqref{e:K-hat_p_expansion}, $\widehat{K}_{p}$ is an appropriately rescaled version of the reproducing kernel \eqref{eq:def Kernel z} (see expressions \eqref{eq:def Kernel lambda} and \eqref{kernel local coordinates}), $K_{\textnormal{Airy}}$ is the Airy kernel (see \eqref{Eq. Def Airy kernel}), and $g_p$ is a \textit{known} function of both the Airy function (see \eqref{Eq. Def Airy function}) and its derivative. In other words, vis-\`{a}-vis \eqref{eq:TW}, relation \eqref{e:K-hat_p_expansion} reduces the $\textnormal{TW}_2$ approximation error order from ${\mathcal O}(p^{-2/3})$ to ${\mathcal O}(p^{-1})$.  Third, it is well known that the Tracy-Widom class is associated with third-order phase transitions between strong and weak coupling regimes (see Gross and Witten \cite{gross:witten:1980}, Nadal and Majumdar \cite{nadal:majumdar:2011}, Majumdar and Schehr \cite{majumdar:schehr:2014}, Mari\~{n}o \cite{marino:2014}). We contribute to this literature by further quantifying the rate of the upper tail deviation as a function of $p$. Namely, we show that, for some constant $k > 0$, \begin{equation}\label{e:upper_tail_intro}
\bbP(\theta_{\max} > a) = {\mathcal O}(e^{-kp}).
\end{equation} 
In other words, the widely observed exponential decay of the upper tail of the largest eigenvalue (e.g., Ben Arous et al.\ \cite{benarous:dembo:guionnet:2001}, Majumdar and Vergassola \cite{majumdar:vergassola:2009}) also takes place in the case of the unbalanced complex Jacobi ensemble (on the lower tail, see Remark \ref{r:lower_tail}).

It is informative to further reexamine our contributions in light of the seminal papers Johnstone \cite{johnstone:2008,johnstone:2009} and Jiang \cite{jiang:2009,jiang:2013} -- in the context of this work, when restricted to $\beta = 2$. First, for the case covered in Johnstone \cite{johnstone:2008} (i.e., $\nu = 0$ in \eqref{eq: parameters us}), the support of the limiting \textit{bulk} spectral measure remains \textit{fixed}. By contrast, $\nu > 0$ in \eqref{eq: parameters us} corresponds to the mathematically challenging case where this support \textit{shrinks}, with both endpoints approaching zero. Second, Johnstone \cite{johnstone:2008} uses a \textit{Liouville-Green approach} taken from Olver \cite{OLVER1974392} for the asymptotic analysis of Jacobi polynomials. By making use of the Riemann-Hilbert approach, our work provides an alternative proof of the main result in Johnstone \cite{johnstone:2008} for the balanced complex Jacobi ensemble. Third, our results provide evidence that the $\textnormal{TW}_2$ approximation to the distribution of the largest eigenvalue is rather robust with respect to several instances of unbalanced design. Hence, with a view towards applications, our results complement the quantile and robustness analysis of the Tracy-Widom distribution provided in Johnstone \cite{johnstone:2009} for the Jacobi ensemble. Fourth, the unbalanced regimes considered in Jiang \cite{jiang:2009,jiang:2013} generally amount to the particular case $n = o(\sqrt{m})$ and $p = o(\sqrt{m})$ in \eqref{eq: parameters us}. Consequently, our results confirm part of the conjecture put forth in Jiang \cite{jiang:2009}, p.\ 228, ``Future problems": namely, we generalize the rates of convergence used in Theorem 7 in that paper.  

The remainder of this paper is organized as follows. Section \ref{s:preliminaries} is dedicated to preliminaries. In Section \ref{subsec RHP intro}, for the sake of illustration we provide examples of both balanced and unbalanced instances of the complex Jacobi ensemble encountered in applications. In Section~\ref{subsec RHP intro}, we introduce the Riemann–Hilbert formulation for the Jacobi polynomials associated with the eigenvalue distribution of the Jacobi ensemble. In Section \ref{s:connecting_RH_top_eigenval}, for the reader's convenience we provide a broad discussion of the connection between the Riemann-Hilbert problem described in Section~\ref{subsec RHP intro} and the distribution of the top eigenvalue of the complex Jacobi ensemble. In Section~\ref{subsec statement results}, we state the main results of the paper, namely, Theorems \ref{thm: Probability Tracy-Widom} and \ref{prop:prop7}, and Corollary~\ref{mainthrm}. In Section~\ref{proofs}, we provide the proofs of the main results, together with some auxiliary results directly needed to establish them. The proofs of these auxiliary results, in turn, are deferred to Section~\ref{proofs auxiliary}. Section~\ref{section RH Analysis} develops the corresponding Riemann–Hilbert analysis and derives the asymptotic behavior of these polynomials. Section \ref{s:conclusion} contains conclusions and a discussion of open problems that this work leads to.  In addition, for the readers' convenience, Section \ref{s:technical_results} is dedicated to auxiliary technical results and Section \ref{s:complex_RVs} recaps basic results on complex-valued random vectors.

\section{Preliminaries}\label{s:preliminaries}

\subsection{The complex Jacobi ensemble: examples from applications}\label{subsec Examples}

For the sake of illustration, in this subsection we provide two examples -- one balanced, one unbalanced -- where the complex Jacobi ensemble is of direct interest in applications.

\begin{example}\label{ex:mesoscopic}
A \textit{mesoscopic device} is a nanoscale electronic device that blends quantum and classical behavior, with quantum effects playing a significant role (Blanter and B\"{u}ttiker  \cite{blanter:buttiker:2000}, Wisne et al.\ \cite{wisne_etal:2024}). In the contexts of chaotic scattering or quantum conductance ({\.Z}yczkowski and Sommers \cite{zyczkowski:sommers:2000}, Forrester \cite{forrester:2006}), consider a mesoscopic device coupled to two leads, each of which supports $N/2$ open channels. The process of scattering can be described by means of the so-called \textit{scattering matrix} of size $N$, i.e., a unitary matrix ${\mathbf S}$ which maps incoming waves to outgoing waves. More precisely, for two input vectors $I_1, I_2 \in \mathbb{C}^{N/2}$, and two output vectors $O_1, O_2 \in \mathbb{C}^{N/2}$, we write
$$
{\mathbf S} \left(\begin{array}{c}
I_1 \\
I_2
\end{array}\right)
= \left(\begin{array}{c}
O_1 \\
O_2
\end{array}\right).
$$
We can now recast 
$$
{\mathbf S} = 
\left(\begin{array}{cc}
{\mathbf R}_{11} & {\mathbf T}_{12} \\
{\mathbf T}_{21} & {\mathbf R}_{22} \\
\end{array}\right), \quad {\mathbf R}_{11}, {\mathbf R}_{22}, {\mathbf T}_{12}, {\mathbf T}_{21} \in {\mathcal M}\big( N/2,N/2,\mathbb{C} \big),
$$
where the submatrices ${\mathbf R}_{\cdot \cdot}$ and ${\mathbf T}_{\cdot \cdot}$ are called, respectively, the \textit{reflection} and \textit{transmission} matrices. In particular, they are identically distributed. Moreover, the eigenvalues of each random matrix ${\mathbf T}^*_{12}{\mathbf T}_{12} $, ${\mathbf T}^*_{21}{\mathbf T}_{21} $, ${\mathbf R}^*_{11}{\mathbf R}_{11}$ and ${\mathbf R}^*_{22}{\mathbf R}_{22}$ coincide, in law, with those of a (balanced) complex Jacobi ensemble (cf.\ Collins \cite{collins:2003}, Jiang \cite{jiang:2009} and Forrester \cite{forrester:2010}, Section 3.8).
\end{example}

\begin{example}\label{ex:time_series}
In classical time series analysis, the estimation of the spectral density $f_{{\mathcal X}}$ of a stationary stochastic process ${\mathcal X}$ is based on Fourier transforms (Brockwell and Davis \cite{brockwell:davis:1991}, chapter 10). Under assumptions, the complex Jacobi ensemble naturally emerges in this context. 

In fact, for $T \in \mathbb{N}$, let 
\begin{equation}\label{e:omega-j}
\omega_j = \frac{2\pi j}{T}, \quad j = 1, \hdots, T,
\end{equation}
be the \textit{Fourier frequencies}. Also, let ${\mathcal X } = \{\mathcal X_t\}_{t = 1,\hdots,T} \subseteq \mathbb{R}^p$ be a sequence of random vectors. The \textit{discrete Fourier transform} (DFT) of ${\mathcal X}$ is given by the mapping
\begin{equation}\label{e:DFT}
j = 1,\hdots,T \hspace{1mm}\mapsto \hspace{1mm}{\mathcal F}_{T}(\omega_j) := \frac{1}{\sqrt{T}} \sum^{T}_{t=1}\mathcal X_t e^{- \Im t \omega_j} \in \mathbb{C}^{p}.
\end{equation}
Suppose, in addition, that all entries in the sequence ${\mathcal X }$ are zero-mean, i.i.d.\ Gaussian random variables (in particular, $p = 1$). Note that, by Proposition \ref{p:DFT_is_complex_Gaussian}, over about $50\%$ of the Fourier frequencies \eqref{e:omega-j} the DFT of univariate Gaussian white noise is a sequence of i.i.d.\ complex Gaussian random variables. As a direct consequence of this fact in dimension $p =1$, for general dimension $p \in \mathbb{N}$ we have
\begin{equation}\label{e:DFTs_are_iid_C-Normal}
\big\{ {\mathcal F}_T(\omega_j)\}_{j = 1,\hdots,\lfloor T /2 \rfloor-1} \stackrel{\mathbf{i.i.d.}}\sim \mathbb{C}{\mathcal N}({\mathbf 0}_{p \times 1},\sigma^2 I_p, {\mathbf 0}_{p \times p}), \quad \sigma^2 > 0.
\end{equation}
Therefore, one can set up a test for the equality of noise spectra in the Fourier domain by analogy to the test for equality of covariance matrices in the real setting (see Johnstone \cite{johnstone:2008,johnstone:2009}).

In fact, let ${\mathcal X}=\{{\mathcal X}_{t}\}_{t = 1,\hdots,T_{{\mathcal X}}} \in \mathbb{R}^p $ and ${\mathcal Y}=\{{\mathcal Y}_{t}\}_{t = 1,\hdots,T_{{\mathcal Y}}} \in \mathbb{R}^p $ be two independent Gaussian white noise time series. Each series is further assumed to display equal variances across its $p$ entries and no instantaneous correlation. Set 
$m =  \lfloor T_{{\mathcal X}} /2 \rfloor-1$, $n =  \lfloor T_{{\mathcal Y}} /2 \rfloor-1$, and consider their respective $p$-variate DFTs
$\big\{ {\mathcal F}_{T,{\mathcal X}}(\omega_j)\}_{j  = 1,\hdots,m}$, $\big\{ {\mathcal F}_{T,{\mathcal Y}}(\omega_j)\}_{j  = 1,\hdots, n} \in \mathbb{C}^p.$ Define the null hypothesis
\begin{equation}\label{e:H0_identically_equal_spectra}
H_0: \forall \omega \in [0,\pi], \hspace{2mm} f_{{\mathcal X}}(\omega) =  f_{{\mathcal Y}}(\omega) \equiv \frac{\sigma^2}{2\pi}\hspace{0.5mm} I_p =: \Sigma \in {\mathcal S}_{> 0}(p,\mathbb{C}).
\end{equation}
In other words, we would like to test whether or not the spectral densities for the two populations are identically equal (they are in case ${\mathcal X}$ and ${\mathcal Y}$ have the same variance structure). Then, we can define
$$
{\mathbf X}= \left( \begin{array}{c}
{\mathcal F}_{T,{\mathcal X}}(\omega_1)^{\top} \\
\vdots \\
{\mathcal F}_{T,{\mathcal X}}(\omega_{m})^{\top}
\end{array}\right) \in {\mathcal M}\big(m,p,\mathbb{C}\big),
\quad
{\mathbf Y}= \left( \begin{array}{c}
{\mathcal F}_{T,{\mathcal Y}}(\omega_1)^{\top} \\
\vdots \\
{\mathcal F}_{T,{\mathcal Y}}(\omega_{n})^{\top}
\end{array}\right) \in {\mathcal M}\big(n,p,\mathbb{C}\big).
$$
Under the null hypothesis \eqref{e:H0_identically_equal_spectra}, it is easy to see that $\mathbf{X}^*\mathbf X \sim \mathbb{C} {\mathcal W}_p(m,\Sigma)$ and $\mathbf{Y}^* \mathbf{Y} \sim \mathbb{C} {\mathcal W}_p(n,\Sigma)$ as a consequence of \eqref{e:DFTs_are_iid_C-Normal}. Since, in addition, $\mathbf{X}^*\mathbf X$ and $\mathbf{Y}^* \mathbf{Y}$ are independent, it can be shown that
$$
\lambda_{\max}\big( (\mathbf{X}^*\mathbf X + \mathbf{Y}^* \mathbf{Y} )^{-1} \hspace{0.5mm} \mathbf{Y}^* \mathbf{Y} \big) \stackrel{d}= \theta_{\mathbb{C}}(p,m, n),
$$
where $\theta_{\mathbb{C}}(p,m,n)$ is the largest eigenvalue of the (unbalanced) complex Jacobi ensemble.
\end{example}

\subsection{The Riemann–Hilbert problem for Jacobi polynomials with singular weights}\label{subsec RHP intro}

As discussed in the Introduction, the analysis of the largest eigenvalue of the matrix $({\mathbf A} +{\mathbf B})^{-1}{\mathbf B}$ relies on the asymptotics of the associated Jacobi polynomials, which in our setting correspond to the modified weights given in \eqref{eq:eqweight}. As shown in \cite{fokas:its:kitaev:1992}, these asymptotics, in turn, can be derived from a Riemann–Hilbert problem (RHP) characterizing the polynomials. The specific RHP for the Jacobi polynomials needed in our work is stated next.

\begin{RHP}\label{RHP Y}
Find a $2 \times 2$ matrix valued function $Y = Y(z) = Y(z;p)$ that satisfies the usual three conditions of analyticity in the upper and lower half-planes, jump relations across the real axis, and prescribed behavior as $z \to \infty$.
\begin{enumerate}
\item $Y(z) $ is analytic in $\mathbb{C} \setminus [0,1]$.
\item $Y(z)$ has smooth (actually analytic) boundary values $Y_{\pm}(x)$ for $ x \in (0,1)$, and satisfies the jump relation
\begin{eqnarray}
Y_{+}(x) = Y_{-}(x) \pmtwo{1}{w_p(x)}{0}{1} \ , \ \ x \in (0,1) \ .
\end{eqnarray}
\item $Y(z)$ has the following asymptotic behavior:
\begin{eqnarray}
Y(z) = \left( I_2 + \mathcal{O} \left( z^{-1} \right) \right) z^{p \sigma_{3}} =  \left( I_2 + \mathcal{O} \left( z^{-1} \right) \right) \pmtwo{z^{p}}{0}{0}{z^{-p}} \ , \ \ z \to \infty \ . 
\end{eqnarray}
\item $Y(z)$ is bounded at $0$ and $1$.
\end{enumerate}
\end{RHP}

The following classical result establishes the existence of a solution for RHP \eqref{RHP Y}. For its proof, see, for example, \cite{Deift2000}.

\begin{proposition}\label{RHP Y sol polynomials}
The unique solution of RHP \eqref{RHP Y} is given by
\begin{center}
  Y(z)\,=$\begin{pmatrix}
       \pi_p(z) & \frac{1}{2\pi i} \int_{0}^1\frac{\pi_p(t)w_p(t)}{t-z}\, dt\\
       -2\pi i \, \kappa_{p-1}^2 \pi_{p-1}(z) & \frac{-2\pi i \, \kappa_{p-1}^2}{2\pi i} \int_{0}^1\frac{\pi_{p-1}(t)w_p(t)}{t-z}\, dt\\
    \end{pmatrix},$
\end{center}
where $\pi_p(z)$ denotes the monic Jacobi polynomial of degree $p$, with respect to the weight $w_p(x)$ in \eqref{eq:eqweight} and $\kappa_{p-1}^2=\left(\int_0^1 \pi_{p-1}^2(s)w_p(s)\,ds\right)^{-1}$. 
\end{proposition}

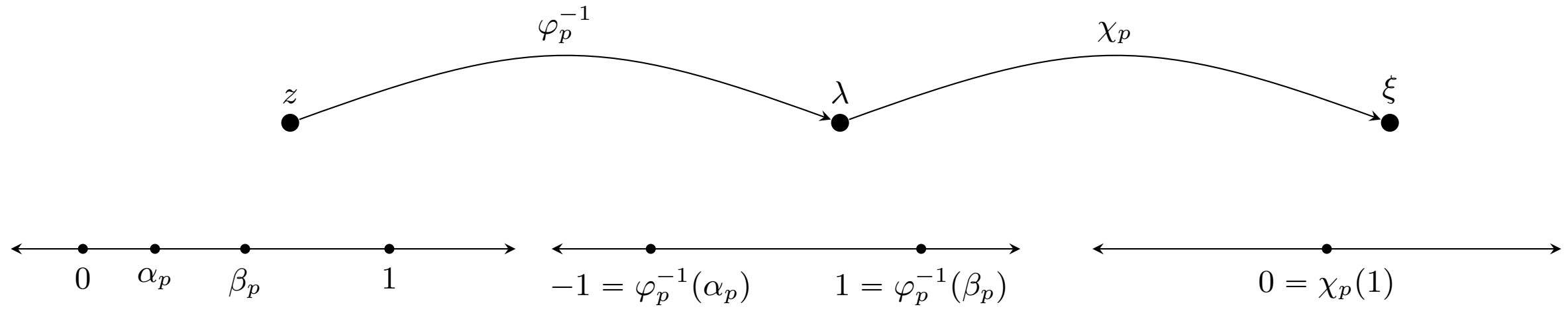
\begin{figure}[htbp]
\centering
\begin{tikzpicture}[scale=0.9, >=stealth, every node/.style={font=\small}]

\def\axislevel{0}

\begin{scope}[xshift=-6cm]
    \draw[<->] (-2.6, \axislevel) -- (3, \axislevel);
    \node[above=4pt] at (0, \axislevel) {};
    
    \node[circle,fill,inner sep=1pt] at (-1.8,\axislevel) {};
    \node[below=2pt] at (-1.8,\axislevel) {$0$};
    
    \node[circle,fill,inner sep=1pt] at (-1,\axislevel) {};
    \node[below=2pt] at (-1,\axislevel) {$\alpha_p$};

        \node[circle,fill,inner sep=1pt] at (0,\axislevel) {};
    \node[below=2pt] at (0,\axislevel) {$\beta_p$};
    
    \node[circle,fill,inner sep=1pt] at (1.6,\axislevel) {};
    \node[below=2pt] at (1.6,\axislevel) {$1$};

    \node[circle,fill,inner sep=1.8pt] (zpoint) at (0.5, 1.4) {};
    \node[above=2pt] at (zpoint) {$z$};
\end{scope}

\begin{scope}[xshift=0cm]
    \draw[<->] (-2.6, \axislevel) -- (2.6, \axislevel);
    \node[above=4pt] at (0, \axislevel) {};
    
    \node[circle,fill,inner sep=1pt] at (-1.5,\axislevel) {};
    \node[below=2pt] at (-1.5,\axislevel) {$-1=\varphi_p^{-1}(\alpha_p)$};
    
    \node[circle,fill,inner sep=1pt] at (1.5,\axislevel) {};
    \node[below=2pt] at (1.5,\axislevel) {$1=\varphi_p^{-1}(\beta_p)$};
    
    \node[circle,fill,inner sep=1.8pt] (lampoint) at (0.6, 1.4) {};
    \node[above=2pt] at (lampoint) {$\lambda$};
\end{scope}

\begin{scope}[xshift=6cm]
    \draw[<->] (-2.6, \axislevel) -- (2.6, \axislevel);
    \node[above=4pt] at (0, \axislevel) {};
    
    \node[circle,fill,inner sep=1pt] at (0,\axislevel) {};
    \node[below=2pt] at (0,\axislevel) {$0=\chi_p(1)$};
    
    \node[circle,fill,inner sep=1.8pt] (xipoint) at (0.7, 1.4) {};
    \node[above=2pt] at (xipoint) {$\xi$};
\end{scope}

\draw[->, ]
    (zpoint) to[out=20, in=160, looseness=1.2]
    node[above, pos=0.5, font=\small] {$\varphi_p^{-1}$}
    (lampoint);

\draw[->, ]
    (lampoint) to[out=20, in=160, looseness=1.2]
    node[above, pos=0.5, font=\small] {$\chi_p$}
    (xipoint);

\end{tikzpicture}
\caption{Conformal mappings $z \mapsto \lambda \mapsto \xi$.} 
\label{fig:combined-maps}
\end{figure}

Following classical Riemann--Hilbert analysis, the solution of RHP~\ref{RHP Y} involves a sequence of transformations $z \mapsto \lambda \mapsto \xi$. These changes of variables, which are particularly relevant to the rescaling of the reproducing kernel, are briefly introduced below and illustrated in Figure~\ref{fig:combined-maps}. Their role in the Riemann–Hilbert analysis is discussed in detail in Section~\ref{section RH Analysis}.

In the terminology of Riemann-Hilbert analysis, RHP~\ref{RHP Y} is said not to be in standard form since $Y(z)$ does not approach the identity matrix as $z \to \infty$. For this reason, a transformation is introduced to convert the problem into an equivalent RHP normalized at infinity. This transformation embeds a  so-called \textit{$g$-function} in the problem associated 
to an equilibrium measure via the following relation:
\begin{eqnarray}
g(z) = \int \log{(z - s)} \psi(s) \, ds, \quad z \in \mathbb{C} \setminus (-\infty, \beta_p)\end{eqnarray}
where $\beta_p$ will be determined from the RHP analysis. Moreover, the $g$-function determines an interval $(\alpha_p, \beta_p) \subset (0,1)$ where the distribution of the bulk eigenvalues is most concentrated as $p\rightarrow \infty$. Further details on the construction of the $g$-function and the determination of these endpoints are provided in Subsection~\ref{g function}. In particular, $\alpha_p$ and $\beta_p$ are given by:
\begin{equation}\label{eq. alpha_beta_explicit}
\begin{split}
  \alpha_p & =\frac{-(a^2_{\infty} p^{2\nu}-b^2_{\infty}-\eta) - \sqrt{(a^2_{\infty} p^{2\nu}-b^2_{\infty}-\eta)^2 - 4\eta\, b^2_{\infty}}}{2\eta},\\
  \beta_p& =\frac{-(a^2_{\infty} p^{2\nu}-b^2_{\infty}-\eta) + \sqrt{(a^2_{\infty} p^{2\nu}-b^2_{\infty}-\eta)^2 - 4\eta\, b^2_{\infty}}}{2\eta},  
\end{split}
\end{equation}
with $\eta = (b_{\infty} + a_{\infty} p^{\nu} + 2)^2$. 

As we will see in Remark~\ref{remark asymp alpha beta}, in the unbalanced regime $\nu>0$, both endpoints $\alpha_p$ and $\beta_p$ tend to zero as $p\to\infty$, and consequently the support of the $g$-function shrinks. So, to capture the full asymptotic analysis of the Jacobi polynomials, we introduce the linear transformation
\begin{equation}\label{eq:eqchangeofvari zlambda}
      z = \varphi_p(\lambda) \coloneqq  \frac{\beta_p + \alpha_p}{2} + \frac{\beta_p - \alpha_p}{2} \, \lambda.
\end{equation}
The map $\varphi_p^{-1}: z \mapsto \lambda$ scales the original $z$-plane as $\lambda \approx p^{\nu}z$. This transforms RHP~\ref{RHP Y} from the $z$-plane to the $\lambda$-plane by mapping $\alpha_p$ and $\beta_p$ to the fixed points $\lambda=-1$ and $\lambda=1$, respectively.

The next step in the analysis is to perform the transformation known as "opening lenses". This deforms the RHP, so that rapidly oscillatory terms appearing in the jump matrix are converted into exponentially decaying ones. Effectively, this allows us to disregard jumps that are exponentially close to the identity matrix, reducing the RHP to one with a jump condition only along the interval $(-1,1)$. Solving this deformed problem provides an approximate solution, referred to as the \textit{outer parametrix}, which is valid away from the endpoints $\lambda = \pm 1$. In the vicinity of these points, however, the outer parametrix ceases to be valid, since the previously neglected jump matrices are no longer close to the identity matrix. Consequently, we must construct a local solution near $\lambda = \pm 1$. 
In the region around $\lambda=1$ (which is related to the largest eigenvalue),
we introduce the local conformal map $\chi_p : \lambda \mapsto \xi$ defined by
\begin{equation}\label{eq:mapchi}
\xi:= \chi_p(\lambda)  =p^{\frac{2}{3}} \left(- \frac{3}{2} \, \,\phi(\lambda)\right)^{2/3}.
\end{equation}
In \eqref{eq:mapchi}, $\phi(\lambda)$ is a holomorphic function defined in a neighborhood of $\lambda = 1$, and is given by the integral representation \eqref{eq:funphi}. This map has local behavior $\xi \approx p^{\frac{2}{3}}(\lambda-1)$ (see~\eqref{eq. chi_p asymtotics}). Hence, it maps a sufficiently small neighborhood of $\lambda=1$ onto a comparatively large neighborhood of $\xi=0$.

Overall, the correspondence between the new local variable \(\xi\) and the original variable \(z\) is given by the composite inverse function
\begin{equation}\label{eq. map xi to z^p}
z:=\varphi_p\big(\chi_p^{-1}(\xi)\big).
\end{equation}
For a comprehensive analysis of RHP~\ref{RHP Y} and its transformations, we refer the reader to Section~\ref{section RH Analysis}.

\subsection[Connecting the RHP to the distribution of
the largest eigenvalue]{
Connecting the RHP~\ref{RHP Y} to the distribution
of $\theta_{\max}$: Fredholm determinants and
rescaled kernels}\label{s:connecting_RH_top_eigenval}

In this section, we provide a succinct description of the connection between, on one hand, RHP \ref{RHP Y} and the change of variables described in Section \ref{subsec RHP intro}, and, on the other hand, the analysis of the distribution of $\theta_{\max}=\theta_{\max,p}$.

Recall from \(\eqref{eq: prob lambda max}\) that the distribution of $\theta_{\max}$ admits the Fredholm determinant representation
\begin{equation}\label{e:P(theta_max<a)_recap}
\mathbb{P}(\theta_{\max}<a)=\det( \mathbb{I}-\mathcal{K}_p)_{L^2([a,1])}, \quad 0 < a < 1,
\end{equation}
where $\mathcal{K}_p$ is the integral operator \eqref{eq:eqdefnoper} with reproducing kernel $K_p$ defined in \eqref{eq:def Kernel z}. As it turns out, we can express \(K_p\) in terms of the entries of the matrix $Y$ that solves RHP \ref{RHP Y}. In fact, starting from relation \eqref{eq:def Kernel z}, an application of the Christoffel–Darboux formula and of Proposition \ref{RHP Y sol polynomials} yields
\begin{equation}\label{eq:kernelp}
K_p(z_1,z_2)  = -\frac{\sqrt{w_p(z_1)w_p(z_2)}}{2\pi i}\left(\frac{Y_{11}(z_1)Y_{21}(z_2)-Y_{11}(z_2)Y_{21}(z_1)}{z_1-z_2} \right).
\end{equation}
We now examine how each step of the change of variables \eqref{eq. map xi to z^p} successively transforms this kernel and its associated integral operator. Based on the first step -- namely, the transformation \(z=\varphi_p(\lambda)\), $\lambda \in [\varphi^{-1}_p(0),\varphi^{-1}_p(1)]$, as in \eqref{eq:eqchangeofvari zlambda} --, expression \eqref{eq:kernelp} becomes
\begin{equation}\label{eq:eq32}
K_p\big(\varphi_p(\lambda_1),\varphi_p(\lambda_2)\big) = -\frac{\sqrt{w_p(\varphi_p(\lambda_1))w_p(\varphi_p(\lambda_2))}}{\pi i(\beta_p-\alpha_p)}\left(\frac{Y_{11}(\varphi_p(\lambda_1))Y_{21}(\varphi_p(\lambda_2))-Y_{11}(\varphi_p(\lambda_2))Y_{21}(\varphi_p(\lambda_1))}{\lambda_1-\lambda_2} \right).
\end{equation}
We can now define the corresponding rescaled kernel in the \(\lambda\)-plane as
\begin{equation}\label{eq:def Kernel lambda}
  \widetilde{K}_p (\lambda_1, \lambda_2): = \frac{\beta_p-\alpha_p}{2} K_p \big( \varphi_p(\lambda_1), \varphi_p(\lambda_2) \big).
\end{equation}
Now let
\[ \widetilde{Y}(\lambda)=\left(\frac{\beta_p-\alpha_p}{2}\right)^{-p\sigma_3}Y(\varphi_p(\lambda))\]
be the solution of RHP \ref{RHP Y} in the $\lambda$-plane. Then, we can express the transformed kernel \eqref{eq:def Kernel lambda} directly in terms of \(\widetilde Y\), i.e.,
\begin{equation}\label{eq:eqlastrelation}
\widetilde{K}_p (\lambda_1, \lambda_2) =-\frac{\sqrt{w_p(\varphi_p(\lambda_1))w_p(\varphi_p(\lambda_2))}}{2 \pi i}\left(\frac{\widetilde{Y}_{11}(\lambda_1)\widetilde{Y}_{21}(\lambda_2)-\widetilde{Y}_{11}(\lambda_2)\widetilde{Y}_{21}(\lambda_1)}{\lambda_1-\lambda_2} \right).
\end{equation}
The associated integral operator \(\widetilde{\mathcal K}_p\), acting on $L^2([\varphi_p^{-1}(a),\varphi_p^{-1}(1)])$, is naturally given by
\begin{equation}\label{eq:eqoperator lambda}
\widetilde{\mathcal{K}}_p f(\lambda_1) = \int_{\varphi_p^{-1}(a)}^{\varphi_p^{-1}(1)} \widetilde{K}_p (\lambda_1, \lambda_2) f(\lambda_2) \, d \lambda_2.
\end{equation}
Then, in terms of $\widetilde{\mathcal{K}}_p$ as in \eqref{eq:eqoperator lambda}, relation \eqref{e:P(theta_max<a)_recap} can be recast in the form
\begin{equation}\label{eq: prob Ktilde}
  \mathbb{P}(\theta_{\max}<a)=\det( \mathbb{I}-\widetilde{\mathcal{K}}_p)_{L^2([\varphi_p^{-1}(a),\varphi_p^{-1}(1)])}, \quad 0 < a < 1. 
\end{equation}
The present work uses~\eqref{eq: prob Ktilde} to study the distribution of $\theta_{\max}$, which, as we will show, is concentrated near the value $\beta_p$ defined in~\eqref{eq. alpha_beta_explicit}. However, as $p \rightarrow \infty$, the behavior of the Fredholm determinant in~\eqref{eq: prob Ktilde} depends on the position of $a$ relative $\beta_p$. In this work, we fix a constant $L>0$ and consider a sufficiently small disk \(\mathcal U_1  \subseteq \mathbb{C}\) centered at \(\lambda=\varphi^{-1}(\beta_p)=1\), where the local map $\chi_p$ as in~\eqref{eq:mapchi} is defined (cf.\ Remark~\ref{remark epsilon}).  We denote by \(\lambda_D\) the rightmost intersection of its boundary with the real axis (see Figure \(\ref{fig Contour E}\)). We consider the following two cases:
\begin{itemize}
\item [$(i)$] the value $a$ is near $\beta_{p}$, as $p\rightarrow\infty$. More precisely, we consider $a \in \left(\varphi_p(\chi_p^{-1}(-L))\,,\varphi_p(\lambda_D)\right)$. In this case we have 
\[
\varphi_p(\chi_p^{-1}(-L))
=\beta_p-C_1 \, p^{-\nu-2/3} + o\left(p^{-\nu-2/3}\right) \qquad\text{and} \qquad \varphi_p(\lambda_D)=\beta_p+ C_2\, p^{-\nu}+o \left( p^{-\nu} \right)
\]
where $C_1$ and $C_2$ are positive constants. In particular this implies that $\varphi_p^{-1}(a) \in \left(\chi_p^{-1}(-L),\lambda_D\right)\subset \mathcal{U}_1$. See Section~\ref{eq near edge} for the analysis.

\item [$(ii)$]  the value $a$ lies to the right of $\beta_{p}$, at a distance of at least order $p^{-\nu}$. More precisely, there is a constant $C>0$, such that for all sufficiently large $p$ we have 
\[
a\ge \beta_p+C\,p^{-\nu}.
\]
In this case, $\varphi_p^{-1}(a) \in \left[\lambda_D, \varphi_p^{-1}(1)\right)\subset \, \mathbb{R}\,\setminus\,\mathcal{U}_1$. See Section \ref{subsec eigval stats case 1} for the analysis.
\end{itemize}

\begin{remark}
Note that cases $(i)$ and $(ii)$ exhaust the case $a\ge\beta_p$: in the first case, $a$ lies to the right of $\beta_p$ at a distance of order $p^{-\nu}$ or smaller, whereas in the second case, it lies at a distance of order $p^{-\nu}$ or larger. On the other hand, item $(i)$ is related to the case where $a$ lies to the left of $\beta_p$ at a distance of order $p^{-\nu-2/3}$ or smaller. The case where $a$ lies to the left of $\beta_p$, at a distance of order bigger than $p^{-\nu-2/3}$ (lower-tail regime) is not considered in this work (see also Remark~\ref{r:lower_tail}).
\end{remark}

It turns out that in case $(ii)$, the Fredholm determinant in \eqref{eq: prob Ktilde} is exponentially small in $p$ (see Lemma~\ref{Lemma: right of lambdaD}). On the other hand, in case $(i)$ the Fredholm determinant can be restricted to the local interval $(\varphi_p^{-1}(a),\lambda_D)$, rather than considering the entire interval $[\varphi_p^{-1}(a),\varphi_p^{-1}(1)]$. More precisely, Lemma~\ref{Prop truncating kernel} shows that there exists a constant $k>0$ such that
\begin{equation}\label{eq truncating K tilde}
    \det\left( \mathbb{I} - \widetilde{\mathcal{K}}_{p} \right)_{L^2([\varphi_p^{-1}(a),\varphi_p^{-1}(1)])} = \det\left( \mathbb{I} - \widetilde{\mathcal{K}}_{p} \right)_{L^2([\varphi_p^{-1}(a),\lambda_D])}+\mathcal{O}\left(e^{-k p}\right),
\end{equation}
where the implicit constant in $ \mathcal{O}\left(e^{-k p}\right)$ depends only on the fixed value $\lambda_D$, or equivalently, on the radius of the disk $\mathcal{U}_1$.

Based on the second (i.e., the last) step of the change of variables \eqref{eq. map xi to z^p}, we now apply the map \(\chi_p(\lambda)\), $\lambda \in {\mathcal U}_1$, as in \eqref{eq:mapchi} to express the kernel \(\widetilde K_p\) as in \eqref{eq:eqlastrelation} in the new local coordinate \(\xi\). After proper renormalization, the resulting kernel can be defined as
\begin{equation}\label{kernel local coordinates}
    \widehat{K}_p(\xi_1,\xi_2) \coloneqq \frac{1}{\sqrt{\chi_p'(\chi_p^{-1}(\xi_1))\,\chi_p'(\chi_p^{-1}(\xi_2))}}\,\widetilde{K}_p\left(\chi_p^{-1}(\xi_1),\chi_p^{-1}(\xi_2)\right), \quad \xi_1, \xi_2\in \chi_p(\mathcal{U}_1).
\end{equation}
The corresponding integral operator \(\widehat{\mathcal K}_p\), acting on $L^2([\chi_p(\varphi_p^{-1}(a)),\chi_p(\lambda_D)])$, is naturally defined by
\begin{equation}\label{def operator K hat}
\widehat{\mathcal{K}}_p f(\xi_1) = \int_{\chi_p \left( \varphi_p^{-1}(a)\right)}^{\chi_p \left( \lambda_D \right)} \widehat{K}_p (\xi_1, \xi_2) f(\xi_2) \, d \xi_2.
\end{equation}
Based on the operator $\widehat{\mathcal{K}}_p$ as in \eqref{def operator K hat}, the Fredholm determinant on the right-hand side of \eqref{eq truncating K tilde} can be recast in the form
\begin{equation}\label{e:Fred_det_K-tilde_to_Fred_det_K-hat}\det\left( \mathbb{I} - \widetilde{\mathcal{K}}_{p} \right)_{L^2([\varphi_p^{-1}(a),\lambda_D])} =\det( \mathbb{I}-\widehat{\mathcal{K}}_p)_{L^2([\chi_p(\varphi_p^{-1}(a)),\chi_p(\lambda_D)])}.
\end{equation}
As a consequence of relations \eqref{eq: prob Ktilde}, \eqref{eq truncating K tilde} and \eqref{e:Fred_det_K-tilde_to_Fred_det_K-hat}, we can now conveniently write
\begin{equation}\label{eq. prob Kp hat}
  \mathbb{P}(\theta_{\max}<a) \;=\; \det\!\left(1-\widehat{\mathcal{K}}_p\right)_{L^2([\chi_p(\varphi_p^{-1}(a)),\,\chi_p(\lambda_D)])} \;+\; \mathcal{O}\left(e^{-kp}\right).
\end{equation}
In other words, relation \eqref{eq. prob Kp hat} establishes the connection between the study of the asymptotic behavior of $\theta_{\max}$ to that of the operator $\widehat{\mathcal{K}}_p$ as in \eqref{def operator K hat}, and hence, to the RHP \ref{RHP Y}.

For the reader's convenience, we further summarize the successive transformations of the kernels and their associated integral operators in Table~\ref{tab:notation}.

\begin{table}[ht]
\centering
\renewcommand{\arraystretch}{2.5} 
\setlength{\tabcolsep}{4pt}       
\begin{tabular}{|p{1.5cm}|p{6.5cm}|c|}
\hline
\textbf{Operator} & \multicolumn{1}{c|}{\textbf{Action}} & \textbf{Kernel} \\
\hline
$\mathcal{K}_p$ & $\mathcal{K}_p f(z_1)=\int_a^{1} K_p(z_1,z_2) f(z_2) dz_2$ \eqref{eq:eqdefnoper} & $K_p(z_1,z_2) = \sqrt{w_p(z_1)w_p(z_2)}\,
\sum\limits_{\ell=0}^{p-1}
\mathbf{p}_\ell(z_1)\mathbf{p}_\ell(z_2)$ \eqref{eq:def Kernel z}\\
\hline
$\widetilde{\mathcal{K}}_p$ & $  \widetilde{\mathcal{K}}_p f(\lambda_1) = \displaystyle \int_{\varphi_p^{-1}(a)}^{\varphi_p^{-1}(1)} \widetilde{K}_p (\lambda_1, \lambda_2) f(\lambda_2) d \lambda_2$  \eqref{eq:eqoperator lambda} & $\widetilde{K}_p (\lambda_1, \lambda_2) = \frac{\beta_p-\alpha_p}{2} K_p \left( \varphi_p(\lambda_1), \varphi_p(\lambda_2) \right)$ \eqref{eq:def Kernel lambda}\\
\hline
$\widehat{\mathcal{K}}_p$ & $\widehat{\mathcal{K}}_p f(\xi_1) = \displaystyle \int_{\chi_p \left( \varphi_p^{-1}(a)\right)}^{\chi_p \left( \varphi_p^{-1}(1)\right)} \widehat{K}_p (\xi_1, \xi_2) f(\xi_2) \, d \xi_2$ \eqref{def operator K hat} & $ \widehat{K}_p(\xi_1,\xi_2) = \frac{\widetilde{K}_p\left(\chi_p^{-1}(\xi_1),\chi_p^{-1}(\xi_2)\right)}{\sqrt{\chi_p'(\chi_p^{-1}(\xi_1))\,\chi_p'(\chi_p^{-1}(\xi_2))}}$ \eqref{kernel local coordinates}\\
\hline
\end{tabular}
\caption{Summary of operators, kernels, and their transformations under the changes of variables.}
\label{tab:notation}
\end{table}


\section{Main results}\label{subsec statement results}

In this section, we present the three main results of the paper. All claims pertain to the random matrix ${\mathbf J}$ as in \eqref{e:(A+B)^(-1)B} under the choice $\beta = 2$ in the eigenvalue p.d.f.\ \eqref{e:joint_eigenvalue_pdf_beta}. Throughout this section, we assume that the conditions in~\eqref{eq: parameters us} hold. Recall also that $\theta_{\max}$ denotes the largest eigenvalue of  ${\mathbf J}$ (cf.\ \eqref{eq: prob lambda max}). The proofs are provided in Section~\ref{proofs} and rely on the asymptotic analysis of Jacobi polynomials developed in Section~\ref{section RH Analysis} and the analysis of the Fredholm determinants characterizing the distribution of $\theta_{\max}$ presented in Section~\ref{proofs auxiliary}.

The first theorem establishes that, in the asymptotic regime~\eqref{eq: parameters us}, after the rescaling determined by the maps $\chi_p$ and $\varphi_p$, the cumulative distribution function (c.d.f.) of the largest eigenvalue admits an asymptotic expansion whose leading term is the type 2 Tracy--Widom distribution. In particular, the theorem provides an explicit correction of order $p^{-2/3}$ in terms of the Painlev\'{e}~II function.

\begin{theorem}\label{thm: Probability Tracy-Widom}
Fix $L>0$ . Then, for all $L_0>-L$, 
\begin{align}\label{eq:mres2}  
    \mathbb{P}\left(\theta_{\max}\le \varphi_p\big(\chi_p^{-1}(L_0)\big)\right) 
      = F_2(L_0) \left( 1 -2\, \widehat{C}_p \, p^{-2/3} \, v(L_0) \right) + \mathcal{O} \left( p^{-1} \right), \qquad p \to \infty,
\end{align}
with the implicit constant in $\mathcal{O} \left( p^{-1} \right)$ independent of $L_0$. Here, $\widehat{C}_p$ denotes the coefficient appearing in \eqref{eq:errt}, and the function $v(s)$ is given explicitly by
\begin{equation} \label{eq: def v}
v(s) = \frac{1}{2} \left( \int_s^{\infty} q^2(t) \, dt \right)^2 - \frac{q^2(s)}{2},
\end{equation}
where $q(s)$ is the same solution of the Painlev\'{e} II equation described in~\eqref{eq. PII}
\end{theorem}
\begin{remark}
Using the definition of the Tracy-Widom law in \eqref{e:F2=exp(integ)} and the definition of $v$ in \eqref{eq: def v}, relation \eqref{eq:mres2} can be recast  more explicitly in the form:
\begin{small}
\begin{equation}\label{TWcorrection}  
\mathbb{P}\left(\theta_{\max}\le \varphi_p\big(\chi_p^{-1}(L_0)\big)\right)  
       = \mathrm{exp} \left[ - \int_{L_0}^{\infty} (y-L_0) \, q^2(y) \, dy + \widehat{C}_p \, p^{-2/3} \left( \frac{1}{2} \left( \int_{L_0}^{\infty} q^2(t) \, dt \right)^2 - \frac{q^2(L_0)}{2} \right) \right] + \mathcal{O} \left( p^{-1} \right).
\end{equation}
\end{small}
\end{remark}

\begin{figure}[htbp]
    \centering
    \includegraphics[width=0.5\textwidth]{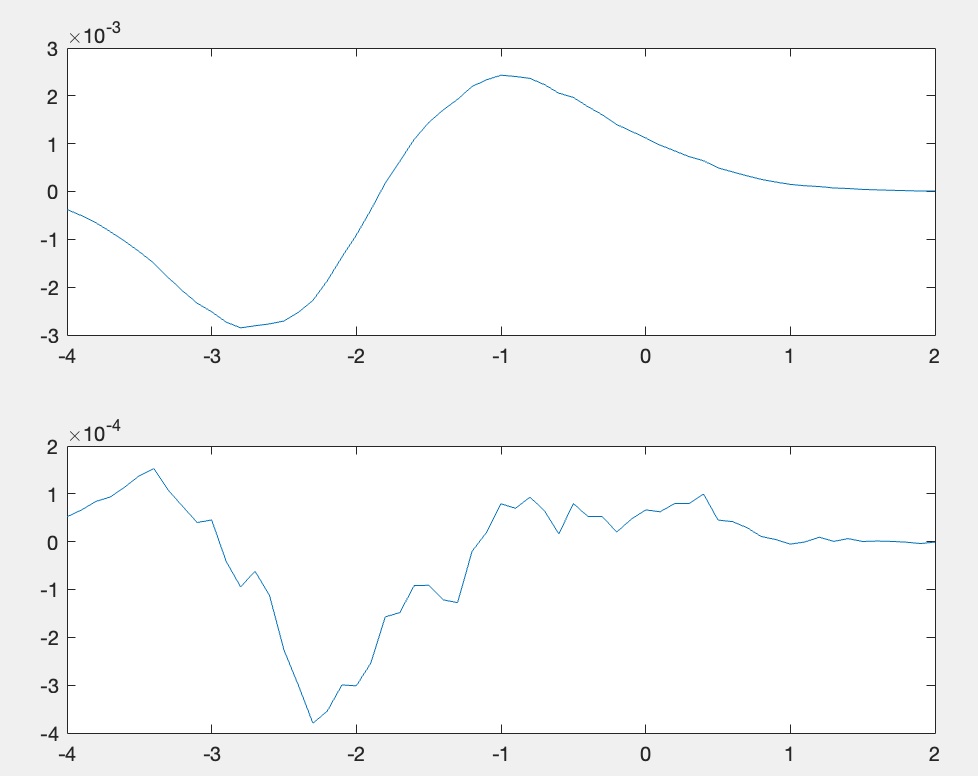}
    \caption{$R = 10{,}000{,}000$ random matrix runs were used to estimate the empirical c.d.f.\ of the largest eigenvalue of the complex Jacobi ensemble under the parameter choices $\nu=0$, $a_{\infty}=b_{\infty}=1$ and $p=100$. The plots display the deviation of the Tracy--Widom law (top) and of the corrected Tracy--Widom approximation (bottom) from the estimated empirical c.d.f. The correction improves the approximation by one order of magnitude.}
    \label{fig:twcomparison}
\end{figure}

To illustrate the claim of Theorem \ref{thm: Probability Tracy-Widom}, Figure \ref{fig:twcomparison} displays a computational study of the classical (balanced) instance $\nu = 0$. The two plots show the deviations of the Tracy--Widom law and of the corrected Tracy--Widom approximation provided in the theorem from the simulated empirical c.d.f. 

The second result establishes the convergence in distribution of $\theta_{\max}$ to the type 2 Tracy--Widom law, and explicitly describes the centering and scaling terms involved in the convergence. The result follows from the analysis of the leading terms in the expansion provided in Theorem~\ref{thm: Probability Tracy-Widom}.

\begin{corollary}\label{mainthrm}
Let $\beta_p$ be defined as in \eqref{eq. alpha_beta_explicit}. Then, as $p\to\infty$, for $\nu\ge 0$:
\begin{equation}\label{eq:main-convergence}
C\,p^{\nu+\frac{2}{3}}
\bigl(\theta_{\max}-\beta_p\bigr)
\xrightarrow{\mathrm{d}}
\mathrm{TW}_2.
\end{equation}
Here:
\begin{equation}
C= \begin{cases}
 \frac{1}{2}\,
\frac{a_{\infty}}
{\sqrt{1+b_{\infty}}}\left(\frac{2\sqrt{2} \, (1+b_{\infty})}{2 + b_{\infty} + 2\sqrt{1+b_{\infty}}}\right)^{2/3}, &\quad \nu>0,\\
\frac{(2+a_\infty+b_{\infty})^2}{2\sqrt{(1+a_{\infty})(1+b_{\infty})(1+a_{\infty}+b_{\infty})}}(-h(1))^{2/3}, & \quad \nu=0,
\end{cases}
\end{equation}
where $h(1)$ is the $p$-independent constant given explicitly in \eqref{eq:h1}. 
\end{corollary}

\begin{remark}\label{remark main thm inverse}
Using an asymptotic expansion of the map $\chi_p\left(\varphi_p^{-1} (\cdot) \right)$ and the fact that $\Omega'(1)=1$, where $\Omega(\lambda)$ is defined in \eqref{eq. def Omega}, we obtain:
\begin{align}\label{eq M0}
M_0:=\chi_p\left(\varphi_p^{-1}\left(\beta_p+\frac{L_0}{C\,p^{\nu+2/3}}\right)\right)=L_0+ \widetilde{c}_2 L_0^2 p^{-2/3}+\mathcal{O}\left(p^{-4/3}\right),
\end{align}
\begin{align*}
\widetilde{c}_2 = \frac{\Omega_0''(1)}{2} \left( \frac{2 \sqrt{2} (1+b_{\infty})}{2+b_{\infty}+2 \sqrt{1+b_{\infty}}}\right)^{-2/3}, \quad \Omega_0(\lambda)=\sum_{j \geq 0} \frac{3 \left[h_p^{(j)}(1)\right]^{(1)}}{(2j+3) j! \left[h_p(1)\right]^{(1)}}(\lambda-1)^j.
\end{align*}
Here the superscript $1$ denotes the leading order term in the large-$p$ asymptotic expansion of $h_p'(1)$ and $h_p(1)$. By comparison to \eqref{eq:main-convergence}, using equation \eqref{eq M0},  we can rewrite relation~\eqref{eq:mres2} in Theorem~\ref{thm: Probability Tracy-Widom} as follows:
\begin{align*}
   \mathbb{P}\Big(C p^{\nu + 2/3} \big(\theta_{\max} - \beta_p \big) \le L_0\Big)&= F_2(M_0)\left( 1 -2\, \widehat{C}_p \, p^{-2/3} \, v(M_0) \right) + \mathcal{O} \left( p^{-1} \right)\\
      &=F_2(M_0)\left\{1+p^{-2/3} \left[ -2 \widehat{C}_p\, v(M_0) + \widetilde{c}_2 \, M_0^2 \left( 2 v(M_0) +q^2(M_0)\right)^{1/2}\right]\right\}+\mathcal{O}\left(p^{-1}\right).
\end{align*}
\end{remark}

The third result shows that the probability of finding eigenvalues "far" to the right from $\beta_p$ (see equation~\eqref{eq:thrm34}) decays exponentially as $p$ goes to infinity.

\begin{theorem}\label{prop:prop7}
Let $\nu\ge0$. For every fixed $c>0$, there exist constants $K,k>0$ such that, for all sufficiently large $p$ and all $\theta\in (0,1)$ satisfying
\begin{equation}\label{eq:thrm34}
\theta-\beta_p \geq c\,p^{-\nu},
\end{equation}
we have
\begin{equation*}
\mathbb{P}\!\left(\theta_{\max}>\theta\right)
\leq K e^{-kp}.
\end{equation*}
\end{theorem}

\begin{remark}\label{r:lower_tail}
In this paper, we do not establish the analogous lower-tail claim to Theorem 
\ref{prop:prop7}, namely, an estimate of the probability of the event $\left\{ \theta_{\max} < \beta_{p} - K p^{- \nu} \right\}$ for some constant $K > 0$. Although the representation of this probability as a Fredholm determinant still holds, the asymptotic analysis of the kernel involves several different asymptotic regimes. This, in turn, precludes the convergence of the Fredholm determinant itself.  Typically, such lower-tail estimates are established using large deviation-type estimates via potential theory (e.g., Dean and Majumdar \cite{dean:majumdar:2006}), which is beyond the scope of this work.
\end{remark}

\begin{remark}In this paper, we only consider the largest eigenvalue $\theta_{\max}$ of the complex Jacobi ensemble.  Nevertheless, the methods present can be adjusted to study the asymptotic behavior of the smallest eigenvalue $\theta_{\min}$ as well. It is expected that this random variable, after suitable rescaling, will also converge in distribution to the Tracy-Widom law.\end{remark}

\begin{remark}
As stated in the theorems, our analysis remains valid in the case $\nu=0$. However, certain arguments require some adjustments. We indicate these at the relevant points throughout the paper.
\end{remark}


\section{Proofs of the main results}\label{proofs}

In this section, we present the proofs of our main results, i.e., Theorems~\ref{thm: Probability Tracy-Widom},~\ref{mainthrm} and~\ref{prop:prop7}. For the reader's convenience, we also provide the statements of some key auxiliary results -- their proofs are deferred to Section~\ref{proofs auxiliary}. The proofs in this section further rely on the Riemann--Hilbert analysis developed in Section~\ref{section RH Analysis}.

Throughout this section $a\in(0,1)$ denotes a fixed value in the original domain of $\theta_{\max}$. So, for the complex Jacobi ensemble, recall the Fredholm determinant representation~\eqref{eq: prob lambda max} of the distribution of the largest eigenvalue:
\begin{equation}\label{e:P(theta_max<a)_recap2}
\mathbb{P}(\theta_{\max}<a)=\det( \mathbb{I}-\mathcal{K}_p)_{L^2([a,1])}, \quad 0 < a < 1
\end{equation}
In our analysis of \eqref{e:P(theta_max<a)_recap2}, we distinguish between the two cases discussed in Subsection~\ref{s:connecting_RH_top_eigenval} related to the position of $a$ relative to the value $\beta_p$.

\subsection[Largest eigenvalue behavior near the end pint]{Largest eigenvalue behavior near the endpoint $\beta_p$}\label{eq near edge}

We start with the case where $a \in \left(\varphi_p(\chi_p^{-1}(-L))\,,\varphi_p(\lambda_D)\right)$. Here, our results rely on the local convergence of the kernel $K_p$ to the so-called Airy kernel (see Propositions~\ref{proposition Khat to K_airy}), which is defined by
\begin{equation}\label{Eq. Def Airy kernel}
K_{\mathrm{Airy}}(\xi,\eta)
=
\frac{\mathrm{Ai}(\xi)\mathrm{Ai}'(\eta)
-\mathrm{Ai}(\eta)\mathrm{Ai}'(\xi)}
{\xi-\eta},
\quad \xi,\eta\in\mathbb{R}.
\end{equation}
Here, $\mathrm{Ai}(\,\cdot\,)$ denotes the \textit{Airy function}, defined by
\begin{equation}\label{Eq. Def Airy function}
\mathrm{Ai}(\xi)
= \frac{1}{\pi}\int_0^\infty
\cos\left(\frac{t^3}{3}+\xi t\right)\,dt,
\qquad \xi\in\mathbb{R}.
\end{equation}
For $s\in\mathbb{R}$, $\mathcal K_{\Airy}$ denotes the associated integral operator on $L^2([s,\infty))$, defined by the relation
\begin{equation}\label{eq: def Airy operator}
\mathcal K_{\Airy}f(x)
=\int_s^\infty K_{\Airy}(x,y)f(y)\,dy.
\end{equation}
The following proposition establishes the convergence of the kernel $\widehat K_p$ defined in~\eqref{kernel local coordinates} to the Airy kernel $K_{\Airy}$. In fact, more generally, it provides an explicit expression for the leading-order correction term of order $p^{-2/3}$. For its proof see Subsection~\ref{proof Prop kernel convergence}.

\begin{proposition}\label{proposition Khat to K_airy}
Fix $L>0$. There exist constants $C,k>0$ and an integer $M>0$, such that, for all $p\ge M$ and all $\xi_1,\xi_2>-L$,
\begin{equation}\label{eq:mainres}
\widehat{K}_p(\xi_1,\xi_2)=
K_{\Airy}(\xi_1,\xi_2)+
\widehat{C}_p \left(\mathrm{Ai}\xi_1)\mathrm{Ai}'(\xi_2)
+\mathrm{Ai}(\xi_2)\mathrm{Ai}'(\xi_1)\right)p^{-2/3}
+\widetilde{R}(\xi_1,\xi_2),
\end{equation}
where
\begin{equation}\label{eq. bounds R tilde prop}
\left|\widetilde{R}(\xi_1,\xi_2)\right|
\le
\begin{cases}
Cp^{-1}e^{-k(\xi_1+\xi_2)},
& (\xi_1,\xi_2)\in(0,\infty)^2,\\[1mm]
Cp^{-1}e^{-k\xi_1},
& (\xi_1,\xi_2)\in(0,\infty)\times(-L,0],\\[1mm]
Cp^{-1}e^{-k\xi_2},
& (\xi_1,\xi_2)\in(-L,0]\times(0,\infty),\\[1mm]
Cp^{-1},
& (\xi_1,\xi_2)\in(-L,0]^2.
\end{cases}
\end{equation}
The coefficient $\widehat{C}_p$ is given explicitly by
\begin{equation}\label{eq:errt}
\widehat{C}_p = \frac{1}{4c_1(p)}
\left( \frac{1}{2} +\frac{2}{5}\frac{h_p'(1)}{h_p(1)} \right), \qquad
c_1(p)=\left(-h_p(1)\right)^{2/3}.
\end{equation}
Here, $h_p$ is the function defined in \eqref{eq. def function h}, and both $c_1(p) $, $\widehat C_p$ are uniformly bounded in $p$.
\end{proposition}

We now turn to the properties of the correction terms appearing on the right-hand side of \eqref{eq:mainres}. In particular, we consider their extensions to an unbounded quadrant where both variables are bounded below by some fixed value $s\in \mathbb{R}$. These properties will allow us to pass from kernel convergence to the convergence of the corresponding Fredholm determinants.  For the proof of the following result see Subsection~\ref{proof Prop trace D0}.

\begin{proposition}\label{prop prop53}
Let $s\in\mathbb{R}$ be fixed, and consider the integral operator $\mathcal{D}_0$ on $L^2([s,\infty))$ with kernel
\begin{equation}\label{eq. kernel D0}
D_0(x,y)
= \widehat{C}_p\,p^{-2/3}
\bigl(\mathrm{Ai}(x)\mathrm{Ai}'(y)
+\mathrm{Ai}'(x)\mathrm{Ai}(y)\bigr).
\end{equation}
Then the operator $(\mathbb I-\mathcal{K}_{\Airy})^{-1}\mathcal D_0$ has kernel
\begin{equation}\label{eq. Kernel PsD}
K(x,y)
= \widehat{C}_p\,p^{-2/3}
\bigl(\mathrm{Ai}'(y)q(x,s)
+\mathrm{Ai}(y)p(x,s)\bigr),
\end{equation}
where
\begin{equation}\label{eq. def p q}
\begin{aligned}
q(x,s) =
\left((\mathbb{I}-\mathcal K_{\Airy})^{-1}
\mathrm{Ai}\right)(x), \quad p(x,s) =
\left((\mathbb{I}-\mathcal K_{\Airy})^{-1}
\mathrm{Ai}'\right)(x).
\end{aligned}
\end{equation}
Moreover, its trace is given by
\begin{equation}\label{eq. Trace PsD}
\mathrm{Tr}\bigl(
(\mathbb I-\mathcal{K}_{\Airy})^{-1}\mathcal D_0
\bigr)
=2\,\widehat{C}_p\,p^{-2/3}\,v(s),
\end{equation}
where $v(s)$ is defined by 
\begin{equation}\label{eq. def v}
v(s)=
\left\langle
\mathrm{Ai},
(\mathbb{I}-\mathcal K_{\Airy})^{-1}
\mathrm{Ai}'
\right\rangle
\end{equation}
and $\langle\cdot,\cdot\rangle$ denotes the inner product on $L^2([s,\infty))$.
\end{proposition}

Now, we fix $L>0$. Let $L_0>-L$ and define
\begin{equation}\label{eq. def L_p}
    L_p:=\chi_p(\lambda_D).
\end{equation}
Then, the operator \(\widehat{\mathcal K}_p\) acts on \(L^2([L_0,L_p])\), as its kernel is defined inside the domain $\chi_p(\mathcal{U}_1)$. Moreover, by Proposition \(\ref{proposition Khat to K_airy}\), for all sufficiently large \(p\), its kernel admits the expression
\[
\begin{aligned}
\widehat K_p(x,y)
={}&K_{\mathrm{Airy}}(x,y)+\widehat C_p p^{-2/3}
\bigl(\mathrm{Ai}(x)\mathrm{Ai}'(y)
+\mathrm{Ai}(y)\mathrm{Ai}'(x)\bigr)
+\widetilde R(x,y),
\end{aligned}
\]
for \((x,y)\in[L_0,L_p]^2\). To work on the unbounded domain $[L_0,\infty)^2$, we extend $\widehat K_p$ to the region outside of $[L_0,L_p]^2$ by setting
\[
\widehat K_p(x,y)=0,
\quad (x,y)\in[L_0,\infty)^2\setminus[L_0,L_p]^2.
\]
We denote by \(\widehat{\mathcal{K}}_p^{(E)}\) the integral operator on \(L^2([L_0,\infty))\) associated with this extended kernel.
Recall also that the kernel \(D_0\), defined in \(\eqref{eq. kernel D0}\), is the analytic extension of the leading correction term of order \(p^{-2/3}\) to the unbounded quadrant. 

The following result shows that, at the level of Fredholm determinants, the contribution of the remaining lower-order terms can be controlled, allowing us to retain only the leading correction \(D_0\). Its proof can be found in Subsection~\ref{proof Prop Fredholm D to D0}

\begin{proposition}\label{lem:bound-D-D0}
Fix $L>0$ and let $L_0>-L$. Define the integral operators
\[
\mathcal{P}:=\left( \mathbb{I}-\mathcal{K}_{Airy}\right)^{-1}, \quad
\mathcal{D}
:=\widehat{\mathcal{K}}_p^{(E)}-\mathcal{K}_{\Airy},
\]
acting on \(L^2([L_0,\infty))\). Then,
\begin{equation}\label{eq. error trunc Ai kernel}
     \det \left( \mathbb{I}- \mathcal{P}\mathcal{D}\right)_{L^2([L_0,\infty))}=\det \left( \mathbb{I}-\mathcal{P} \mathcal{D} _0\right)_{L^2([L_0,\infty))}+\mathcal{O}\left(p^{-1}\right),
\end{equation}
with the implicit constant in $\mathcal{O}\left(p^{-1}\right)$ is independent of $L_0$, for all $L_0>-L$.
\end{proposition}


\begin{proof} [\textbf{Proof of Theorem~\ref{thm: Probability Tracy-Widom}}]

Fix $L>0$ and let $L_0>-L$. Recall that the map $\chi_p$ has local behavior $\xi=\chi_p(\lambda) \approx p^{\frac{2}{3}}(\lambda-1)$ (see~\eqref{eq. chi_p asymtotics}). Then, as $p$ goes to infinity, the image of the disk $\mathcal{U}_1$ covers the entire complex plane. Hence, there exists $p_0 \in \mathbb{N}$ such that, for all $p \geq p_0$,
\[ 
\lambda_0^{(p)} := \chi_p^{-1}(L_0)\qquad \text{and} \qquad\lambda_L^{(p)} := \chi_p^{-1}(-L),
\]
are well defined. Moreover, for all such $p$, the value $L_0^{(p)}:=\varphi_p(\lambda_0^{(p)})$ is also well defined, as $\varphi_p$ is a linear map (see~\eqref{eq:eqchangeofvari zlambda}). Hence, for $p \geq p_0$, relation~\eqref{eq. prob Kp hat} gives
\begin{align}
    \mathbb{P}(\theta_{\max} \leq L_0^{(p)})&= \det(\mathbb{I}-\widehat{\mathcal{K}}_p)_{L^2( [L_0,\,\chi_p(\lambda_D) ])}+\mathcal{O}\left(e^{-kp}\right) \nonumber\\
    &=\det(\mathbb{I}-\widehat{\mathcal{K}}^{(E)}_p)_{L^2( [L_0,\,\infty))}+\mathcal{O}\left(e^{-kp}\right) \nonumber\\   
     &= \det \left( \mathbb{I}-\mathcal{K}_{Airy} \right)_{L^2([L_0,\infty))} \, \det \left( \mathbb{I}-\left( \mathbb{I}-\mathcal{K}_{Airy} \right)^{-1} \mathcal{D} \right)_{L^2([L_0,\infty))}+\mathcal{O}\left(e^{-kp}\right)\\   
     &= \det \left( \mathbb{I}-\mathcal{K}_{Airy} \right)_{L^2([L_0,\infty))} \, \det \left( \mathbb{I}-\left( \mathbb{I}-\mathcal{K}_{Airy} \right)^{-1} \mathcal{D}_0 \right)_{L^2([L_0,\infty))}+\mathcal{O}\left(p^{-1}\right). \label{eq:probtheta}
\end{align} 
In the system \eqref{eq:probtheta}, the third equality uses the multiplicative property of the Fredholm determinant~\cite[Theorem~3.8]{simon1977notes}, and the last equality follows from Proposition~\ref{lem:bound-D-D0}. Moreover, we can write
\begin{align}
\det \left( \mathbb{I}-\left( \mathbb{I}-\mathcal{K}_{Airy} \right)^{-1} \mathcal{D}_0 \right)_{L^2([L_0,\infty))} &= \mathrm{exp} \left[ \log \left( \det \left( \mathbb{I}-\left( \mathbb{I}-\mathcal{K}_{Airy} \right)^{-1} \mathcal{D}_0 \right)_{L^2([L_0,\infty))} \right) \right] \label{e:log_det}\\
&= \mathrm{exp} \left[ \mathrm{Tr} \left( \log \left( \mathbb{I}-\left( \mathbb{I}-\mathcal{K}_{Airy} \right)^{-1} \mathcal{D}_0  \right)_{L^2([L_0,\infty))} \right) \right].\label{e:log_det_Tr}
\end{align}
In \eqref{e:log_det}, the logarithm operator is defined by
\begin{equation}\label{eq. Taylor log}
\log\left(
\mathbb{I}-
\left(\mathbb{I}-\mathcal{K}_{\mathrm{Airy}}\right)^{-1}
\mathcal{D}_0
\right)
=
-\sum_{j=1}^{\infty}\frac{1}{j}
\left[
\left(\mathbb{I}-\mathcal{K}_{\mathrm{Airy}}\right)^{-1}
\mathcal{D}_0
\right]^j.
\end{equation}
Note that \eqref{eq. Taylor log} is well defined for $p$ sufficiently large so that $\rVert \left(\mathbb{I}-\mathcal{K}_{\mathrm{Airy}}\right)^{-1}
\mathcal{D}_0\rVert_{op}<1$. Next recall that, by Proposition~\ref{prop prop53},
\begin{equation}\label{e:Tr(I-K-Airy)^(-1)D0}
\mathrm{Tr} \left( \left( \mathbb{I}-\mathcal{K}_{Airy} \right)^{-1} \mathcal{D}_0 \right) = 2\,\widehat{C}_p\, p^{-2/3}\,v(s).
\end{equation}
By combining \eqref{e:log_det_Tr}, \eqref{eq. Taylor log} and \eqref{e:Tr(I-K-Airy)^(-1)D0}, we obtain
\begin{align*}
\det \left( \mathbb{I}-\left( \mathbb{I}-\mathcal{K}_{Airy} \right)^{-1} \mathcal{D}_0 \right)_{L^2([L_0,\infty))} &= \mathrm{exp} \left[-2\,\widehat{C}_p\, p^{-2/3}\,v(s) +  \mathcal{O} \left( p^{-4/3} \right) \right].
\end{align*}
Therefore, by relation~\eqref{eq:probtheta},
\begin{align*}
  \mathbb{P}(\theta_{\max} \leq L_0^{(p)})&= \det \left( \mathbb{I}-\mathcal{K}_{Airy} \right)_{L^2([L_0,\infty))} \, \mathrm{exp} \left(-2\, \widehat{C}_p\, p^{-2/3} \,v(L_0) + \mathcal{O} \left( p^{-4/3} \right)\right) +\mathcal{O}\left(p^{-1}\right)\\
  &=F_2(L_0) \left( 1 -2\, \widehat{C}_p \, p^{-2/3} v(L_0) \right) + \mathcal{O}\left( p^{-1} \right).
\end{align*}
This completes the proof of Theorem~\ref{thm: Probability Tracy-Widom}.
\end{proof}

\begin{remark}
The function \(v(s)\) appearing in the statement of the Theorem \(\ref{thm: Probability Tracy-Widom}\) coincides with the function defined through the inner product in \(\eqref{eq. def v}\). Its explicit connection to the Painlevé II equation will be established in Section~\ref{Proof Prop 5.4}.
\end{remark}


\begin{proof}  [\textbf{Proof of Corollary~\ref{mainthrm}}]
Fix an arbitrary $s\in\mathbb{R}$ and consider $L>0$ such that $s>-L$. We start by deriving an asymptotic expression for
\[
s^{(p)}:= \varphi_p\bigl(\chi_p^{-1}(s)\bigr).
\]
Recall the definition of the local map in \eqref{eq:eqchangeofvari zlambda}:
\[\chi_p(\lambda)=p^{2/3} \left( -\frac{3}{2} \phi(\lambda)\right)^{2/3}.\]
The local behavior of the function $\phi(\lambda)$ near $\lambda=1$ (see \eqref{eq local}) gives:
\begin{align*}
\chi_p(\lambda)= p^{2/3}\left( -h_p(1)\right)^{2/3} (\lambda-1)\left( 1+\mathcal{O}(\lambda-1)\right)^{2/3}.
\end{align*}
Inverting the last one gives
\begin{align}\label{eq invere}
\chi_p^{-1}\left( s\right) = 1+\frac{s}{p^{2/3}\left( -h_p(1)\right)^{2/3} }\left( 1+\mathcal{O}(s-1)\right)^{-2/3},
\end{align}
Since $\chi_p'(1)=p^{2/3}(-h_p(1))^{2/3}$, with $h_p(1)$ bounded away from zero for sufficiently large $p$,
the local inverse satisfies $s-1=\mathcal{O}(p^{-2/3})$, and therefore we can write
\begin{align}
\chi_p^{-1}\left( s\right) = 1+\frac{s}{p^{2/3}\left( -h_p(1)\right)^{2/3} }\left( 1+\mathcal{O}\left(p^{-2/3}\right)\right).
\end{align}
Next, we apply the transformation $z=\varphi_p(\lambda)$ (see~\eqref{eq:eqchangeofvari zlambda}) to \eqref{eq invere}. This gives:
\begin{align}\label{inv map}
\varphi_p(\chi_p^{-1}\left( s\right))=\beta_p + \check C_p p^{-2/3} s \left( 1+\mathcal{O}\left(p^{-2/3}\right)\right), \quad \check C_p = \frac{\beta_p - \alpha_p}{2 \left( -h_p(1)\right)^{2/3}},
\end{align}
which is valid for $\nu \geq 0$. The coefficient $\check C_p$ depends on $p$ through $\alpha_p$, $\beta_p$, and $h_p(1)$. Its asymptotic behavior differs depending on whether $\nu>0$ or $\nu=0$. We consider these two cases separately.
\begin{itemize}
\item[(i)] $\nu >0$: substituting the asymptotic expansions of the endpoints $\alpha_p$, $\beta_p$ given in \eqref{eq:eq237}, and of the quantity $h_p(1)$ (see \eqref{eq:consh1}):
\begin{align*}
h_p(1) = h_{\infty} + \mathcal{O}(p^{-\nu}), \quad h_{\infty}=-\frac{2 \sqrt{2}(1+b_{\infty})}{2+b_{\infty}+2\sqrt{1+b_{\infty}}}
\end{align*}    
we get the following asymptotic result for $\check C_p$:
\begin{align*}
\check C_p = \frac{2 \sqrt{1+b_{\infty}}}{a_{\infty}} \, \left(-h_{\infty}\right)^{-2/3} \, p^{-\nu} + \mathcal{O}(p^{-2 \nu}), \quad p \to \infty.
\end{align*}
We then rewrite \eqref{inv map} as follows:
\begin{align}\label{inv map 1}
\varphi_p(\chi_p^{-1}\left( s\right))=\beta_p + \frac{s}{C \, p^{\nu+2/3}} + \mathcal O(p^{-2\nu-2/3})+\mathcal O(p^{-\nu-4/3})+\mathcal O(p^{-2\nu-4/3}).
\end{align}
Here, \(C\) is the constant given in Corollary~\ref{mainthrm} for $\nu>0$, and is the reciprocal of the leading coefficient in the asymptotic expansion of \(\check C_p\).
\item[(ii)] $\nu=0$: in this case $\check C_p$ is a fixed constant. We drop the $p$-dependence and call it $\check C$:
\[\check C=\frac{2\sqrt{(1+a_{\infty})(1+b_{\infty})(1+a_{\infty}+b_{\infty})}}{(2+a_\infty+b_{\infty})^2} \, \left( -h(1)\right)^{-2/3}\]
where $h(1)$ is the fixed constant given in equation \eqref{eq:h1}. We rewrite \eqref{inv map} as follows:
\begin{equation}
\label{inv map 0}
\varphi_p(\chi_p^{-1}\left( s\right))=\beta + \frac{s}{C \, p^{2/3}} + \mathcal{O}(p^{-4/3}), \quad p \to \infty,
\end{equation}
where $C$ is the constant provided in Corollary~\ref{mainthrm} for $\nu=0$, and is the reciprocal of $\check C$.
\end{itemize}
In conclusion, from~\eqref{inv map 1} and \eqref{inv map 0} we obtain:
\begin{align}\label{inv map 2}
\varphi_p(\chi_p^{-1}\left(s\right))=\beta_p+ \frac{s}{C\,p^{\nu+2/3}} +\begin{cases}
\mathcal{O}\left(p^{-\nu-2/3-\min\{\nu,2/3\}}\right),
& \nu>0,\\[2mm]
\mathcal{O}\left(p^{-4/3}\right),
& \nu=0,
\end{cases}
\end{align}
where $C$ is the constant appearing in Corollary~\ref{mainthrm}.

Now, let \(\delta>0\) be sufficiently small so that $\chi_p^{-1}$ is well defined for $s\pm\delta$. It follows from \eqref{inv map 2} that, for all sufficiently large \(p\),
\[
(s-\delta)^{(p)}
<
\beta_p+\frac{s}{C\,p^{\nu+2/3}}
<
(s+\delta)^{(p)}.
\]
Therefore,
\begin{equation} \label{eq: TW pm delta}
  \mathbb{P}\bigl(\theta_{\max}\leq(s-\delta)^{(p)}\bigr)
\leq
\mathbb{P}\left(
C\,p^{\nu+2/3}(\theta_{\max}-\beta_p)\leq s
\right)
\leq
\mathbb{P}\bigl(\theta_{\max}\leq (s+\delta)^{(p)}\bigr).  
\end{equation}
On the other hand, for $L_0>-L$, Theorem~\ref{thm: Probability Tracy-Widom} implies, in particular, that
\begin{equation}
\mathbb{P}\bigl(\theta_{\max}\leq \varphi_p\bigl(\chi_p^{-1}(L_0)\bigr)\bigr) = F_2(L_0)+\mathcal{O}(p^{-2/3}),
\quad p\to\infty.
\end{equation}
Hence, taking \(p\to\infty\) in \eqref{eq: TW pm delta}, we obtain
\[
F_2(s-\delta)
\leq
\liminf_{p\to\infty}
\mathbb{P}\left(
C\,p^{\nu+2/3}(\theta_{\max}-\beta_p)\leq s
\right)
\leq
\limsup_{p\to\infty}
\mathbb{P}\left(
C\,p^{\nu+2/3}(\theta_{\max}-\beta_p)\leq s
\right)
\leq
F_2(s+\delta).
\]
Since the c.d.f.\ $F_2$ is continuous, letting \(\delta\to0\), we conclude that
\[
\lim_{p\to\infty}
\mathbb{P}\left(
C\,p^{\nu+2/3}(\theta_{\max}-\beta_p)\leq s
\right)
=
F_2(s).
\]
\end{proof}


\subsection[Largest eigenvalue away from the endpoint]{Largest eigenvalue away from the endpoint $\beta_p$}\label{subsec eigval stats case 1}

We now consider the case where $a>\varphi_p(\lambda_D))$. We start by proving that in this case, the probability of finding eigenvalues in the interval $(a,1)$ decays exponentially fast as $p\to\infty$.

\begin{lemma}\label{Lemma: right of lambdaD}
Suppose $a \in \mathbb{R}$ is such that $a>\varphi_p(\lambda_D))$.  Then, there exists a constant $k>0$ such that
\begin{equation}\label{e:P(theta_max>a_p)}
\mathbb{P}\!\left(\theta_{\max}>a\right)
=\mathcal{O}\left(e^{-kp}\right),
\end{equation}
where $k$ and the implicit constant in
$\mathcal{O}\left(e^{-kp}\right)$ are uniform for all
$a\in\bigl(\varphi_p(\lambda_D),1\bigr)$.
\end{lemma}

\begin{proof}
From relation~\eqref{eq: prob Ktilde} we have
\begin{gather}\label{eq: prob lambda max greater}
     \mathbb{P}(\theta_{\max}>a)=1-\det( \mathbb{I}-\widetilde{\mathcal{K}}_p)_{L^2([\varphi_p^{-1}(a),\varphi_p^{-1}(1)])}
\end{gather}
where
\begin{gather}\label{eq: Fredholm determinant}
\det\left( \mathbb{I} - \widetilde{\mathcal{K}}_{p} \right)_{L^2([\varphi_p^{-1}(a),\varphi_p^{-1}(1)])} = 1 + \sum_{j=1}^{p} \frac{(-1)^{j}}{j!} \int_{\varphi_p^{-1}(a)}^{{\varphi_p^{-1}(1)}} \cdots \int_{\varphi_p^{-1}(a)}^{{\varphi_p^{-1}(1)}} \det_{i,k=1,2,\ldots, j} \left[ \widetilde{K}_{p} \left( \lambda_i, \lambda_k \right) \right] d^j\lambda.\hspace{1cm} 
\end{gather}
We now consider the integration with respect to $\lambda_1$ in \eqref{eq: Fredholm determinant}. That is,
\begin{align}\label{e:eqdet}
 \int_{\varphi_p^{-1}(a)}^{{\varphi_p^{-1}(1)}} & \det_{i,k=1,2,\ldots, j} \left[ \widetilde{K}_{p} \left( \lambda_i, \lambda_k \right) \right] d \lambda_1 = \int_{\varphi_p^{-1}(a)}^{{\varphi_p^{-1}(1)}} \det \begin{vmatrix}
\widetilde{K}_p(\lambda_1,\lambda_1) & \cdots & \widetilde{K}_p(\lambda_1,\lambda_j)\\
\vdots & \cdots & \vdots \\
\widetilde{K}_p(\lambda_j,\lambda_1) & \cdots & \widetilde{K}_p(\lambda_j,\lambda_j)
\end{vmatrix} \, d \lambda_1  \nonumber \\
& =  \sum_{i=1}^{j} (-1)^{1+i} \int_{\varphi_p^{-1}(a)}^{\varphi_p^{-1}(1)} \widetilde{K}_p(\lambda_1,\lambda_i) \det \begin{vmatrix}
\widetilde{K}_p(\lambda_1,\lambda_1) & \cdots & \widetilde{K}_p(\lambda_1,\lambda_j)\\
\vdots & \cdots & \vdots \\
\widetilde{K}_p(\lambda_j,\lambda_1) & \cdots & \widetilde{K}_p(\lambda_j,\lambda_j)
\end{vmatrix}^{(1,i)} \, d \lambda_1,
\end{align}
where the subscript $(1,i)$ on the determinant of the last equation indicates that the first row and the $i$-th column are omitted.

Note that $\lambda_1, \lambda_i \in (\varphi_p^{-1}(a), \varphi_p^{-1}(1))\subset(\lambda_D,\varphi_p^{-1}(1))$, with the length of the former interval being of order $p^{\nu}$. By Lemma~\ref{lem:lemma2}, for $\lambda_1, \lambda_2 \in (\lambda_D, \varphi_p^{-1}(1))$, there exists a constant $k > 0$ such that
\begin{equation}\label{e:K-tilde_p(lambda1,lambda2)=O(exp(-kp))}
\widetilde{K}_p(\lambda_1,\lambda_2)=  \mathcal{O}(e^{-k p}).
\end{equation}
Then, relation \eqref{e:P(theta_max>a_p)} follows immediately from \eqref{e:K-tilde_p(lambda1,lambda2)=O(exp(-kp))} applied in \eqref{e:eqdet}, together with \eqref{eq: prob lambda max greater} and \eqref{eq: Fredholm determinant}. 
\end{proof}

\begin{proof} [\textbf{Proof of Theorem~\ref{prop:prop7}}]
Let $c$ be any arbitrary positive constant, and suppose that for $p$ sufficiently large, $\theta \in (0,1)$ satisfies $$\theta-\beta_p\geq c\,p^{-\nu}.$$
Recall that $\lambda_D$ denotes the rightmost intersection of the boundary of $\mathcal{U}_1$ with the real axis. Here, $\mathcal{U}_1$ is a disk of radius $\epsilon_0$, centered at $1$, on which the local map $\chi_p$ is defined (see Remark~\ref{remark epsilon}). Also, recall from \eqref{eq:eq237} that
\[p^\nu(\beta_p-\alpha_p)\to \frac{4\sqrt{1+b_\infty}}{a_\infty}, \qquad\text{as}\quad p \to \infty.\] 

There are two cases to be considered.

\begin{itemize}
    \item If $c>\epsilon_0\,\frac{2\sqrt{1+b_\infty}}{a_\infty},$ then, for all suficiently large $p$
    \[ \varphi_p^{-1}(\theta) =  1+\frac{2(\theta-\beta_p)}{\beta_p-\alpha_p} > 1+\epsilon_0 = \lambda_D.\]

    \item  If $c\leq\epsilon_0\,\frac{2\sqrt{1+b_\infty}}{a_\infty}$, the  radius of the disk $\mathcal U_1$ can be redefined as $\epsilon_1 :=\frac{c\,a_\infty}{4\sqrt{1+b_\infty}}<\epsilon_0$ (see Remark~\ref{remark epsilon}). Let $\widetilde{\lambda}_D:=1+\epsilon_1$ be the new rightmost intersection of the boundary of $\mathcal U_1$ with the real axis.
    As in the previous case, for all sufficiently large $p$,
    \[\varphi_p^{-1}(\theta) = 1+\frac{2(\theta-\beta_p)}{\beta_p-\alpha_p} > 1+\epsilon_1 = \widetilde{\lambda}_D.\]
\end{itemize}
Thus, in either case, $\varphi_p^{-1}(\theta)$ lies to the right of the rightmost intersection of the boundary of the disk centered at $1$, on which the local map $\chi_p$ is defined. Hence,  Lemma~\ref{Lemma: right of lambdaD} implies that there exist
constants $K,k>0$ such that, for all sufficiently large $p$,
\[
\mathbb{P}\!\left(\theta_{\max}>\theta\right)
\leq
Ke^{-kp}.
\]
\end{proof}

\section{Riemann-Hilbert analysis}\label{section RH Analysis}

In this section, we provide a detailed analysis of RHP~\ref{RHP Y}, with the primary goal of deriving the asymptotic behavior, as $p \to \infty$, of the Jacobi polynomials associated to the weight in \eqref{eq:eqweight}. We follow the classical Riemann–Hilbert  analysis techniques from \cite{deift:kriecherbauer:mclaughlin:venakides:zhou:1999, kuijlaars2004riemann,Deift2000}. These asymptotics will be applied to study the distribution of the largest eigenvalue. We begin with a sequence of transformations that reduce RHP~\ref{RHP Y} to a model problem, whose solution provides the outer approximation valid away from the endpoints $\lambda=\pm 1$. Near these endpoints, we construct local approximations. The global asymptotics are then obtained by controlling the error between the global solution and these approximations through a small-norm RHP.


\subsection[The g-function]{The $g$-function}
\label{g function}

We recall that RHP~\ref{RHP Y} is not in standard form, since the function $Y(z)$ does not approach the identity matrix as $z \rightarrow \infty$. To address this issue, define
\begin{equation}\label{eq: transformation T}
    T(z) = e^{-\frac{p \ell}{2} \sigma_3} Y(z) e^{-p \left( g(z) - \frac{\ell}{2} \right) \sigma_3},
\end{equation}
where the function $g$ and the constant $\ell$ will be determined below (see \eqref{Eqtn g} and \eqref{e:eq22}). Then, $T(z)$ solves the following RHP.
\begin{RHP}
Find a $2 \times 2$ matrix-valued function $T(z)$ that satisfies the following conditions.
\begin{enumerate}
    \item $T(z)$ is analytic in $\mathbb{C} \setminus (0,1)$.
    \item $T(z)$ satisfies the jump relation 
    \begin{equation}
       T_{+}(x) = T_{-}(x) J_{T}(x), \quad J_T(x)=\pmtwo{e^{-p(g_+-g_-)}}{e^{p(g_{+} + g_{-} - \ell ) + p \, b_{\infty} \log x + p^{\nu +1} \, a_{\infty} \log(1-x)}}{0}{e^{p(g_+-g_-)}}, \quad x \in (0,1).
    \end{equation}
    \item $T(z)$ has the following large $z$ asymptotic behavior:
    \begin{equation}\label{e:asymT}
        T(z) = I_2 + \mathcal{O}(z^{-1}), \quad z \to \infty.
    \end{equation}
\end{enumerate}
\end{RHP}
Following \cite{DKMVZ3}, we seek $0<\alpha_p<\beta_p<1$,  and a function $g$ satisfying the following conditions:
\begin{subequations}\label{e:eq7}
\begin{equation}\label{e:eq7a}
g_{+} (x) + g_{-}(x) = - b_{\infty} \, \log{x} -a_{\infty} \, p^{\nu}\log{ (1 - x) }  + \ell \ \ , \ x \in (\alpha_p, \beta_p) \,
\end{equation}
\begin{equation}\label{e:eq7b}
g_{+}(x) - g_{-}(x) = 0, \ \ x > \beta_p, \ \ \ \ g_{+}(x) - g_{-}(x) = 2 \pi i, \ \ x < \alpha_p 
\end{equation}
\end{subequations}
and the inequality
\begin{equation}\label{e:eq9}
g_{+} (x) + g_{-}(x) \le  -  b_{\infty} \, \log{x} - a_{\infty} \, p^{\nu} \log{ (1 - x) }  + \ell, \end{equation}
holding for $x$ in the entire interval $(0,1)$, and satisfying the large-$z$ asymptotic behavior
 \begin{equation}\label{Eqtn asymptotics g}
g(z) = \log{z} + \mathcal{O}(z^{-1}), \quad z \to \infty.
\end{equation}
Note that we require condition \eqref{Eqtn asymptotics g} so that \eqref{e:asymT} follows.
\begin{proposition}
Let $\alpha_p$ and $\beta_p$ be the two real numbers defined in \eqref{eq. alpha_beta_explicit}. Consider the function 
\begin{equation*}
 R(z)=(z-\alpha_p)^{1/2}(z-\beta_p)^{1/2}, \quad z\in \mathbb{C}\,\setminus\,[\alpha_p,\beta_p],   
\end{equation*}
where the branch is chosen so that $R$ is analytic in $\mathbb{C}\,\setminus\,[\alpha_p,\beta_p]$ and satisfies $R(z)\sim z$ as $z\rightarrow \infty$ . Then, the $g$-function defined as 
\begin{equation}\label{Eqtn g}
g(z) =c+ \frac{1}{2} \int_{\beta_p}^z \Bigg(  -\frac{b_{\infty}}{s} + \frac{a_{\infty} \, p^{\nu}}{1-s} + R(s)\left(\frac{b_{\infty}}{s \, R(0)} + \frac{a_{\infty} \, p^{\nu}}{(s-1) \, R(1)}\right) \Bigg)\,ds
\end{equation}
satisfies conditions \eqref{e:eq7}-\eqref{Eqtn asymptotics g}. Here, the constant $c$ is given by
\begin{equation}
c = \log \beta_p - \int_{\beta_p}^{\infty} \left( g'(s) - \frac{1}{s} \right) \,ds,
\end{equation}
and the parameter $\ell$ is given by
\begin{equation}\label{e:eq22}
\ell=b_{\infty}\log \beta_p+ a_{\infty} \, p^{\nu} \log(1-\beta_p)+2c.
\end{equation}
\end{proposition}

\begin{proof}
Let $g$ be a function satisfying conditions \eqref{e:eq7}-\eqref{Eqtn asymptotics g}. Then $g'$ satisfies the conditions
\begin{subequations}\label{e:eq10}
\begin{equation}\label{e:eq10a}
    g'_{+}(x)+g'_{-}(x) = - \frac{b_{\infty}}{x} + \frac{a_{\infty} \, p^{\nu}}{1-x}, \quad  x \in (\alpha_p, \beta_p) 
\end{equation}
\begin{equation}\label{e:eq10b}
    g'(z) = \frac{1}{z}, \quad z \to \infty.
\end{equation}
\end{subequations}
Define $g'(z) = F(z)R(z)$. Since for $x \in (\alpha_p, \beta_p)$, the function $R$ satisfies $R_{+}(x) = -R_{-}(x)$,  we have 
\begin{equation*}
   g'_{+}(x)+g'_{-}(x) =R_+(x)(F_+(x)-F_-(x)),\quad x \in (\alpha_p, \beta_p).
\end{equation*}
Hence, by equation \eqref{e:eq10a}
\begin{eqnarray}
    F_+(x)-F_-(x)=\frac{1}{R_+(x)}\left(- \frac{b_{\infty}}{x} + \frac{a_{\infty} \, p^{\nu}}{1-x}\right)\ , \ x \in (\alpha_p, \beta_p) .
\end{eqnarray}
Also, since $F(z) = \frac{g'(z)}{R(z)}$, eq.~\eqref{e:eq10b} implies that
\begin{eqnarray}
F(z) = \frac{1}{z^2}, \quad z \to \infty.
\end{eqnarray}
From the Sokhotski-Plemelj formula, it then follows that for $z\in \mathbb{C}\setminus (\alpha_p,\beta_p)$
\begin{align*}
    F(z)&=\frac{1}{2\pi i} \int_{\alpha_p}^{\beta_p} \frac{1}{R_+(s)}\left(- \frac{b_{\infty}}{s} + \frac{a_{\infty} \, p^{\nu}}{1-s}\right)\frac{1}{s-z} \,ds\\
    &=\frac{1}{4\pi i} \int_{\Gamma} \frac{1}{R(s)}\left(- \frac{b_{\infty}}{s} + \frac{a_{\infty} \, p^{\nu}}{1-s}\right)\frac{1}{s-z}\,ds,
\end{align*}
where $\Gamma$ is a closed contour shown in Figure \ref{fig Contour Gamma}.
\begin{figure}[htbp]
    \centering
    \begin{tikzpicture}[scale=0.8]
        \draw[<->] (-4,0) -- (4,0) node[below]{};
        
        \node at (-3,0) [circle,fill,inner sep=1pt,label=below:0] {};
        \node at (3,0) [circle,fill,inner sep=1pt,label=below:1] {};
        \node at (-2,0) [circle,fill,inner sep=1pt,label=below: $\alpha_p$] {};
        \node at (1,0) [circle,fill,inner sep=1pt,label=below: $\beta_p$] {};
        \node at (0.5,1.2) [circle,fill,inner sep=1pt,label=above: $z$] {};
        
        \node at (-1.5,0.7) [above] {$\Gamma$};

        \draw[line width=1pt,blue!50!, decoration={
            markings, mark=at position 0.5 with {\arrow[line width=1.2pt]{>} } }, postaction={decorate}]
        (-2,0) to[out=15,in=165] (1,0);
        
        \draw[line width=1pt,blue!50!, decoration={
            markings, mark=at position 0.5 with {\arrow[line width=1.2pt]{<}} }, postaction={decorate}]
        (-2,0) to[out=-15,in=-165] (1,0); 
        
    \end{tikzpicture}
    
    \caption{Contour $\Gamma$.}
    \label{fig Contour Gamma}
\end{figure}

The contour $\Gamma$ can be deformed to a contour $\Gamma^*$ consisting of small loops around the poles at $0$, $1$, $z$ and a large loop enclosing all three poles (see Figure~\ref{fig Contour Gamma*}). 
\begin{figure}[htbp]
    \centering
    \begin{tikzpicture}[scale=0.6]
        \draw[<->] (-4,0) -- (4,0) node[below]{};
        
        \node at (-3,0) [circle,fill,inner sep=1pt,label=below:0] {};
        \node at (3,0) [circle,fill,inner sep=1pt,label=below:1] {};
        \node at (-2,0) [circle,fill,inner sep=1pt,label=below: $\alpha_p$] {};
        \node at (1,0) [circle,fill,inner sep=1pt,label=below: $\beta_p$] {};
        \node at (0.5,1.2) [circle,fill,inner sep=1pt,label=above: $z$] {};
        
        \node at (-1.5,0.7) [above] {$\Gamma^*$};
        
        \draw[line width=1pt, blue!50!, decoration={
            markings,
            mark=at position 0.25 with {\arrow[line width=1.2pt]{>}},
        }, postaction={decorate}] 
        (-3,0) circle (0.3);

        \draw[line width=1pt, blue!50!, decoration={
            markings,
            mark=at position 0.25 with {\arrow[line width=1.2pt]{>}},
        }, postaction={decorate}] 
        (3,0) circle (0.3);

        \draw[line width=1pt,blue!50!, decoration={
            markings,
            mark=at position 0.75 with {\arrow[line width=1.2pt]{>}},
        }, postaction={decorate}] 
        (0.5,1.2) circle (0.3);

       \draw[line width=1pt,blue!50!, decoration={
            markings,
            mark=at position 0.2 with {\arrow[line width=1.2pt]{<}},
            mark=at position 0.8 with {\arrow[line width=1.2pt]{<}},
        }, postaction={decorate}] 
        (0,0) circle (3.5);

    \end{tikzpicture}
    
    \caption{Contour $\Gamma^*$.}
    \label{fig Contour Gamma*}
\end{figure}

By the residue theorem and noting that the integral over the large loop vanishes as its radius tends to infinity, we obtain
\begin{eqnarray}
F(z) =\frac{1}{2} \Bigg( \frac{1}{R(z)} \left( -\frac{b_{\infty}}{z} + \frac{a_{\infty} \, p^{\nu}}{1-z} \right) + \frac{b_{\infty}}{z \, R(0)} + \frac{a_{\infty} \, p^{\nu}}{(z-1) \, R(1)} \Bigg).
\end{eqnarray}
Hence,  the definition of $g'(z)$ yields
\begin{equation}\label{e:eq17}
g'(z) = F(z)R(z) =\frac{1}{2} \Bigg(  -\frac{b_{\infty}}{z} + \frac{a_{\infty} \, p^{\nu}}{1-z} + R(z)\left(\frac{b_{\infty}}{z \, R(0)} + \frac{a_{\infty} \, p^{\nu}}{(z-1) \, R(1)}\right) \Bigg).
\end{equation}
After direct integration, this gives \eqref{Eqtn g}. 

Here, we need to verify that $g'$ satisfies the condition $g'(z)=1/z$, as $z\rightarrow\infty$ (see eq.~\eqref{e:eq10b}). This condition determines the choice of the values of $\alpha_p$ and $\beta_p$. Indeed, large-$z$ asymptotics of eq. \eqref{e:eq17} gives
\begin{eqnarray}
g'(z) = \frac{1}{2} \left( \frac{b_{\infty}}{R(0)} + \frac{a_{\infty} \, p^{\nu}}{R(1)} \right) - \frac{1}{4} \left[ \frac{b_{\infty}}{R(0)} \left( 2 R(0) + \alpha_p + \beta_p \right) + \frac{a_{\infty} \, p^{\nu}}{R(1)} \left( -2 + 2 R(1) + \alpha_p + \beta_p \right) \right] \frac{1}{z} + \mathcal{O}(z^{-2}).
\end{eqnarray}
Then, to match the asymptotics of $g'$ in eq. \eqref{e:eq10b} we require
\begin{subequations}
\begin{equation}\label{EQ for alpha_beta1}
\frac{b_{\infty}}{R(0)} + \frac{a_{\infty} \, p^{\nu}}{R(1)} = 0,
\end{equation}
\begin{equation}\label{EQ for alpha_beta2}
- \frac{1}{4} \left[ \frac{b_{\infty}}{R(0)} \left( 2 R(0) + \alpha_p + \beta_p \right) + \frac{a_{\infty} \, p^{\nu}}{R(1)} \left( -2 + 2 R(1) + \alpha_p + \beta_p \right) \right] = 1.
\end{equation}
\end{subequations}
This system of equations determines the values $\alpha_p$ and $\beta_p$:
\begin{equation*}
\begin{split}
  \alpha_p & =\frac{-(a^2_{\infty} p^{2\nu}-b^2_{\infty}-\eta) - \sqrt{(a^2_{\infty} p^{2\nu}-b^2_{\infty}-\eta)^2 - 4\eta\, b^2_{\infty}}}{2\eta},\\
  \beta_p& =\frac{-(a^2_{\infty} p^{2\nu}-b^2_{\infty}-\eta) + \sqrt{(a^2_{\infty} p^{2\nu}-b^2_{\infty}-\eta)^2 - 4\eta\, b^2_{\infty}}}{2\eta},  
\end{split}
\end{equation*}
where $\eta = (b_{\infty} + a_{\infty} p^{\nu} + 2)^2$.  

Now, we compute the constant $c$ in \eqref{Eqtn g}. For this we use the asymptotics of $g$ as $z \to \infty$. Integrating \(g'(s) - 1/s\) from \(\beta_p\) to \(z\) yields
\begin{align*}
g(z) &= g(\beta_p) + \int_{\beta_p}^z \left( g'(s) - \frac{1}{s} \right) ds + \log z - \log \beta_p .\\
&= g(\beta_p) + \log z - \log \beta_p + \int_{\beta_p}^{\infty} \left( g'(s) - \frac{1}{s} \right) \,ds - \int_{z}^{\infty} \left( g'(s) - \frac{1}{s} \right) \,ds.
\end{align*}
To match \eqref{Eqtn asymptotics g} we need 
\[g(\beta_p) = \log \beta_p -\int_{\beta_p}^{\infty} \left( g'(s) - \frac{1}{s} \right) \,ds,\]
which according to \eqref{Eqtn g}, must be exactly the value of $c$. Finally, to compute the value of $\ell$, note that for \( x \in (\alpha_p,\beta_p) \) equation \eqref{Eqtn g} yields
\begin{align*}
    g_+(x) + g_-(x) &= 2c + \int_{\beta_p}^x \left(-\frac{b_{\infty}}{s} + \frac{a_{\infty} \, p^{\nu}}{1-s}\right) \, ds \\
    &=- b_{\infty} \log x- a_{\infty} \, p^{\nu} \log(1-x)+ b_{\infty} \log \beta_p + a_{\infty} \, p^{\nu} \log(1-\beta_p)+2c.
\end{align*}
Comparing with \eqref{e:eq7a} we obtain equation \eqref{e:eq22}.
\end{proof}

\begin{remark}\label{remark asymp alpha beta}
If $\nu > 0$, the endpoints $\alpha_p$ and $\beta_p$ of the support of the $g$-function approach zero as $p \to \infty$. Moreover, large-$p$ asymptotics for $\alpha_p$ and $\beta_p$ yield
 \begin{equation}\label{eq:eq237}
\alpha_p = \frac{2 + b_{\infty} - 2 \sqrt{1+b_{\infty}}}{a_{\infty} p^{\nu}} + \mathcal{O}(p^{-2 \nu}), \quad \beta_p = \frac{2 + b_{\infty} + 2 \sqrt{1+b_{\infty}}}{a_{\infty} p^{\nu}} + \mathcal{O}(p^{-2 \nu}).
\end{equation}
\end{remark}

\begin{remark}
When $\nu=0$, the endpoints $\alpha_p$ and $\beta_p$ are independent of $p$, and we therefore omit the subscript $p$ in their notation. They are given explicitly by
\begin{equation}\label{eq:abzero}
\begin{split}
  \alpha & = \frac{2+(a_{\infty}+b_{\infty})(2+b_{\infty})-2\sqrt{(1+a_{\infty})(1+b_{\infty})(1+a_{\infty}+b_{\infty})}}{(2+a_{\infty}+b_{\infty})^2},\\
  \beta & = \frac{2+(a_{\infty}+b_{\infty})(2+b_{\infty})+2\sqrt{(1+a_{\infty})(1+b_{\infty})(1+a_{\infty}+b_{\infty})}}{(2+a_{\infty}+b_{\infty})^2}.
\end{split}
\end{equation}
Moreover, the fixed distance between $\alpha$ and $\beta$ is:
\[
\mathrm{dist}(\alpha,\beta) = \frac{4 \sqrt{(1+a_{\infty})(1+b_{\infty})(1+a_{\infty}+b_{\infty})}}{(2+a_{\infty}+b_{\infty})^2}.
\]
\end{remark}


\subsection[Reduction of the jump matrix and the phi-function]{Reduction of the jump matrix and the $\phi$-function}

In this section, following the approach of \cite{Deift2000}, we further simplify the jump matrix associated with RHP~\ref{RHP Y}, thereby preparing the problem for the application of the steepest-descent method and the subsequent construction of local parametrices. For $x \in (0,\alpha_p) \cup (\beta_p,1)$, we get:
\begin{align*}
g_{+}(x) + g_{-}(x) - \ell &= 2c - \ell + \int_{\beta_p}^{x} \left[ -\frac{b_{\infty}}{s} + \frac{a_{\infty} p^{\nu}}{1-s} + R(s) \left( \frac{b_{\infty}}{s R(0)} + \frac{a_{\infty} p^{\nu}}{(s-1) R(1)} \right) \right] \, ds\\
&= - b_{\infty} \log z - a_{\infty} p^{\nu} \log (1-z) + \int_{\beta_p}^{x} R(s) \left( \frac{b_{\infty}}{s R(0)} + \frac{a_{\infty} p^{\nu}}{(s-1) R(1)} \right) \, ds
\end{align*}
where we used eq. \eqref{e:eq22}. Equivalently this gives
\begin{eqnarray}
g_{+}(x) + g_{-}(x) - \ell + b_{\infty} \log x + a_{\infty} p^{\nu} \log (1-x) = \int_{\beta_p}^{x} R(s) \left( \frac{b_{\infty}}{s R(0)} + \frac{a_{\infty} p^{\nu}}{(s-1) R(1)} \right) \, ds, \quad x \in (0,\alpha_p) \cup (\beta_p,1).
\end{eqnarray}
Moreover, eq. \eqref{e:eq7b} implies that
\begin{eqnarray}
e^{g_{+}(x) - g_{-}(x)} = 1, \quad x \in (0,\alpha_p) \cup (\beta_p,1). 
\end{eqnarray}
Now, eq.~\eqref{e:eq7a} implies that
\begin{equation*}
  g_{+}(x) + g_{-}(x) - \ell + b_{\infty} \log x + a_{\infty} p^{\nu} \log (1-x)=0, \quad x \in (\alpha_p, \beta_p).  
\end{equation*}
Finally, using the definition of $g$ in \eqref{Eqtn g} and the fact that for $x \in (\alpha_p, \beta_p)$, $R_{+}(x) = -R_{-}(x)$. Then 
\begin{eqnarray}
 g_{+}(x) - g_{-}(x) = \int_{\beta_p}^{x} R_{+}(s) \left( \frac{A}{s R(0)} + \frac{B}{(s-1) R(1)} \right) \, ds, \quad x \in (\alpha_p, \beta_p).
 \end{eqnarray}
Therefore, the jump matrix for $T$ reduces into:
\begin{eqnarray}
J_{T}(x) = \begin{cases}
\begin{pmatrix}
1 & e^{p \int_{\beta_p}^{x} R(s) \left(\frac{b_{\infty}}{s R(0)} + \frac{a_{\infty} p^{\nu}}{(s-1) R(1)} \right) \, ds} \\
0 & 1
\end{pmatrix},  \quad  x \in (0,\alpha_p)\cup  (\beta_p, 1)\\ \\[1pt]
\begin{pmatrix}
e^{-p \int_{\beta_p}^{x} R_{+}(s) \left(\frac{b_{\infty}}{s R(0)} + \frac{a_{\infty} p^{\nu}}{(s-1) R(1)} \right) \, ds} & 1 \\
0 & e^{p \int_{\beta_p}^{x} R_{+}(s) \left(\frac{b_{\infty}}{s R(0)} + \frac{a_{\infty} p^{\nu}}{(s-1) R(1)} \right) \, ds}
\end{pmatrix},  \quad \alpha_p < x < \beta_p
\end{cases}
\end{eqnarray}
For brevity, we introduce the function
\begin{align}\label{e:phi(z)_def}
\phi (z) = \phi_p (z) := \frac{1}{2} \int_{\beta_p}^{z} R(s) \left(\frac{b_{\infty}}{s R(0)} + \frac{a_{\infty} p^{\nu}}{(s-1) R(1)} \right) \, ds, \quad z \in \mathbb{C} \setminus [-\infty, \beta_p]
\end{align}
where the contour of integration is chosen to avoid the interval $(-\infty,\beta_p]$. The matrix $J_{T}$ then becomes
\begin{equation}\label{eq:eq236}
J_{T}(x) = \begin{cases}
\begin{pmatrix}
1 & e^{2 p \phi (x)} \\
0 & 1
\end{pmatrix},  \quad x \in (0,\alpha_p)\cup  (\beta_p, 1)\\ \\[1pt]
\begin{pmatrix}
e^{-2 p \phi_{+}(x)} & 1 \\
0 & e^{-2 p \phi_{-}(x)}
\end{pmatrix},  \quad \alpha_p < x < \beta_p
\end{cases}.
\end{equation}
The function $\phi(z)$ plays a crucial role in the forthcoming RHP analysis, as will become evident in the next section. Recall from equation \eqref{eq:eq237} that both endpoints $\alpha_p$ and $\beta_p$ collapse to zero as $p \to \infty$. To proceed with the analysis, we recall the linear map $\varphi_p^{-1}: z \longmapsto \lambda$, defined in \eqref{eq:eqchangeofvari zlambda} by: 
\begin{equation}\label{eq:rescale}
    z \coloneqq \varphi_p(\lambda) = \frac{\beta_p + \alpha_p}{2} + \frac{\beta_p - \alpha_p}{2} \, \lambda,
\end{equation}
 which maps $\alpha_p$ and $\beta_p$ to $-1$ and $1$ respectively. To write the function  $\phi$ in terms of $\lambda$, we set \[s(\zeta) = \frac{\beta_p + \alpha_p}{2} + \frac{\beta_p - \alpha_p}{2} \, \zeta.\] Then, $\phi$ becomes 
\begin{equation}\label{eq:funphi}
\phi(\lambda) = \frac{1}{2} \left( \frac{\beta_p - \alpha_p}{2} \right)^2 \int_{1}^{\lambda} (\zeta + 1)^{1/2} (\zeta - 1)^{1/2} \left( \frac{b_{\infty}}{R(0) s(\zeta)} + \frac{a_{\infty} p^{\nu}}{R(1) \big( s(\zeta) -1\big)} \right) \, d \zeta.
\end{equation}
Eq. \eqref{eq:eq237} and the definition of $s$ in terms of $\zeta$ imply that
\begin{subequations}
\begin{equation}\label{eq:eq238a}
\frac{b_{\infty}}{R(0) s(\zeta)} \approx - \frac{a_{\infty}^2}{2 \sqrt{1+b_{\infty}} \left( \zeta + \frac{2+b_{\infty}}{2 \sqrt{1+b_{\infty}}} \right) } p^{2 \nu} + \mathcal{O}(p^{\nu})
\end{equation}
\begin{equation}\label{eq:eq238b}
\frac{a_{\infty} p^{\nu}}{R(1) \big( s(\zeta) -1\big)} \approx - \frac{a_{\infty}^2 p^{2 \nu}}{a_{\infty} p^{\nu}-(2+b_{\infty})+\mathcal{O}(p^{-\nu})} = - a_{\infty} p^{\nu} + \mathcal{O}(1).
\end{equation}
\end{subequations}
Using eqs. \eqref{eq:eq237}, \eqref{eq:eq238a} and \eqref{eq:eq238b}, one has the following large-$p$ asymptotics of $\phi(\lambda)$:
\begin{gather}\label{eq:eq239}
\phi(\lambda) = -\sqrt{1+b_{\infty}} \int_{1}^{\lambda} \frac{(\zeta+1)^{1/2}(\zeta-1)^{1/2}}{ \zeta + \frac{2+b_{\infty}}{2 \sqrt{1+b_{\infty}}}} \, d \zeta + \mathcal{O}(p^{-\nu}).
\end{gather} 

\begin{remark}\label{eq:rem44}
When $0 < \nu < 1$, the endpoints remain at a fixed distance from zero under the transformation \eqref{eq:rescale}, which is essential for carrying out the asymptotic analysis of the polynomials in subsection~\ref{eq:poly}. In contrast, when $\nu=0$, 
the endpoints already remain at a fixed distance from zero in the original variable. While one can still introduce the mapping \eqref{eq:rescale} for convenience (after replacing $\alpha_p$ and $\beta_p$ with $\alpha$ and $\beta$, respectively), in this case, it is not a rescaling but merely a straightforward change of variables.
\end{remark}
\begin{remark}
When $\nu=0$, the function $\phi(\lambda)$ is given by:
\begin{equation}
\phi(\lambda) = \frac{1}{2} \left( \frac{\beta - \alpha}{2} \right)^2 \int_{1}^{\lambda} (\zeta + 1)^{1/2} (\zeta - 1)^{1/2} \left( \frac{b_{\infty}}{R(0) s(\zeta)} + \frac{a_{\infty}}{R(1) \big( s(\zeta) -1\big)} \right) \, d \zeta
\end{equation}
where $\alpha$ and $\beta$ are defined in equation~\eqref{eq:abzero}. In this case, $\phi(\lambda)$ is independent of $p$, and does not admit the asymptotic expansion given in \eqref{eq:eq239}.
\end{remark}


\subsection{Opening lenses and jump matrix factorization}\label{e:sec23}

We next introduce a transformation, referred to as ``opening lenses'', whose main objective is to deform the RHP along contours (referred to as the lenses). This will make rapidly oscillatory terms become exponentially decaying ones, and the corresponding jump matrices approach the identity matrix. This transformation is invertible, so that knowledge of the solution of this new RHP is equivalent to knowledge of the original RHP.

Recall that the matrix $T$, defined in \eqref{eq: transformation T}, satisfies the jump condition given in equation \eqref{eq:eq236}. In the $\lambda$-plane this jump takes the form:
\begin{gather}
J_{T}(\lambda) = \begin{cases}
\begin{pmatrix}
1 & e^{2 p \phi (\lambda)} \\
0 & 1
\end{pmatrix},  \quad \lambda \in (\varphi_p^{-1}(0),-1) \cup (1, \varphi_p^{-1}(1)),\\ \\[1pt]
\begin{pmatrix}
e^{-2 p \phi_{+}(\lambda)} & 1 \\
0 & e^{-2 p \phi_{-}(\lambda)}
\end{pmatrix},  \quad -1 < \lambda < 1,
\end{cases}
\end{gather}
where $\varphi_p^{-1}(0)$ is of order $\mathcal{O}(1)$, while $\varphi_p^{-1}(1)$ is of order $\mathcal{O}(p^{\nu})$.

Now we study the sign of $\Re(\phi(z))$ around the interval \((\alpha_p, \beta_p)\). To do this, write $z=x+iy$. Then, for $ z\in (\alpha_p,\beta_p)$:
\begin{equation}\label{eq:monotone}
-\frac{d}{dx} \Im(\phi_+(z))=-\frac{d}{dx} \Im(\phi_+(x))=-\Im\left( \frac{1}{2} R_+(x) \left(\frac{b_{\infty}}{x R(0)} + \frac{a_{\infty} p^{\nu}}{(x-1) R(1)} \right)\right)>0.
\end{equation}
Hence, by the Cauchy-Riemann equations, we have \(\frac{\partial}{\partial y}\Re(\phi(z)) > 0\), which implies that \(\Re(\phi(z))\) increases in the upward direction, and thus it must be positive in a small neighborhood above the interval \((\alpha_p, \beta_p)\). By a similar argument, we can show that  \(\Re(\phi(z))\) must be positive in a small neighborhood below the interval \((\alpha_p, \beta_p)\).

Using the jump condition of $\phi(\lambda)$ on $(-1, 1)$, we can factorize $J_{T}$ as follows:
\begin{eqnarray}
J_{T}(\lambda) = L(\lambda) \, M \, U(\lambda),  \quad \lambda \in (-1, 1),
\end{eqnarray}
where
\begin{equation}\label{eq:factor}
L(\lambda) = \begin{pmatrix}
1 & 0\\
e^{-2 p \phi_{-}(\lambda)} & 1\\
\end{pmatrix}, \quad M = \begin{pmatrix}
0 & 1\\
-1 & 0
\end{pmatrix}, \quad U(\lambda) = \begin{pmatrix}
1 & 0\\
e^{-2 p \phi_{+}(\lambda)} & 1\\
\end{pmatrix}.
\end{equation}
Our aim is to arrive at an RHP with a constant jump on the interval $(-1,1)$ and a ``small'' jump on the contour $\Sigma_1 \cup \Sigma_2 \cup (\varphi_p^{-1}(0),-1) \cup (1, \varphi_p^{-1}(1))$, where $\Sigma_j$'s, $j=1,2$ are constructed according to the following proposition.

\begin{proposition}\label{prop:small_norm}
Let $\mathcal{U}_1$ and $\mathcal{U}_{-1}$ be small disks (see Remark~\ref{remark disks}) centered at $1$ and $-1$, respectively, and define \(\cV=\mathbb{C}\setminus\bigl(\mathcal{U}_1\cup\mathcal{U}_{-1}\bigr).\) Consider the matrices $L$ and $U$ defined in \eqref{eq:factor}. There exist contours $\Sigma_1$ and $\Sigma_2$, and constants $\kappa>0$ and $p_0\in\mathbb{N}$, such that, for all $p\geq p_0$,
\begin{equation}
        U(\lambda) = \begin{pmatrix}
            1 & 0 \\
            0 & 1
        \end{pmatrix} + O(e^{-p \kappa}),  \quad \lambda \in \Sigma_1 \cap \cV; \quad L(\lambda) = \begin{pmatrix}
            1 & 0 \\
            0 & 1
        \end{pmatrix} + O(e^{-p \kappa}), \quad \lambda \in \Sigma_2 \cap \cV.
    \end{equation}
\end{proposition}

\begin{proof}
We focus on $U(\lambda)$, and the proof for $L(\lambda)$ is analogous. Due to the monotonicity of $\Re(\phi(\lambda))$, we have that $\Re(\phi(\lambda)) > 0$ above the interval $(-1,1)$ (see Eq.~\eqref{eq:monotone}). Therefore, we can choose a contour $\Sigma_1$ lying above $(-1,1)$ along which $\Re(\phi(\lambda)) > 0$. Since $\Re(\phi(\lambda))$ is continuous on $\Sigma_1 \cap \cV$, there exists a constant $\kappa > 0$ such that
\[ \min_{\lambda \in \Sigma_1 \cap \cV} \Re(\phi(\lambda))  = \kappa.\]
This implies that there exists $p_0\in \mathbb{N},$ such that for all $p>p_0$ 
\begin{equation}
e^{-2 p \phi(\lambda)} =O(e^{-p\kappa}), \quad \lambda \in \Sigma_1 \cap \cV,
\end{equation}
proving the result.
\end{proof}

\begin{figure}[h]
    \centering
    \begin{tikzpicture}[
        myarrow/.style={
            thick,
            postaction={decorate,
            decoration={markings,mark=at position #1 with {\arrow[scale=1.5]{stealth}}}}
        },
        every matrix/.style={scale=0.6} ]
        \draw[myarrow=0.5] (-3,0) -- (3,0);
        \draw[myarrow=0.5] (-6,0) -- (-3,0);
        \draw[myarrow=0.5] (3,0) -- (6,0);
        \draw[myarrow=0.65] 
            (-3,0) .. controls (0,2) .. (3,0); 
        \draw[myarrow=0.65] 
            (-3,0) .. controls (0,-2) .. (3,0); 

        \matrix (m1) at (-0.5,0.5) [matrix of math nodes, left delimiter=(, right delimiter=), inner sep=0pt, column sep=10pt] {
            0& 1\\ \\[2pt]
            -1 &0\\};

        \matrix (m2) at (4.5,0.5) [matrix of math nodes, left delimiter=(, right delimiter=), inner sep=0pt, column sep=10pt] {
            1 & e^{2 p \phi(\lambda)}  \\ \\[2pt]
           0& 1\\};

         \matrix (m3) at (-4.5,0.5) [matrix of math nodes, left delimiter=(, right delimiter=), inner sep=0pt, column sep=10pt] {
            1 & e^{2 p \phi(\lambda)}  \\ \\[2pt]
           0& 1\\};

         \matrix (m4) at (2.5,1.5) [matrix of math nodes, left delimiter=(, right delimiter=), inner sep=0pt, column sep=10pt] {
            1 & 0  \\ \\[2pt]
           e^{-2 p \phi(\lambda)}& 1\\};

        \matrix (m5) at (2.5,-1.5) [matrix of math nodes, left delimiter=(, right delimiter=), inner sep=0pt, column sep=10pt] {
            1 & 0  \\ \\[2pt]
           e^{-2 p \phi(\lambda)}& 1\\};
           
        \node (eq) at (-3,-0.3) {$-1$};
        \node (eq) at (3,-0.3) {$1$};
        \node (eq) at (-5.5,-0.3) {$\varphi_p^{-1}(0)$};
        \node (eq) at (5.5,-0.3) {$\varphi_p^{-1}(1)$};
        \node (eq) at (1,0.3) {$\Omega_1$};
        \node (eq) at (1,-0.3) {$\Omega_2$};
        \node (eq) at (-1.5,1.4) {$\Sigma_1$};
        \node (eq) at (-1.5,-1.4) {$\Sigma_{2}$};
        
        \draw[dashed] (3,0) circle (0.5);
        \draw[dashed] (-3,0) circle (0.5);
        \node (eq) at (-3.5,-0.7) {$U_{-1}$};
        \node (eq) at (3.5,-0.7) {$U_{1}$};

    \end{tikzpicture}
    \caption{Jump matrices for $\widetilde{T}(\lambda)$.}
    \label{fig Contour T tilde}
\end{figure}
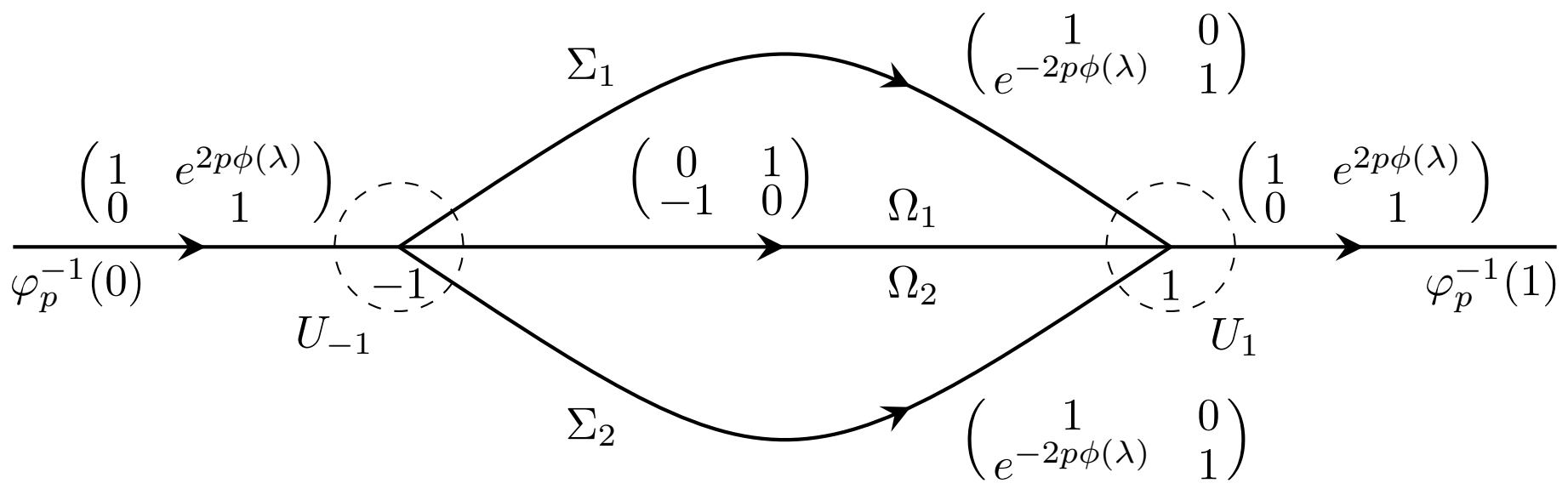

\begin{remark}\label{remark disks}
The contours $\Sigma_1$ and $\Sigma_2$ are upper and lower lenses connecting $-1$ and $1$, with $\Sigma_1$ lying above $(-1,1)$ and $\Sigma_2$ lying below $(-1,1)$, see Fig.~\ref{fig Contour T tilde}. Apart from these conditions, their precise shapes are left unspecified at this stage and will be fixed in Section \ref{e:sec24}. Moreover, the disks are taken to have a common radius $\epsilon>0$, chosen sufficiently small so that they are disjoint; see Fig.~\ref{fig Contour T tilde}. The value of $\epsilon$ will be fixed later so as to ensure the existence of a local map defined in Section \ref{e:sec24}.
\end{remark}
Given $\Sigma_1$, and $\Sigma_2$, we define the function:
\begin{eqnarray}\label{eq:eq243}
\widetilde{T}(\lambda) = \begin{cases}
T(\lambda) \, U^{-1}(\lambda), \quad & \lambda \in \Omega_1 \\
T(\lambda) \, L(\lambda), \quad & \lambda \in \Omega_2 \\
T(\lambda), \quad & \text{elsewhere},
\end{cases}
\end{eqnarray}
where $\Omega_1$ and $\Omega_2$ are the two regions shown in Fig.~\ref{fig Contour T tilde}. Then $\widetilde{T}$ satisfies the following RHP.
\begin{RHP}\label{e:RHP25}
Find a $2 \times 2$ matrix-valued function $\widetilde{T}(\lambda)$ that satisfies the following conditions.
\begin{enumerate}
\item $\widetilde{T} (\lambda)$ is analytic in $\mathbb{C} \setminus \Sigma, \quad \Sigma = \Sigma_1 \cup \Sigma_2 \cup (\varphi_p^{-1}(0),-1) \cup (-1, 1) \cup (1, \varphi_p^{-1}(1))$.
\item $\widetilde{T} (\lambda)$ satisfies the jump relation $\widetilde{T}_{+}(\lambda) = \widetilde{T}_{-}(\lambda) \, J_{\widetilde{T}}(\lambda)$ with the jump matrix given by:
\begin{align*}
J_{\widetilde{T}}(\lambda) = \begin{cases}
\begin{pmatrix}
1 & e^{2 p \phi (\lambda)} \\
0 & 1
\end{pmatrix},  \quad & \lambda \in (\varphi_p^{-1}(0),-1) \cup (1, \varphi_p^{-1}(1)) \\ \\[0.5pt]
U(\lambda), \quad & \lambda \in \Sigma_1 \\ \\[0.5pt]
L(\lambda), \quad & \lambda \in \Sigma_2 \\ \\[0.5pt]
M, \quad & \lambda \in (-1, 1)
\end{cases}
\end{align*}
\item $\widetilde{T} (\lambda)$ has the following asymptotic behavior:
\begin{eqnarray}
\widetilde{T} (\lambda) = I_2 + \mathcal{O}(\lambda^{-1}), \quad \lambda \to \infty.
\end{eqnarray}
\end{enumerate}     
\end{RHP}

\begin{remark}\label{e:remark22}
Equation \eqref{eq:funphi} implies that $\phi(\lambda) < 0$, for $\lambda \in (\varphi_p^{-1}(0),-1) \cup (1, \varphi_p^{-1}(1))$. A similar result to Proposition~\ref{prop:small_norm} can also be applied to $J_{\widetilde{T}}$ on $(\varphi_p^{-1}(0),-1) \cup (1, \varphi_p^{-1}(1))$. To arrive at an approximation to $\widetilde{T}$, we will temporarily ignore the jumps on $(\varphi_p^{-1}(0),-1) \cup (1, \varphi_p^{-1}(1))$, as well as on $\Sigma_1 \cup \Sigma_2$.
\end{remark}


\subsection{Outer parametrix}

In this subsection, we formulate and solve the model problem that yields the outer approximate solution to RHP~\ref{e:RHP25} for $\lambda$ away from $\pm 1$. The following RHP originates from the small norm approximation presented in the previous section (see Remark~\ref{e:remark22}).

\begin{RHP}\label{e:RHP24}
Find a $2 \times 2$ matrix-valued function $\hat{T}(\lambda)$ that satisfies the following conditions. 
\begin{enumerate}
    \item $\hat{T}(\lambda)$ is analytic in $\mathbb{C} \setminus (-1, 1)$.
    \item $\hat{T}(\lambda)$ satisfies the jump relation
    \begin{equation}\label{eq:ejumpThat}
    \hat{T}_{+}(\lambda) = \hat{T}_{-}(\lambda) \begin{pmatrix}
        0 & 1\\
        -1 & 0
    \end{pmatrix}, \quad \lambda \in (-1, 1).
    \end{equation}
    \item $\hat{T}(\lambda) = I_2 + \mathcal{O}(\lambda^{-1}), \quad \lambda \to \infty$.
    \item $\hat{T}(\lambda)=\mathcal{O}(|\lambda+1|^{-1/4}), \quad \lambda \to -1$,\, and \,\, $\hat{T}(\lambda)=\mathcal{O}(|\lambda-1|^{-1/4}), \quad \lambda \to 1$.
\end{enumerate}
\end{RHP}
This Riemann-Hilbert problem can be solved explicitly by diagonalizing the jump matrix and reducing it to two scalar Riemann-Hilbert problems, see \cite[p. 1520]{deift1999strong}. The unique solution is given by
\begin{equation}\label{e:That}
    \hat{T}(\lambda) = \begin{pmatrix}
    1 & 1\\
    i & -i 
\end{pmatrix} \, F(\lambda) \, \begin{pmatrix}
1 & 1\\
i & -i 
\end{pmatrix}^{-1}=\frac{1}{2} \begin{pmatrix}
F_{11}(\lambda) + \frac{1}{F_{11}(\lambda)} & - i F_{11}(\lambda) + \frac{i}{F_{11}(\lambda)} \\
i F_{11}(\lambda) - \frac{i}{F_{11}(\lambda)} & F_{11}(\lambda) + \frac{1}{F_{11}(\lambda)}
\end{pmatrix},
\end{equation}
where 
\begin{equation}
F(\lambda) = \left( F_{11}(\lambda) \right)^{\sigma_3}, \quad F_{11}(\lambda) = \left( \frac{\lambda - 1}{\lambda + 1} \right)^{1/4}.
\end{equation}


\subsection{Local parametrix}\label{e:sec24}

\begin{figure}[ht]
    \centering
    \begin{tikzpicture}[
        myarrow/.style={
            thick,
            postaction={decorate,
            decoration={markings,mark=at position #1 with {\arrow[scale=1.5]{stealth}}}}
        },
        every matrix/.style={scale=0.4} 
    ]
        \draw[myarrow=0.5] (-2.5,-2.5) -- (0,0);
        \draw[myarrow=0.5] (-3.5,0) -- (0,0);
        \draw[myarrow=0.5] (0,0) -- (3,0);
        \draw[myarrow=0.5] (-2.5,2.5) -- (0,0);
        
        \matrix (m1) at (-0.5,2) [matrix of math nodes, left delimiter=(, right delimiter=), inner sep=0pt, column sep=10pt] {
            1 & 0\\ \\[2pt]
            1 & 1 \\
        };
        \matrix (m2) at (-3,0.5) [matrix of math nodes, left delimiter=(, right delimiter=), inner sep=0pt, column sep=10pt] {
            0 & 1\\ \\[2pt]
            -1 & 0 \\
        };
         \matrix (m3) at (-0.5,-2) [matrix of math nodes, left delimiter=(, right delimiter=), inner sep=0pt, column sep=10pt] {
            1 & 0\\ \\[2pt]
            1 & 1 \\
        };        
        \matrix (m4) at (2,0.5) [matrix of math nodes, left delimiter=(, right delimiter=), inner sep=0pt, column sep=10pt] {
            1 & 1  \\ \\[2pt]
           0& 1\\
        };
        \node (eq) at (-1.5,2) {$\gamma_1$};
        \node (eq) at (-1.5,-2) {$\gamma_3$};
        \node (eq) at (0.8,-0.3) {$\gamma_4$};
        \node (eq) at (-1.2,-0.3) {$\gamma_2$};
        \node (eq) at (0,-0.3) {$0$};
        \node (eq) at (1,2) {$I$};
        \node (eq) at (-3,2) {$II$};
        \node (eq) at (-3,-2) {$III$};
        \node (eq) at (1,-2) {$IV$};
    \end{tikzpicture}
    \caption{Jump relations for $\Psi$.}
    \label{fig Contour Psi1(xi)}
\end{figure}
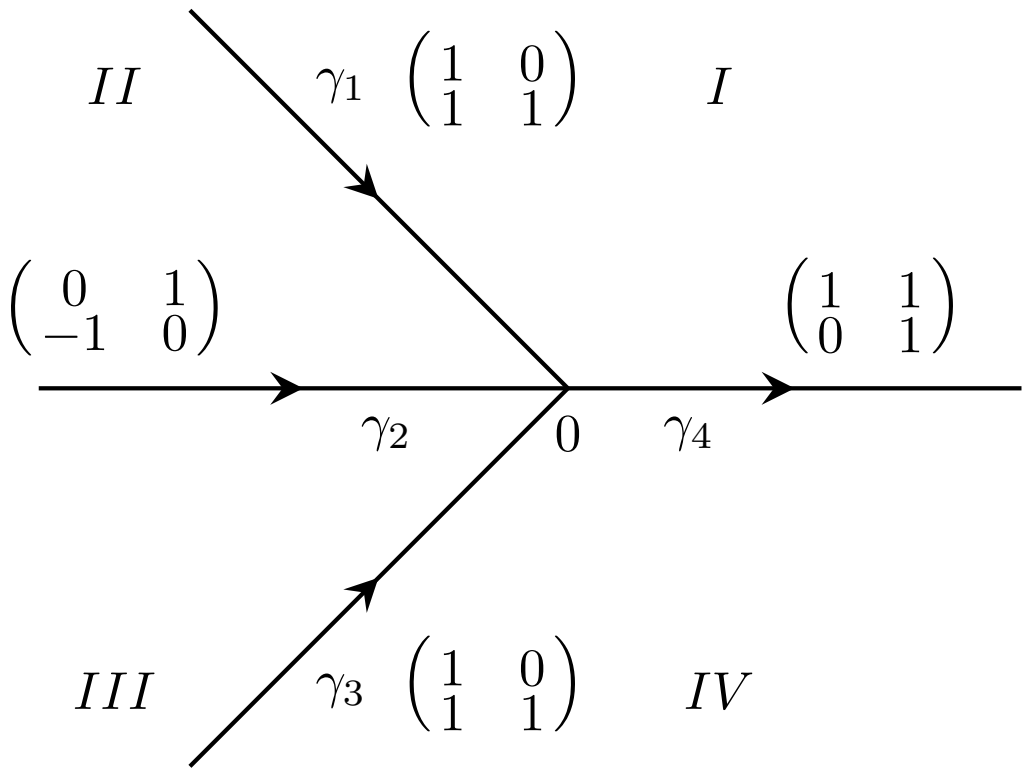

In the previous section, we derived a solution to the RHP for $\hat{T}(\lambda)$ which serves as an outer approximation to $\widetilde{T}(\lambda)$  away from $\lambda = \pm 1$. Near these points, however, this approximation is no longer accurate, as the jump matrices that we ignored are not close to the identity matrix. Therefore, we must construct a solution in a small neighborhood of $\lambda=\pm 1$ to better approximate $\widetilde{T}(\lambda)$ locally. We follow \cite{deift1999strong} to construct these local parametrices. Recall that $\mathcal{U}_1$ and $\mathcal{U}_{-1}$ are two disks of radius $\epsilon$ around $\pm 1$, respectively. We begin with the local parametrix at $\lambda=1$, denoted by $\widetilde{T}_{1}$, which satisfies the following jump relation within $\mathcal{U}_1$:
\begin{align}\label{eq:eq255}
\left( \widetilde{T}_{1}(\lambda) \right)_{+} = \left( \widetilde{T}_{1}(\lambda) \right)_{-} \,J_{\widetilde{T}_{1}}(\lambda), \quad 
J_{\widetilde{T}_{1}}(\lambda) = \begin{cases}
\begin{pmatrix}
1 & e^{2 p \phi (\lambda)} \\
0 & 1
\end{pmatrix},  \quad & \lambda \in (1,1+\epsilon) \\ \\[0.5pt]
\begin{pmatrix}
1 & 0\\
e^{-2 p \phi(\lambda)} & 1\\
\end{pmatrix}, \quad & \lambda \in \Sigma_1 \\ \\[0.5pt]
\begin{pmatrix}
1 & 0\\
e^{-2 p \phi(\lambda)} & 1\\
\end{pmatrix}, \quad & \lambda \in \Sigma_2 \\ \\[0.5pt]
\begin{pmatrix}
0 & 1\\
-1 & 0
\end{pmatrix}, \quad & \lambda \in (-1, 1) 
\end{cases}
\end{align}
We construct the solution $\widetilde{T}_1$ using Airy functions (see \cite{deift1999strong} for further details). We define the function
\begin{align}\label{eq: sol RHP for Psi1}
    \Psi(\zeta)=\begin{cases}
        \begin{pmatrix}
            \mathrm{Ai}(\zeta)&\mathrm{Ai}(\omega^2 \zeta)\\
            \mathrm{Ai}'(\zeta)&\omega^2 \mathrm{Ai}'(\omega^2 \zeta)
        \end{pmatrix} e^{-\frac{\pi i}{6}\sigma_3}, \quad & \zeta \in I,
          \vspace{0.1cm} \\
           \begin{pmatrix}
            \mathrm{Ai}(\zeta)&\mathrm{Ai}(\omega^2 \zeta)\\
            \mathrm{Ai}'(\zeta)&\omega^2 \mathrm{Ai}'(\omega^2 \zeta)
        \end{pmatrix} e^{-\frac{\pi i}{6}\sigma_3}
        \begin{pmatrix}
            1 &0\\
            -1& 1
        \end{pmatrix}, \quad & \zeta \in II,
         \vspace{0.1cm}\\
         \begin{pmatrix}
            \mathrm{Ai}(\zeta)&-\omega^2 \mathrm{Ai}(\omega \zeta)\\
            \mathrm{Ai}'(\zeta)&- \mathrm{Ai}'(\omega \zeta)
        \end{pmatrix} e^{-\frac{\pi i}{6}\sigma_3}
        \begin{pmatrix}
            1 &0\\
            1& 1
        \end{pmatrix}, \quad & \zeta \in III,
         \vspace{0.1cm}\\
         \begin{pmatrix}
            \mathrm{Ai}(\zeta)&-\omega^2 \mathrm{Ai}(\omega \zeta)\\
            \mathrm{Ai}'(\zeta)&-\mathrm{Ai}'(\omega \zeta)
        \end{pmatrix} e^{-\frac{\pi i}{6}\sigma_3}, \quad & \zeta \in IV,
        \end{cases}
\end{align}
where $\omega=e^{\frac{2\pi i}{3}}$. In \cite{deift1999strong}, the authors show that $\Psi$ satisfies the jumps shown in Fig.~\ref{fig Contour Psi1(xi)}.

\begin{proposition}\label{eq:pro26}
The function $\widetilde{T}_1$ given by the expression
\begin{equation}\label{eq: T tilde +1}
\widetilde{T}_1(\lambda) = Q(\lambda) \Psi(\chi_p(\lambda)) e^{-p \, \phi(\lambda) \sigma_3}
\end{equation}
satisfies jump relation \eqref{eq:eq255}. Here, $\Psi$ is as in \eqref{eq: sol RHP for Psi1}, $Q$ is an arbitrary invertible analytic function inside $\mathcal{U}_{1}$ to be determined, and the map $\chi_p$ (recall also \eqref{eq:mapchi}) is defined by \eqref{eq:eqmao}.
\end{proposition}
\begin{proof}
Analysis of the local behavior of $\phi(\lambda)$ near $\lambda=1$ yields:
\begin{align}\label{eq local}
\phi(\lambda) & = \frac{1}{2} \left( \frac{\beta_p - \alpha_p}{2} \right)^2 \int_{1}^{\lambda} (\zeta + 1)^{1/2} (\zeta - 1)^{1/2} \left( \frac{b_{\infty}}{R(0) s(\zeta)} + \frac{a_{\infty} p^{\nu}}{R(1) \big( s(\zeta) -1\big)} \right) \, d \zeta \nonumber \\ 
 & = \int_{\lambda}^{1} (\zeta - 1)^{1/2} h_p(\zeta) \, d \zeta \nonumber \\
& = \frac{2}{3} h_p(1) (\lambda -1)^{3/2} + \mathcal{O}(\lambda-1)^{5/2} ,
\end{align}
where $h_p$ is analytic in a neighborhood of $\lambda=1$ and is given by
\begin{equation}\label{eq. def function h}
  h_p(\zeta) = \frac{1}{2} \left( \frac{\beta_p - \alpha_p}{2} \right)^2 (\zeta + 1)^{1/2}  \left( \frac{b_{\infty}}{R(0) s(\zeta)} + \frac{a_{\infty} p^{\nu}}{R(1) \big( s(\zeta) -1\big)} \right).  
\end{equation}
More generally, for $\nu>0$, the expression in \eqref{eq. def function h} yields the following for $p \to \infty$:
\begin{align*}
h_p(\zeta)=h_\infty(\zeta)+O(p^{-\nu}), \quad h_{\infty}(\zeta)=-\frac{2(1+b_{\infty})(\zeta+1)^{1/2}}{2 \sqrt{1+b_{\infty}} \zeta +2 + b_{\infty}}+O(p^{-\nu}),
\end{align*}
uniformly in a sufficiently small fixed neighborhood of $\zeta=1$. Since $h_\infty(1)\neq 0$, this provides uniform control of the local conformal mapping $\chi_p$ for sufficiently large $p$. In particular, the value $h_p(1)$ has the following asymptotic expansion in $p$:
\begin{equation}\label{eq:consh1}
h_p(1) = \frac{-2\sqrt{2} \, (1+b_{\infty})}{2 + b_{\infty} + 2\sqrt{1+b_{\infty}}} + \mathcal{O} \left( p^{-\nu} \right), \quad p \to \infty.
\end{equation}
From the local expansion of $\phi$ at $\lambda=1$, the function
\begin{equation}\label{eq:eqmao}
\chi_p(\lambda)
=
p^{\frac{2}{3}}
\left(-\frac{3}{2}\phi(\lambda)\right)^{2/3}
\end{equation}
is holomorphic in a neighborhood of $\lambda=1$ and satisfies $\chi_p(1)=0
$, $\chi_p'(1) = p^{2/3}(-h_p(1))^{2/3}\neq0$. Hence, $\chi_p$ is biholomorphic in a sufficiently small neighborhood of $\lambda=1$.

\begin{remark}\label{remark epsilon}
    \begin{enumerate}
        \item The radius $\epsilon$ of the disk $\mathcal{U}_1$ introduced in Proposition~\ref{prop:small_norm} was left unspecified. We now choose $\epsilon=\epsilon_0$ so that $\mathcal{U}_1$ is contained in the neighborhood on which the local map $\chi_p$ is defined. The choice of the disk determines the regions in which the solution of the RHP is approximated either by Airy functions or by simpler exponential functions. In the overlapping regions, the large-argument asymptotics of the Airy functions exhibit the exponential behavior expected from the latter approximation, illustrating the consistency between the two representations. Consequently, any disk centered at $\lambda=1$ with a fixed radius $0<\epsilon<\epsilon_0$ can also be used without affecting the validity of the asymptotic analysis.
        
        \item As $p\rightarrow\infty$, the map $\chi_p$ sends a sufficiently small neighborhood of $\lambda=1$ onto a comparatively large neighborhood of $\xi=0$, while mapping the portion of the interval $(-1,1)$ lying inside $\mathcal{U}_1$ onto a segment of the negative real axis in the $\xi$-plane.
        
        \item Here, we also choose the contours $\Sigma_1$ and $\Sigma_2$ so that, under $\chi_p$, the portions lying inside $\mathcal{U}_1$ are mapped precisely onto the rays $\gamma_1$ and $\gamma_3$, as illustrated in Fig.~\ref{fig Contour Psi1(xi)}.
    \end{enumerate}
\end{remark}

Now, consider the jump relations satisfied by $\widetilde{T}_1 $ in the $\xi$-plane:
\begin{align*}
\left( \widetilde{T}_{1} \left( \chi_p^{-1}(\xi)\right) \right)_{+} = \left( \widetilde{T}_{1} \left( \chi_p^{-1}(\xi)\right) \right)_{-} \,
\begin{cases}
\begin{pmatrix}
1 & e^{- \frac{4}{3} \xi^{3/2}} \\
0 & 1
\end{pmatrix},  \quad & \xi \in (0,+\infty) \\ \\[0.5pt]
\begin{pmatrix}
1 & 0\\
e^{\frac{4}{3} \xi^{3/2}} & 1\\
\end{pmatrix}, \quad & \xi \in \gamma_1\\ \\[0.5pt]
\begin{pmatrix}
1 & 0\\
e^{\frac{4}{3} \xi^{3/2}} & 1\\
\end{pmatrix}, \quad & \xi \in \gamma_3\\ \\[0.5pt]
\begin{pmatrix}
0 & 1\\
-1 & 0
\end{pmatrix}, \quad & \xi \in (-\infty, 0) ,
\end{cases}
\end{align*}

and define
\begin{equation}\label{eq:psifun}
\Psi(\xi) = \widetilde{T}_1(\chi_p^{-1}(\xi)) \, e^{-\frac{2}{3}\xi^{3/2}\sigma_3}.
\end{equation}
A direct calculation shows that $\Psi$ satisfies the expected jump conditions given on Fig.~\ref{fig Contour Psi1(xi)}. Then, using relation \eqref{eq:psifun} and the definition of the conformal (and hence invertible) map $\chi_p(\lambda)$, we define $\widetilde{T}_1(\lambda)$ back in the $\lambda$-plane, as desired (modulo the analytic function $Q(\lambda)$ which will be determined in Section \ref{eq:sec26}, and does not affect the jump relation).
\end{proof}

\begin{remark}
When $\nu=0$, the function $\phi(\lambda)$ has the following behavior  as $\lambda\to1$:
\[\phi(\lambda) = \frac{2}{3} h(1) (\lambda -1)^{3/2}+ \mathcal{O}(\lambda-1)^{5/2}.\]
Here $h(1)$ is the fixed constant:
\begin{equation}\label{eq:h1}
h(1) = \frac{1}{\sqrt{2}} \left( \frac{\beta - \alpha}{2} \right)^2 \left( \frac{b_{\infty}}{R(0) s(1)} + \frac{a_{\infty}}{R(1) \big( s(1) -1\big)} \right),
\end{equation}
independent of $p$ (see~\eqref{eq. def function h}).
\end{remark}

Next, we construct the local parametrix at $\lambda=-1$, denoted by $\widetilde{T}_{-1}$, which satisfies the following jump relation within $\mathcal{U}_{-1}$:
\begin{align}\label{eq:eqjump}
\left( \widetilde{T}_{-1}(\lambda) \right)_{+} = \left( \widetilde{T}_{-1}(\lambda) \right)_{-} \,J_{\widetilde{T}_{-1}}(\lambda), \quad 
J_{\widetilde{T}_{-1}}(\lambda) = \begin{cases}
\begin{pmatrix}
1 & e^{2 p \phi (\lambda)} \\
0 & 1
\end{pmatrix},  \quad & \lambda \in (-1-\epsilon,-1) \\ \\[0.5pt]
\begin{pmatrix}
1 & 0\\
e^{-2 p \phi(\lambda)} & 1\\
\end{pmatrix}, \quad & \lambda \in \Sigma_1 \\ \\[0.5pt]
\begin{pmatrix}
1 & 0\\
e^{-2 p \phi(\lambda)} & 1\\
\end{pmatrix}, \quad & \lambda \in \Sigma_2 \\ \\[0.5pt]
\begin{pmatrix}
0 & 1\\
-1 & 0
\end{pmatrix}, \quad & \lambda \in (-1, 1) 
\end{cases}
\end{align}

\begin{proposition}
The function $\widetilde{T}_{-1}$ given by the following expression
\begin{equation}\label{eq: T tilde -1}
\widetilde{T}_{-1}(\lambda) = \widetilde{Q}(\lambda) \, \sigma_3 \, \Psi \left( - \chi_p(\lambda) \right) \, \sigma_3 \, e^{-p \, \phi(\lambda) \sigma_3}
\end{equation}
satisfies jump relation \eqref{eq:eqjump}. Here, $\widetilde{Q}(\lambda)$ is an arbitrary invertible analytic function inside $\mathcal{U}_{-1}$ to be determined, $\Psi$ and $\chi_p$ are as in Proposition~\ref{eq:pro26}.
\end{proposition}

\begin{proof}
For brevity and for the purposes of the proof, we introduce the notation
\begin{equation}\label{eq:psitil}
\widetilde{\Psi}\left(\chi_p(\lambda) \right) = \sigma_3 \, \Psi \left( - \chi_p(\lambda) \right) \, \sigma_3.
\end{equation}
One can show that if $\Psi$ satisfies the jump conditions shown in Fig.~\ref{fig Contour Psi1(xi)}, then $\widetilde{\Psi}$ defined in equation \eqref{eq:psitil} satisfies the jump conditions given in Fig.~\ref{fig Contour Psi_{-1}(xi)}.

Therefore, the relation
\begin{equation}
\widetilde{\Psi}(\xi) = \widetilde{T}_{-1}(\chi_p^{-1}(\xi)) \, e^{-\frac{2}{3}\xi^{3/2}\sigma_3}
\end{equation}
defines $\widetilde{T}_{-1}(\lambda)$ back in the $\lambda$-plane, as desired (modulo the analytic function $\widetilde{Q}(\lambda)$ which  will be determined in Section \ref{eq:sec26}, and does not affect the jump relation).
\end{proof}

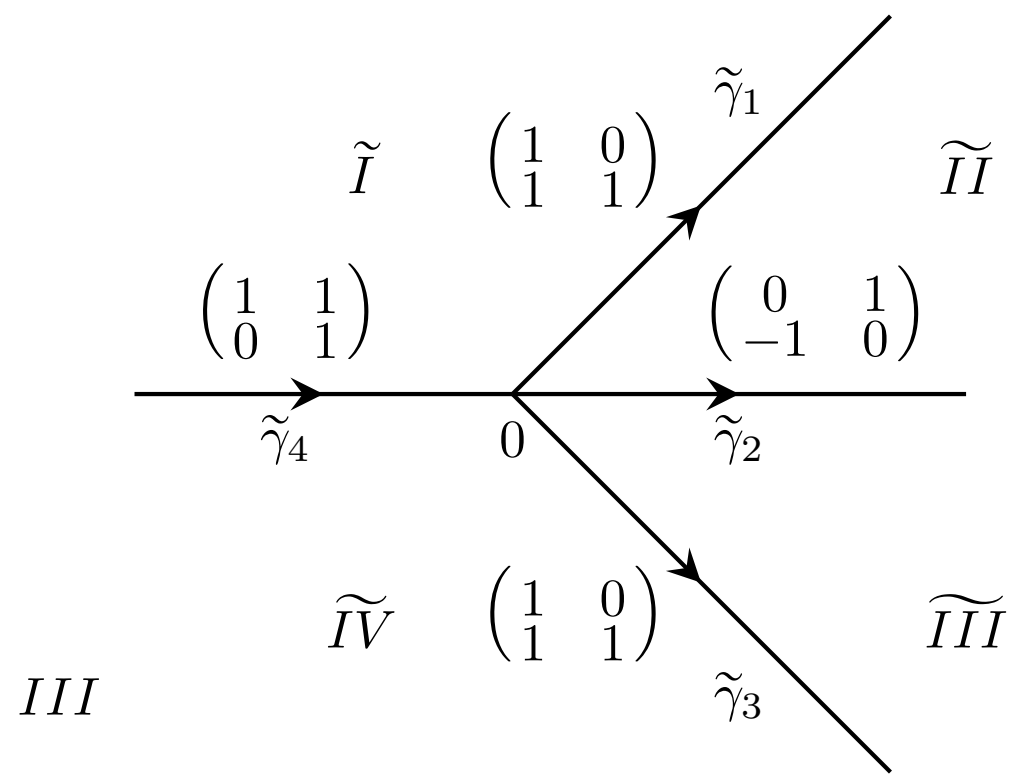
\begin{figure}[h!]
    \centering
    \begin{tikzpicture}[
        myarrow/.style={
            thick,
            postaction={decorate,
            decoration={markings,mark=at position #1 with {\arrow[scale=1.5]{stealth}}}}
        },
        every matrix/.style={scale=0.4} 
    ]
        \draw[myarrow=0.5]  (0,0) -- (2.5,2.5);
        \draw[myarrow=0.5] (-2.5,0) -- (0,0);
        \draw[myarrow=0.5] (0,0) -- (3,0);
        \draw[myarrow=0.5]  (0,0) -- (2.5,-2.5);
        
        \matrix (m1) at (0.4,1.5) [matrix of math nodes, left delimiter=(, right delimiter=), inner sep=0pt, column sep=10pt] {
            1 & 0\\ \\[2pt]
            1 & 1 \\
        };
        \matrix (m2) at (-1.5,0.5) [matrix of math nodes, left delimiter=(, right delimiter=), inner sep=0pt, column sep=10pt] {
         1 & 1  \\ \\[2pt]
           0& 1\\
        };
         \matrix (m3) at (0.4,-1.5) [matrix of math nodes, left delimiter=(, right delimiter=), inner sep=0pt, column sep=10pt] {
            1 & 0\\ \\[2pt]
            1 & 1 \\
        };        
        \matrix (m4) at (2,0.5) [matrix of math nodes, left delimiter=(, right delimiter=), inner sep=0pt, column sep=10pt] {
           0 & 1\\ \\[2pt]
            -1 & 0 \\
        };
        \node (eq) at (1.5,2) {$\widetilde{\gamma}_1$};
        \node (eq) at (1.5,-2) {$\widetilde{\gamma}_3$};
        \node (eq) at (1.5,-0.3) {$\widetilde{\gamma}_2$};
        \node (eq) at (-1.5,-0.3) {$\widetilde{\gamma}_4$};
        \node (eq) at (0,-0.3) {$0$};
        \node (eq) at (3,1.5) {$\widetilde{II}$};
        \node (eq) at (3,-1.5) {$\widetilde{III}$};
        \node (eq) at (-1,1.5) {$\widetilde{I}$};
        \node (eq) at (-3,-2) {$III$};
        \node (eq) at (-1,-1.5) {$\widetilde{IV}$};
    \end{tikzpicture}
    \caption{Jump matrices for $\widetilde{\Psi}$.}
    \label{fig Contour Psi_{-1}(xi)}
\end{figure}


\subsection{The error function and the global parametrix}\label{eq:sec26}

In this subsection, we define the global approximation to $\widetilde{T}$, and introduce the error function $E$. The global approximation is defined separately within each of the disks $\mathcal{U}_{1}$ and $\mathcal{U}_{-1}$, and in the exterior region outside both disks:
\begin{gather}\label{eq:Tapprx}
\widetilde{T}_{\mathrm{approx}}(\lambda) = 
\begin{cases}
\widetilde{T}_{1}(\lambda), & \quad \lambda \in \mathcal{U}_{1}, \\
\widetilde{T}_{-1}(\lambda), & \quad \lambda \in \mathcal{U}_{-1} \\
\hat{T}(\lambda), & \quad \lambda \in \mathbb{C} \setminus \{\mathcal{U}_{1} \cup \mathcal{U}_{-1}\}
\end{cases}
\end{gather}
where $\hat{T}$, $\widetilde{T}_{1}$ and $\widetilde{T}_{-1}$ are known explicitly (see sections~\ref{e:sec23} and \ref{e:sec24}). 

The quantity $\widetilde{T}_{\mathrm{approx}}$ is a piecewise analytic matrix-valued function defined in the entire complex plane, satisfying jump relationships across most of the contours shown in Figure~\ref{fig Contour T tilde}.  More precisely:
\begin{itemize}
\item Inside the disks centered at $\pm 1$, $\widetilde{T}_{\mathrm{approx}}(\lambda)$ is piecewise analytic with jumps across the contours $\Sigma_{1}$, $\Sigma_{2}$, and the real axis, and satisfies the same jump relations there as $\widetilde{T}(\lambda)$.
\item Outside the disks, $\widetilde{T}_{\mathrm{approx}}(\lambda)$ is analytic in $\mathbb{C} \setminus [-1,1]$, and satisfies the same jump relations across $[-1,1]$ as $\widetilde{T}$.  
\item Outside the disks, $\widetilde{T}_{\mathrm{approx}}(\lambda)$ is analytic (and has no jumps) across the contours $\Sigma_{1}$, $\Sigma_{2}$, and the rest of the real axis.  
\end{itemize}

We then define $E(\lambda)$:
\begin{gather}\label{eq:def E}
E(\lambda) = \widetilde{T}(\lambda) \, \left( \widetilde{T}_{\mathrm{approx}}(\lambda) \right)^{-1}.
\end{gather}

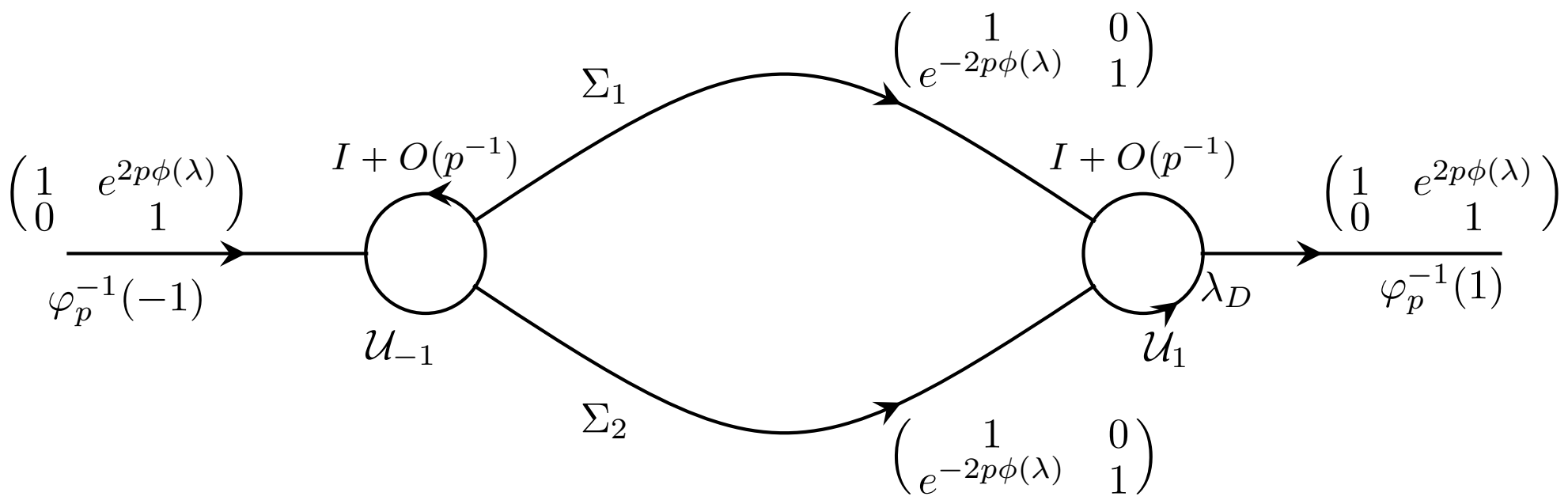
\begin{figure}[ht]
    \centering
    \begin{tikzpicture}[
        myarrow/.style={
            thick,
            postaction={decorate,
            decoration={markings,mark=at position #1 with {\arrow[scale=1.5]{stealth}}}}
        },
        every matrix/.style={scale=0.6} ]
        
        \draw[myarrow=0.5] (-6,0) -- (-3,0);
        \draw[myarrow=0.5] (3,0) -- (6,0);
        
        \draw[myarrow=0.65] 
            (-3,0) .. controls (0,2) .. (3,0); 
        \draw[myarrow=0.65] 
            (-3,0) .. controls (0,-2) .. (3,0); 
      
        \fill[white] (-3,0) circle [radius=0.48];
        \fill[white] (3,0) circle [radius=0.48];

        \draw[thick, black, postaction={decorate,decoration={markings,mark=at position 0.25 with {\arrow[scale=1.5]{stealth}}}}] 
    (-3,0) circle [radius=0.5];
\draw[thick, black, postaction={decorate,decoration={markings,mark=at position 0.85 with {\arrow[scale=1.5]{stealth}}}}] 
    (3,0) circle [radius=0.5];

        \matrix (m2) at (5.5,0.5) [matrix of math nodes, left delimiter=(, right delimiter=), inner sep=0pt, column sep=10pt] {
            1 & e^{2 p \phi(\lambda)}  \\ \\[2pt]
           0& 1\\};

         \matrix (m3) at (-5.5,0.5) [matrix of math nodes, left delimiter=(, right delimiter=), inner sep=0pt, column sep=10pt] {
            1 & e^{2 p \phi(\lambda)}  \\ \\[2pt]
           0& 1\\};

         \matrix (m4) at (2,1.7) [matrix of math nodes, left delimiter=(, right delimiter=), inner sep=0pt, column sep=10pt] {
            1 & 0  \\ \\[2pt]
           e^{-2 p \phi(\lambda)}& 1\\};

        \matrix (m5) at (2,-1.7) [matrix of math nodes, left delimiter=(, right delimiter=), inner sep=0pt, column sep=10pt] {
            1 & 0  \\ \\[2pt]
           e^{-2 p \phi(\lambda)}& 1\\};
           
        \node (eq) at (-3.2,-0.8) {$\mathcal{U}_{-1}$};
        \node (eq) at (3.2,-0.8) {$\mathcal{U}_{1}$};
         \node (eq) at (3.7,-0.3) {$\lambda_{D}$};
        
        \node (eq) at (-5.5,-0.4) {$\varphi_p^{-1}(-1)$};
        \node (eq) at (5.5,-0.3) {$\varphi_p^{-1}(1)$};
        \node (eq) at (-1.5,1.4) {$\Sigma_1$};
        \node (eq) at (-1.5,-1.4) {$\Sigma_{2}$};
        \node (eq2) at (-3,0.8) [black] {\small{$I+O(p^{-1})$}};
        \node (eq2) at (3,0.8) [black] {\small{$I+O(p^{-1})$}};

    \end{tikzpicture}
    \caption{Jump matrices for $E$.}
    \label{fig Contour E}
\end{figure}

It is straightforward to verify that $E(\lambda)$ is analytic for $\lambda$ in every contour where $\widetilde{T}$ and $\widetilde{T}_{\mathrm{approx}}$ satisfy the same jump relation.  For example, $E(\lambda)$ satisfies $E_{+}(\lambda) = E_{-}(\lambda)$ for all $\lambda \in (-1,1)$, and so $E(\lambda)$ is analytic in $\mathbb{C} \setminus \Sigma_{E}$, where $\Sigma_{E}$ is the collection of contours shown in Fig.~\ref{fig Contour E}.

It is customary to state a new Riemann-Hilbert problem for $E$:
\begin{enumerate}
\item $E(\lambda) $ is analytic for $\lambda \in \mathbb{C} \setminus \Sigma_{E}$.
\item $E(\lambda) = I + \mathcal{O} \left( \frac{1}{\lambda} \right), \ \ \lambda \to \infty$.
\item The boundary values of $E$ satisfy 
\begin{align}
E_{+}(\lambda) = E_{-} (\lambda) J_{E}(\lambda), \ \ \lambda \in \Sigma_{E}.
\end{align}
\end{enumerate}

Moreover, standard calculations now show (see, for example, \cite[Section 7.2]{deift1999strong}) that $J_{E}(\lambda)$  is analytic at any point $\lambda$ on the interior of any of the contours comprising $\Sigma_{E}$, and 
\begin{align}
\left| J_{E}(\lambda) - I \right| \le C_{0} e^{ - C_{1} p}
\end{align}
for all $\lambda$ in $\Sigma_{1}$, $\Sigma_{2}$, and the two real intervals shown in Fig.~\ref{fig Contour E}.

We have not yet defined $Q(\lambda)$ and $\widetilde{Q}(\lambda)$. These quantities were introduced in~\eqref{eq: T tilde +1} and~\eqref{eq: T tilde -1} as arbitrary analytic matrices defined in $\mathcal{U}_{1}$ and $\mathcal{U}_{-1}$ respectively.  We define them as follows:
\begin{equation}\label{eq:eqQ}
Q(\lambda) = \hat{T}(\lambda) \left[\frac{e^{\frac{\pi i}{12}}}{2\sqrt{\pi}} \left( \chi_p(\lambda) \right)^{-\frac{\sigma_3}{4}}
    \begin{pmatrix}
1 & 1\\
-1 & 1
\end{pmatrix}e^{-\frac{\pi i}{4} \sigma_3}  
  \right]^{-1} \ , \ \ \lambda \in \mathcal{U}_{1},
\end{equation}
and
\begin{eqnarray}
\widetilde{Q}(\lambda) = \hat{T}(\lambda) \sigma_3 \left[\frac{e^{\frac{\pi i}{12}}}{2\sqrt{\pi}} \left( - \chi_p(\lambda) \right)^{-\frac{\sigma_3}{4}}
    \begin{pmatrix}
1 & 1\\
-1 & 1
\end{pmatrix}e^{-\frac{\pi i}{4} \sigma_3}  
  \right]^{-1} \sigma_3\ , \ \ \lambda \in \mathcal{U}_{-1}.
\end{eqnarray}

\begin{proposition}\label{prop Q, Qtilde}
The matrix valued function $Q$ defined in~\eqref{eq:eqQ} is analytic inside $\mathcal{U}_{1}$. Moreover, $Q$ is such that 
\begin{align}
J_{E}(\lambda)  = I + \mathcal{O}\left( p^{-1} \right), \ \ \mbox{uniformly for } \lambda \in \partial \mathcal{U}_{1}.
\end{align}
Similarly, the matrix valued function $\widetilde{Q}$ is an analytic function inside $\mathcal{U}_{-1}$.  Moreover, $\widetilde{Q}$ is such that 
\begin{align}
J_{E}(\lambda)  = I + \mathcal{O}\left( p^{-1} \right), \ \ \mbox{uniformly for } \lambda \in \partial \mathcal{U}_{-1}.
\end{align}
\end{proposition}

\begin{proof}
We present the proof for $Q$. The proof for $\widetilde{Q}$ follows by an analogous argument, which we omit.

The jump condition of $E$ on $\partial\mathcal{U}_{1}$ is given by:
\begin{align*}
J_{E}(\lambda) 
= \hat{T}(\lambda) \, \widetilde{T}_{1}^{-1} (\lambda) 
= \hat{T}(\lambda) \, \left[\Psi(\chi_p(\lambda))\,e^{\frac{2}{3} \left( \chi_p (\lambda) \right)^{3/2} \sigma_3}\right]^{-1}\,Q^{-1}(\lambda).
\end{align*}
Using the definition~\eqref{eq:eqQ} of $Q$, and rearranging terms, we have 
\begin{equation}\label{eq. J_E}
  J_E(\lambda)=   \hat{T}(\lambda)\left\{2\sqrt{\pi}\,e^{-\frac{\pi i}{12}} e^{\frac{\pi i}{4}\sigma_3}
   \begin{pmatrix}
1 & -1 \\
1 & 1
\end{pmatrix}
\chi_p(\lambda)^{\sigma_3/4}
\Psi(\chi_p(\lambda))\,
  e^{\frac{2}{3}\chi_p(\lambda)^{3/2} \sigma_3}
  \right\}^{-1} \hat{T}^{-1}(\lambda).
\end{equation}
Next, we utilize the asymptotic behavior of the function $\Psi$ defined in \eqref{eq: sol RHP for Psi1}, as presented in \cite[p. 1526]{deift1999strong} 
\begin{gather}\label{eq:asympsi}
\Psi(\zeta) e^{u \, \sigma_3} = \frac{e^{\frac{\pi i}{12}}}{2 \sqrt{\pi}} \begin{pmatrix}
\zeta^{-1/4} & 0 \\
0 & \zeta^{1/4} 
\end{pmatrix} \sum_{k=0}^{\infty} \begin{pmatrix}
    (-1)^k s_k & s_k \\
    - (-1)^k t_k & t_k 
\end{pmatrix} e^{\frac{- \pi i}{4} \sigma_3} u^{-k}, 
\end{gather}
as $\zeta \to \infty$. Here,
\[s_0 = t_0 = 1, \quad s_k = \frac{\Gamma \left( 3k + \frac{1}{2} \right)}{54^k k! \,\Gamma \left( k + \frac{1}{2} \right)}, \quad t_k = -\frac{6k+1}{6k-1} s_k\]
and $u=\frac{2}{3} \zeta^{3/2}$. 

Note that this is an asymptotic series. So, for example, using only the first term:
\begin{equation}\label{eq.Psi asymtotics}
\Psi(\xi)e^{u\sigma_3}=\frac{e^{\frac{\pi i}{12}}}{2\sqrt{\pi}}  \xi^{-\frac{\sigma_3}{4}}\left[\begin{pmatrix}
    1&1\\
    -1&1
\end{pmatrix}
e^{\frac{-\pi i}{4}\sigma_3}+\mathcal{O}\left(u^{-1}\right),
\right]
\end{equation}
where the notation $\mathcal{O}\left(u^{-1}\right)$ is uniform in this sense: There is $R>0$ and $C>0$ so that for all $|u|>R$,
\[
\left|\mathcal{O}\left(u^{-1}\right)\right|\le \frac{C}{u}.
\]
Now using~\eqref{eq. J_E},~\eqref{eq.Psi asymtotics} and the fact that for $\lambda\in \partial\mathcal{U}_1$, $\chi_p(\lambda)=\mathcal{O}(p)$, we have 
\begin{equation}
    J_E(\lambda)=   \hat{T}(\lambda)
    \left\{I+\mathcal{O}\left(\chi_p(\lambda)^{-3/2}\right)
  \right\}^{-1} \hat{T}^{-1}(\lambda)=I+\mathcal{O}\left(p^{-1}\right).
\end{equation}
To see that $Q$ is analytic inside $\mathcal{U}_1$, we just need to show that it is analytic on the interval $(-1,1)$. Indeed, note that for $\lambda \in (-1,1)$ its boundary values are given by:
\begin{subequations}
\begin{eqnarray}
Q_{+}(\lambda)=2\,\sqrt{\pi}\,e^{-\frac{\pi i}{12}} \,
 \hat{T}_{+}(\lambda)\,e^{\frac{\pi i}{4} \sigma_3}\begin{pmatrix}
1 & 1\\
-1 & 1
\end{pmatrix}^{-1} \left( \chi_{p +} (\lambda) \right)^{\frac{\sigma_3}{4}}
\end{eqnarray}
\begin{eqnarray}
Q_{-}^{-1}(\lambda)=\frac{e^{\frac{\pi i}{12}}}{2 \sqrt{\pi}}\, \left( \chi_{p -} (\lambda) \right)^{\frac{-\sigma_3}{4}} \, \begin{pmatrix}
1 & 1\\
-1 & 1
\end{pmatrix}\,e^{\frac{-\pi i}{4} \sigma_3} \,
 \hat{T}_-^{-1}(\lambda).
\end{eqnarray}
\end{subequations}
So that we have for $\lambda\in (-1,1)\cap\mathcal{U}_1$ 
\[
Q_{-}^{-1}Q_+=\left(\chi_p(\lambda)^{\frac{-\sigma_3}{4}}\right)_{-}
\begin{pmatrix}
1 & 1\\
-1 & 1
\end{pmatrix}
e^{\frac{-\pi i}{4}\sigma_3}T_{-}^{-1}T_{+} e^{\frac{\pi i}{4}\sigma_3}
\begin{pmatrix}
1 & 1\\
-1 & 1
\end{pmatrix}^{-1}
\left(\chi_p(\lambda)^{\frac{-\sigma_3}{4}}\right)_{+}.
\]
Now, using the jump relation of $\hat{T}$ on $(-1,1)$ (cf. \eqref{eq:ejumpThat}), and the jump of $(\chi_p(\lambda))^{\frac{-\sigma_3}{4}}$ for $\lambda\in\mathcal{U}_1\cap(-1,1)$:
\[\left( \chi_{p +} (\lambda) \right)^{\frac{\sigma_3}{4}} = \left( \chi_{p -} (\lambda) \right)^{\frac{\sigma_3}{4}} \, \begin{pmatrix}
i & 0\\
0 & -i
\end{pmatrix},\]
one sees that 
\begin{equation}
  Q_{-}=Q_{+},\quad\text{for}\quad   \lambda\in \mathcal{U}_1\cap(-1,1).
\end{equation}
Using Morera's Theorem, $Q$ is analytic in $\mathcal{U}_1$, with the possible exception of an isolated singularity at $\lambda=1$. However, from its definition the entries of $Q$ are bounded by $C/\sqrt{\lambda}$ as $\lambda\rightarrow 0$, which precludes an isolated singularity, and hence $Q$ is analytic at $\lambda=0$ as well.
\end{proof}

Summarizing, the quantity $E(\lambda)$ satisfies a \textit{small norm Riemann-Hilbert Problem}, because 
\begin{align}
\left| J_{E}(\lambda) - I  \right| &\le C_{0}e^{- C_{1} p},  \ \mbox{ for all } \lambda \in \Sigma_{E} \setminus \left( \partial \mathcal{U}_{1} \  \cup  \partial \mathcal{U}_{-1} \right), \\
\left| J_{E}(\lambda) - I  \right| &\le \frac{C_{2}}{p}, \ \mbox{ for all } \lambda \in \partial \mathcal{U}_{1} \cup \partial \mathcal{U}_{-1}.
\end{align}

Once one is in the small-norm setting, standard techniques from the analysis of Riemann-Hilbert Problems provide an independent existence, uniqueness, and asymptotic description of $E(\lambda)$.  See, for example, the statement and proof of Theorem 7.10 on p. 1532 of \cite{deift1999strong}.  In the present setting, we conclude that $E(\lambda)$ exists, and satisfies
\begin{align}\label{eq. Error matrix}
E(\lambda) = I + \mathcal{O} \left( \frac{1}{p ( 1 + | \lambda | ) } \right)\ \ \mbox{ for all } \lambda \in \mathbb{C}.
\end{align}
\begin{remark}\label{remark error matrix}
First, note that the implicit constant in $\mathcal{O}\left( \frac{1}{p ( 1 + | \lambda | ) }\right)$ is independent of $\lambda$, and only depends on the contour $\Sigma_{E}$ and the radius of the disks $\mathcal{U}_{1}$ and $\mathcal{U}_{2}$ in Figure~\ref{fig Contour E}.  Second, note that $E(\lambda)$ is analytic for $\lambda \in \mathcal{U}_{1}$, and satisfies the above bound uniformly, which can be used in a standard way to obtain asymptotic estimates for derivatives of $E(\lambda)$ via Cauchy's theorem.
\end{remark}


\subsection[Solution on the lambda-plane and Jacobi polynomials]{
Solution on the $\lambda$-plane and Jacobi polynomials}
\label{eq:poly}

We now use the global approximation to $\widetilde{T}$ given in  \eqref{eq:Tapprx} to obtain the explicit solution to RHP~\ref{RHP Y} in the $\lambda$-plane. This is done by translating RHP~\ref{RHP Y} to the variable $\lambda$ and unraveling the sequence of transformations $Y\rightarrow T\rightarrow \widetilde{T}$. First, define
\begin{equation*}
    \widetilde{Y}(\lambda)=\left(\frac{\beta_p-\alpha_p}{2}\right)^{-p\sigma_3}Y(\varphi_p(\lambda)),
\end{equation*}
where $\varphi_p$ is the map defined in \eqref{eq:eqchangeofvari zlambda}. One can check that $\widetilde{Y}$ satisfies the following RHP, which is the natural translation of RHP \ref{RHP Y} to the variable $\lambda$, so that the entries of $\widetilde{Y}$ are related with the monic Jacobi polynomials $\mathbf{P}_p(\lambda):=\left(\frac{\beta_p-\alpha_p}{2}\right)^{-p}\pi_p(\varphi_p(\lambda))$.

\begin{RHP}\label{RHP Y tilde}
Find a $2 \times 2$ matrix valued function $\widetilde{Y} = \widetilde{Y}(\lambda;p)$ satisfying:
\begin{enumerate}
\item $\widetilde{Y}$ is analytic in $\mathbb{C} \setminus (-1,1)$.
\item $\widetilde{Y}$ has smooth (actually analytic) boundary values, $\widetilde{Y}_{\pm}(\lambda)$, on $(-1,1)$ which satisfy the jump relation
\begin{eqnarray}
\widetilde{Y}_{+}(\lambda) = \widetilde{Y}_{-}(\lambda) \pmtwo{1}{w_{p}(\varphi_p(\lambda))}{0}{1} \ , \ \ \lambda \in (-1,1) \ .
\end{eqnarray}
\item $\widetilde{Y}$ has the following asymptotic behavior:
\begin{eqnarray}
\widetilde{Y}(\lambda) = \left( I_2 + \mathcal{O} \left( \lambda^{-1}\right) \right) \lambda^{p \sigma_{3}} \ , \ \ \lambda \to \infty . 
\end{eqnarray}
\end{enumerate}
\end{RHP}
Now we proceed to reconstruct the matrix $T$ in the $\lambda$-plane. By the definition of $E$ in \eqref{eq:def E}, we obtain the explicit expression for $\widetilde{T}$:
\begin{gather}
\widetilde{T}(\lambda) =E(\lambda) \, \widetilde{T}_{\mathrm{approx}}(\lambda),
\end{gather}
where $\widetilde{T}_{\mathrm{approx}}$ is given in \eqref{eq:Tapprx}, and the matrices $\hat{T}$, $\widetilde{T}_1$, $\widetilde{T}_{-1}$ are all known explicitly (see \eqref{e:That}, \eqref{eq: T tilde +1} and \eqref{eq: T tilde -1}). Since $\widetilde{T}$ has now been determined, relation \eqref{eq:eq243} immediately yields an explicit expression for $T(\lambda)$. Recall that $\widetilde{Y}(\lambda)$ and $T$ are related via \eqref{eq: transformation T}, therefore 
\begin{equation}\label{e:ytil}
    \widetilde{Y}(\lambda)= e^{\frac{p \ell }{2}\sigma_3} T(\lambda) e^{p \left( g( \varphi_p(\lambda))-\frac{\ell}{2}\right)\sigma_3}
\end{equation}
provides an explicit expression for $\widetilde{Y}(\lambda)$.
\begin{remark}\label{e:remark21}
Recall that the $(1,1)$ entry of $\widetilde{Y}(\lambda)$ corresponds to the Jacobi polynomials $\mathbf{P}_p(\lambda)$, and the expression for $\mathbf{P}_{p-1}(\lambda)$ can be obtained by dividing the $(2,1)$ entry of $\widetilde{Y}(\lambda)$ by the factor $- 2 \pi i \, \kappa_{p-1}^2$ (see Proposition~\ref{RHP Y sol polynomials} for further details).
\end{remark}

In the following three propositions, we present the expressions for the Jacobi polynomials to the left and to the right of $\lambda=1$.

\begin{proposition}\label{Prop: Pm inside U1}
For values of $\lambda$ to the right of $\lambda=1$ and inside the disk $\mathcal{U}_1$, the expressions for the Jacobi polynomials $\mathbf{P}_p(\lambda)$ and $\mathbf{P}_{p-1}(\lambda)$ are given by
\begin{align*}
\mathbf{P}_p(\lambda) &= \frac{\sqrt{\pi} e^{\frac{p \ell}{2}}}{\sqrt{w_{p}(\varphi_p(\lambda))}} E_{11}(\lambda) \left( \chi_p^{1/4}(\lambda) \left( \frac{\lambda +1}{\lambda-1} \right)^{1/4} \mathrm{Ai}(\chi_p(\lambda)) - \chi_p^{-1/4}(\lambda) \left( \frac{\lambda -1}{\lambda+1} \right)^{1/4} \mathrm{Ai}'(\chi_p(\lambda)) \right) \nonumber \\ & + \frac{\sqrt{\pi} e^{\frac{p \ell}{2}}}{\sqrt{w_{p}(\varphi_p(\lambda))}} E_{12}(\lambda) \left( -i \chi_p^{1/4}(\lambda) \left( \frac{\lambda +1}{\lambda-1} \right)^{1/4} \mathrm{Ai}(\chi_p(\lambda)) - i \chi_p^{-1/4}(\lambda) \left( \frac{\lambda -1}{\lambda+1} \right)^{1/4} \mathrm{Ai}'(\chi_p(\lambda)) \right),
\end{align*}
and
\begin{align*}
\mathbf{P}_{p-1}(\lambda) &= -\frac{1}{2 \pi i \kappa^2_{p-1}}\frac{\sqrt{\pi} e^{\frac{p \ell}{2}}}{\sqrt{w_{p}(\varphi_p(\lambda))}} E_{21}(\lambda) \left( \chi_p^{1/4}(\lambda) \left( \frac{\lambda +1}{\lambda-1} \right)^{1/4} \mathrm{Ai}(\chi_p(\lambda)) - \chi_p^{-1/4}(\lambda) \left( \frac{\lambda -1}{\lambda+1} \right)^{1/4} \mathrm{Ai}'(\chi_p(\lambda)) \right)\\ & -\frac{1}{2 \pi i \kappa^2_{p-1}}\frac{\sqrt{\pi} e^{\frac{p \ell}{2}}}{\sqrt{w_{p}(\varphi_p(\lambda))}} E_{22}(\lambda) \left( -i \chi_p^{1/4}(\lambda) \left( \frac{\lambda +1}{\lambda-1} \right)^{1/4} \mathrm{Ai}(\chi_p(\lambda)) - i \chi_p^{-1/4}(\lambda) \left( \frac{\lambda -1}{\lambda+1} \right)^{1/4} \mathrm{Ai}'(\chi_p(\lambda)) \right).
\end{align*}
In the expressions for $\mathbf{P}_p(\lambda)$ and $\mathbf{P}_{p-1}(\lambda)$, the error matrix $E(\lambda)$ in~\eqref{eq. Error matrix} is being denoted by
\begin{equation}\label{eq:eqfunE}
E(\lambda)=\begin{pmatrix}
E_{11}(\lambda) & E_{12}(\lambda)\\ E_{21}(\lambda) & E_{22}(\lambda)
\end{pmatrix}.
\end{equation}
\end{proposition}

\begin{proof}
To the right of $\lambda=1$, we have that $T(\lambda) = \widetilde{T}(\lambda) \equiv E(\lambda) \, \widetilde{T}_{1}(\lambda)$. Therefore, after replacing $\widetilde{T}_1$ and the expression of $\Psi$ in region $I$ into \eqref{e:ytil}, we get:
\begin{align}\label{e:eqYtild}
\widetilde{Y}(\lambda) 
    &=e^{\frac{p \ell}{2}\sigma_3} E(\lambda) \, Q(\lambda) \, \mathcal{A}(\lambda) \, e^{-\frac{\pi i}{6}\sigma_3} e^{-\frac{p \ell}{2} \sigma_3} e^{p\left(g(\varphi_p(\lambda))-\phi(\lambda) \right)\sigma_3}
\end{align} 
where we used the notation  
\begin{equation}\label{eq:Airy}
    \mathcal{A}(\lambda):=\begin{pmatrix}
            \mathrm{Ai}(\chi_p(\lambda)) & \mathrm{Ai}(\omega^2 \chi_p(\lambda))\\
            \mathrm{Ai}'(\chi_p(\lambda)) & \omega^2 \mathrm{Ai}'(\omega^2 \chi_p(\lambda))
        \end{pmatrix}.
\end{equation}
Next, we replace the explicit expression for $\hat{T}(\lambda)$ (see \eqref{e:That}) into the definition of $Q(\lambda)$ \eqref{eq:eqQ}, which yields:
\begin{equation}\label{eq:Q}
    Q(\lambda)=\sqrt{\pi}\, e^{\frac{i\pi}{6}}\begin{pmatrix}
        \left(\frac{\lambda+1}{\lambda-1}\right)^{1/4} \left( \chi_p (\lambda) \right)^{1/4} & -\left(\frac{\lambda-1}{\lambda+1}\right)^{1/4} \left( \chi_p (\lambda) \right)^{-1/4}\\
        -i\,\left(\frac{\lambda+1}{\lambda-1}\right)^{1/4} \left( \chi_p (\lambda) \right)^{1/4} &-i\,\left(\frac{\lambda-1}{\lambda+1}\right)^{1/4} \left( \chi_p (\lambda) \right)^{-1/4}
    \end{pmatrix}. 
\end{equation}
Notice that to the right of $\lambda=1$ the functions $g$ and $\phi$ are analytic, and from their definitions we get
\begin{gather*}
g(\varphi_p(\lambda)) - \phi(\lambda) = c + \frac{1}{2} \left( - b_{\infty} \ln(z) + b_{\infty} \ln(\beta_p) - a_{\infty} \, p^{\nu} \ln(1-z) + a_{\infty} \, p^{\nu} \ln(1-\beta_p) \right).
\end{gather*}
Using equation \eqref{e:eq22}, the last one simplifies into
\begin{gather}\label{eq:eq54}
g(\varphi_p(\lambda)) - \phi(\lambda) = \frac{\ell}{2} - \frac{1}{2} \left( b_{\infty} \ln(z) + a_{\infty} p^{\nu} \ln(1-z) \right),
\end{gather}
which in turn implies the following simplification
\begin{equation}\label{eq:eqsimplification}
e^{-\frac{i \pi}{6}\sigma_3} e^{-\frac{p \ell}{2} \sigma_3} e^{p(g(\varphi_p(\lambda))-\phi(\lambda)) \sigma_3} = e^{-\frac{i \pi}{6}\sigma_3} w_{p}^{-\sigma_3/2}(\varphi_p(\lambda)).
\end{equation}
This reduces equation \eqref{e:eqYtild} into
\begin{equation}\label{eq:eytil}
  \widetilde{Y}(\lambda)= e^{\frac{p \ell}{2}\sigma_3} E(\lambda) \, Q(\lambda)\,\mathcal{A}(\lambda)\,e^{-\frac{\pi i}{6}\sigma_3} (w_p(\varphi_p(\lambda)))^{-\frac{\sigma_3}{2}}.
\end{equation}
Finally, we replace equations \eqref{eq:eqfunE}, \eqref{eq:Airy}, and \eqref{eq:Q} into equation \eqref{eq:eytil}. A straightforward calculation gives the $(1,1)$ and $(2,1)$ entries of $\widetilde{Y}(\lambda)$:
\begin{align}\label{eq:eqYtilde11}
    \widetilde{Y}_{11}(\lambda) = e^{\frac{p \ell}{2}} \Big(E_{11}(\lambda)\left(Q \mathcal{A}\right)(\lambda)_{11} + E_{12}(\lambda)  \left(Q \mathcal{A}\right)(\lambda)_{21}  \Big)\,e^{-\frac{i\pi}{6}}w_p^{-1/2}(\varphi_p(\lambda)) \nonumber \\
    \widetilde{Y}_{21}(\lambda) = e^{-\frac{p \ell}{2}} \Big(E_{21}(\lambda) \left(Q \mathcal{A}\right)(\lambda)_{11} + E_{22}(\lambda)  \left(Q \mathcal{A}\right)(\lambda)_{21}  \Big)\,e^{-\frac{i\pi}{6}}w_p^{-1/2}(\varphi_p(\lambda)) 
\end{align} 
where $\left(Q \mathcal{A}\right)(\lambda)_{11}$ and $\left(Q \mathcal{A}\right)(\lambda)_{21}$ denote the $(1,1)$ and $(2,1)$ entries of the matrix product $\left( Q \mathcal{A}\right)(\lambda)$, and are given explicitly by:
\begin{align}\label{eq QA}
  \left(Q \mathcal{A}\right)(\lambda)_{11} & \coloneqq \sqrt{\pi} \, e^{\frac{i\pi}{6}} \, \left( \chi_p^{1/4}(\lambda) \left( \frac{\lambda +1}{\lambda-1} \right)^{1/4} \mathrm{Ai}(\chi_p(\lambda)) - \chi_p^{-1/4}(\lambda) \left( \frac{\lambda -1}{\lambda+1} \right)^{1/4} \mathrm{Ai}'(\chi_p(\lambda)) \right) \nonumber \\
  \left(Q \mathcal{A}\right)(\lambda)_{21} & \coloneqq \sqrt{\pi} \, e^{\frac{i\pi}{6}} \, \left( -i \chi_p^{1/4}(\lambda) \left( \frac{\lambda +1}{\lambda-1} \right)^{1/4} \mathrm{Ai}(\chi_p(\lambda)) - i \chi_p^{-1/4}(\lambda) \left( \frac{\lambda -1}{\lambda+1} \right)^{1/4} \mathrm{Ai}'(\chi_p(\lambda)) \right)  .
\end{align}
The desired result then follows from Remark \ref{e:remark21}. 
\end{proof}
 
The following result establishes the asymptotics of Jacobi polynomilas on the left of $\lambda=1$ and inside the disk $\mathcal{U}_1$. Recall from Section~\ref{e:sec23} that in this case we must consider the asymtotic expresions for $Y$ on the region $\Omega_1$ (see Figure \ref{fig Contour T tilde}).
\begin{proposition}
For $\lambda \in \mathcal{U}_1\cap \Omega_1$, with $\lambda<1$,  the expressions for the Jacobi polynomials $\mathbf{P}_p(\lambda)$ and $\mathbf{P}_{p-1}(\lambda)$ are identical to those obtained in Proposition~\ref{Prop: Pm inside U1} for values of $\lambda$ to the right of $\lambda=1$ and inside the disk $\mathcal{U}_1$.
\end{proposition}

\begin{proof}
Let us now consider $\lambda \in \mathcal{U}_1\cap \Omega_1$, $\lambda<1$. In this case, we have that $T(\lambda) = \widetilde{T}(\lambda) \, U(\lambda) \equiv E(\lambda) \, \widetilde{T}_{1}(\lambda) \, U(\lambda)$. Therefore, after  replacing $\widetilde{T}_1$, $U$ and the expression of $\Psi$ in region $II$ into \eqref{e:ytil}, we get:
\begin{align}\label{e:eq265}
\widetilde{Y}(\lambda) 
    & = e^{\frac{p \ell}{2}\sigma_3}E(\lambda) Q(\lambda) \mathcal{A}(\lambda) \, e^{-\frac{\pi i}{6}\sigma_3}
        \begin{pmatrix}
            1 &0\\
            -1& 1
        \end{pmatrix}e^{-p \, \phi(\lambda) \sigma_3} \begin{pmatrix} 
        1 & 0\\ 
        e^{-2 p \phi(\lambda)} & 1\\
\end{pmatrix}e^{p\left(g(\varphi_p(\lambda))-\frac{\ell}{2}\right)\sigma_3} \nonumber \\
    &=e^{\frac{p \ell}{2}\sigma_3}E(\lambda) \, Q(\lambda) \mathcal{A}(\lambda) \, e^{-\frac{\pi i}{6}\sigma_3} e^{-\frac{p \ell}{2} \sigma_3} e^{p\left(g(\varphi_p(\lambda))-\phi(\lambda) \right)\sigma_3},
\end{align}
where $\mathcal{A}(\lambda)$ is as in \eqref{eq:Airy}. However, notice that for $\lambda \in (-1,1)$, the functions $g$ and $\phi$ are not analytic. Thus, in equation \eqref{e:eq265}, we need to account for their boundary values. From their definitions and that of the weight function, we obtain:
\begin{align*}
e^{p \left( g_{+}(\varphi_p(\lambda))-\phi_{+}(\lambda) \right) \sigma_3} & = e^{p \left[ g_{+}(\varphi_p(\lambda)) - \frac{1}{2} \left( g_{+}(\varphi_p(\lambda)) - g_{-}(\varphi_p(\lambda)) \right) \right] \sigma_3}\\
& = e^{ \frac{p \ell}{2} \sigma_3} e^{-\frac{p}{2} \left(b_{\infty} \ln(z) + a_{\infty} p^{\nu} \ln(1-z) \right)} \\
& = e^{ \frac{p \ell}{2} \sigma_3} w_{p}^{-\sigma_3/2} \left( \varphi_p(\lambda) \right).
\end{align*}
Substituting this relation and the definition of $Q(\lambda)$ into \eqref{e:eq265}, we derive the desired expression.
\end{proof}

The following result establishes the asymptotics of Jacobi polynomilas outside the disk $\mathcal{U}_1$ and outside $\Omega_1 \cup \Omega_2$ (see Figure \ref{fig Contour T tilde}).
\begin{proposition}
The expression for the Jacobi polynomials $\mathbf{P}_p(\lambda)$ and $\mathbf{P}_{p-1}(\lambda)$, valid for values of $\lambda \in \mathbb{C}\setminus(\Omega_1\cup \Omega_2)$ is given by:
\begin{align}
\mathbf{P}_p(\lambda) =\frac{i}{2}e^{p\,g(\varphi_p(\lambda))}\left(-if(\lambda) E_{11}(\lambda) + \widetilde{f}(\lambda) E_{12}(\lambda)\right)
\end{align}
\begin{align}
\mathbf{P}_{p-1}(\lambda) = - \frac{1}{4 \pi \kappa^2_{p-1}} e^{p(g(\varphi_p(\lambda))-l)}\left(-i f(\lambda) E_{21}(\lambda) + \widetilde{f}(\lambda) E_{22}(\lambda)\right)
\end{align}
where $f$ and $\widetilde{f}$ are as follows:
\begin{gather}\label{eq:eqfftil}
f(\lambda)=\left(\frac{\lambda-1}{\lambda+1}\right)^{1/4}+\left(\frac{\lambda-1}{\lambda+1}\right)^{-1/4}, \quad \widetilde{f}(\lambda)=\left(\frac{\lambda-1}{\lambda+1}\right)^{1/4}-\left(\frac{\lambda-1}{\lambda+1}\right)^{-1/4}.
\end{gather}
\end{proposition}
\begin{proof}
Finally, let us consider $\lambda \in \mathbb{C}\setminus(\Omega_1\cup \Omega_2)$. In this case, we have that $T(\lambda) = \widetilde{T}(\lambda) \equiv E(\lambda) \, \hat{T}(\lambda)$. Therefore, after replacing the definition of $\hat{T}(\lambda)$ into \eqref{e:ytil}, we get:
\begin{align*}
    \widetilde{Y}(\lambda)
    &=\frac{i}{2}
    \begin{pmatrix}
        e^{p\,g(\varphi_p(\lambda))}\left(-i f(\lambda) E_{11}(\lambda) + \widetilde{f}(\lambda) E_{12}(\lambda)\right) & e^{-p (g(\varphi_p(\lambda))-\ell)} \left(-\widetilde{f}(\lambda) E_{11}(\lambda) -i f(\lambda) E_{12}(\lambda)\right)\\
        e^{p(g(\varphi_p(\lambda))-\ell)}\left(-i f(\lambda) E_{21}(\lambda) + \widetilde{f}(\lambda) E_{22}(\lambda)\right) & e^{-p\,g(\varphi_p(\lambda))}\left(-\widetilde{f}(\lambda) E_{21}(\lambda) -i f(\lambda) E_{22}(\lambda)\right)
    \end{pmatrix},
\end{align*}
where $f(\lambda)$, $\widetilde{f}(\lambda)$ are given in \eqref{eq:eqfftil}. This immediately yields:
\begin{gather}\label{eq: Ytilde outside discs}
\widetilde{Y}_{11}(\lambda) =\frac{i}{2}e^{p\,g(\varphi_p(\lambda))}\left(-if(\lambda) E_{11}(\lambda) + \widetilde{f}(\lambda) E_{12}(\lambda)\right) \\
\widetilde{Y}_{21}(\lambda) =\frac{i}{2}e^{p(g(\varphi_p(\lambda))-l)}\left(-i f(\lambda) E_{21}(\lambda) + \widetilde{f}(\lambda) E_{22}(\lambda)\right)
\end{gather}
which completes the proof.
\end{proof}

\section{Proofs of auxiliary results}\label{proofs auxiliary}

In this section, we provide the proofs of the auxiliary results stated in Section~\ref{proofs}, which were used to establish our main results.

\subsection{Proof of Proposition~\ref{proposition Khat to K_airy}}
\label{proof Prop kernel convergence}
For clarity, we divide the proof of this proposition into several steps.
\par
\textbf{Notation.} In this step, we recall some equations that will be used in the proof and fix some notation. First, we recall equation \eqref{eq:eytil}:
\begin{equation*}
  \widetilde{Y}(\lambda)= e^{\frac{p \ell}{2}\sigma_3} E(\lambda) \, Q(\lambda)\,\mathcal{A}(\lambda)\,e^{-\frac{\pi i}{6}\sigma_3} (w_p(\lambda))^{-\frac{\sigma_3}{2}},
\end{equation*}
where $\mathcal{A}(\lambda)$ is defined in \eqref{eq:Airy}, and $Q(\lambda)$ is defined in \eqref{eq:Q}. It is convenient to express the $(i,j)$-entries of $Q(\lambda)$ in terms of the conformal map $\chi_p$. To do this, recall the function:
\begin{align*}
  \phi(\lambda) &= \int_{1}^{\lambda} (\zeta - 1)^{1/2} h_p(\zeta) \, d \zeta, \quad  h_p(\zeta) = \frac{1}{2} \left( \frac{\beta_p - \alpha_p}{2} \right)^2 (\zeta + 1)^{1/2}  \left( \frac{b_{\infty}}{R(0) s(\zeta)} + \frac{a_{\infty} p^{\nu}}{R(1) \big( s(\zeta) -1\big)} \right).
\end{align*}
Define
\begin{equation}\label{eq. def constant c1}
    c_1=c_1(p):=\left( -h_p(1)\right)^{2/3}>0
\end{equation}
(see also \eqref{eq:consh1}). The map $\chi_p$ in \eqref{eq:mapchi} can be written as 
\begin{equation}\label{eq. chi_p asymtotics}
\xi:= \chi_p(\lambda) =p^{\frac{2}{3}} \left(- \frac{3}{2} \, \,\phi(\lambda)\right)^{2/3}\equiv p^{\frac{2}{3}} c_1 \,
\Omega(\lambda),
\end{equation}
where 
\begin{equation}\label{eq. def Omega}
\Omega(\lambda)
=
(\lambda-1)
\left(
\sum_{j=0}^{\infty}
\frac{3h_p^{(j)}(1)}
{(2j+3)j!\,h_p(1)}
(\lambda-1)^j
\right)^{\frac{2}{3}}
\end{equation}
satisfies $\Omega(1)=0$ and $\Omega'(1)=1\neq 0$. Hence, for $|\xi|\le c_1 p^{2/3}$ we can write:
\begin{equation}\label{eq:eq515}
    \lambda = \chi_p^{-1}(\xi)=\Omega^{-1}\left(\frac{\xi}{c_1p^{2/3}}\right)=\sum_{j=0}^{\infty}a_j\left(\frac{\xi}{c_1p^{2/3}}\right)^j, \quad a_0=1=a_1.
\end{equation}
Therefore, we obtain:
\begin{equation}\label{eq:equationq}
    \frac{\lambda-1}{\xi}=\frac{1}{c_1p^{2/3}} \,q\left(\frac{\xi}{c_1p^{2/3}}\right), \quad q\left(\frac{\xi}{c_1p^{2/3}}\right)=\sum_{j=0}^{\infty}a_{j+1}\left(\frac{\xi}{c_1p^{2/3}}\right)^j,
\end{equation}
where $q$ is an analytic function in a neighborhood of $\xi=0$.

Substituting \eqref{eq:equationq} into \eqref{eq:Q} and setting $c_2=c_2(p):=\left(-h_p(1)\right)^{1/6}$, the entries of $Q(\lambda)$ can be written as follows:
\begin{align*}
 Q_{11}(\chi_p^{-1}(\xi))=\sqrt{\pi}\,e^{\frac{i\pi}{6}}\left(\frac{\lambda+1}{\lambda-1}\right)^{1/4}\xi^{1/4}=c_2\,p^{1/6}\sqrt{\pi}\,e^{\frac{i\pi}{6}}(\lambda+1)^{1/4}\left[q\left(\frac{\xi}{c_1p^{2/3}}\right)\right]^{-1/4}
\end{align*}
\begin{align*}
 Q_{21}(\chi_p^{-1}(\xi))=-i\,\sqrt{\pi}\,e^{\frac{i\pi}{6}}\left(\frac{\lambda+1}{\lambda-1}\right)^{1/4}\xi^{1/4}=-i\,c_2\,p^{1/6}\sqrt{\pi}\,e^{\frac{i\pi}{6}}(\lambda+1)^{1/4}\left[q\left(\frac{\xi}{c_1p^{2/3}}\right)\right]^{-1/4}
\end{align*}
\begin{align*}
 Q_{12}(\chi_p^{-1}(\xi))=-\sqrt{\pi}\,e^{\frac{i\pi}{6}}\left(\frac{\lambda-1}{\lambda+1}\right)^{1/4}\xi^{-1/4}=\frac{-1}{c_2\,p^{1/6}}\sqrt{\pi}\,e^{\frac{i\pi}{6}}(\lambda+1)^{-1/4}\left[q\left(\frac{\xi}{c_1p^{2/3}}\right)\right]^{1/4}
\end{align*}
\begin{align*}
 Q_{22}(\chi_p^{-1}(\xi)) =-i\sqrt{\pi}\,e^{\frac{i\pi}{6}}\left(\frac{\lambda-1}{\lambda+1}\right)^{1/4}\xi^{-1/4}=\frac{-i}{c_2\,p^{1/6}}\sqrt{\pi}\,e^{\frac{i\pi}{6}}(\lambda+1)^{-1/4}\left[q\left(\frac{\xi}{c_1p^{2/3}}\right)\right]^{1/4}.
\end{align*}
It will also be convenient to introduce the notation $\mathcal{B}(\lambda):=E(\lambda)Q(\lambda)$, in terms of which the entries of equation~\eqref{e:eq265} can be rewritten as
\begin{equation}\label{Eq. Ytilde BA}
    \widetilde{Y}_{11} (\lambda) = e^{\frac{p \ell}{2}}(\mathcal{B}\mathcal{A}(\lambda))_{11} \,e^{-\frac{i\pi}{6}}w_p^{-1/2}(\lambda),\quad     
    \widetilde{Y}_{21} (\lambda) = e^{-\frac{p \ell}{2}}(\mathcal{B}\mathcal{A} (\lambda))_{21} \,e^{-\frac{i\pi}{6}}w_p^{-1/2}(\lambda).
\end{equation}

Finally, we note that throughout the proof, we switch between the variables $\lambda$ and $\xi$ depending on the task at hand. In particular, whenever uniform bounds are required, we work in the $\lambda$-variable since it is more convenient for deriving uniform estimates (see for example \eqref{eq:est}).

We start with the definition of $\widehat{K}_p$ in \eqref{kernel local coordinates}, namely,
\begin{align}\label{e:Khat_p(xi1,xi2)}
\widehat{K}_p (\xi_1, \xi_2) 
&=-\frac{\sqrt{w_p(\varphi_p(\chi_p^{-1}(\xi_1))) w_p(\varphi_p(\chi_p^{-1}(\xi_2)))}}{2 \pi i\,\sqrt{\chi_p'(\chi_p^{-1}(\xi_1))\,\chi_p'(\chi_p^{-1}(\xi_2))}}\\
&\times \left(\frac{\widetilde{Y}_{11}(\chi_p^{-1}(\xi_1))\widetilde{Y}_{21}(\chi_p^{-1}(\xi_2))-\widetilde{Y}_{11}(\chi_p^{-1}(\xi_2))\widetilde{Y}_{21}(\chi_p^{-1}(\xi_1))}{\chi_p^{-1}(\xi_1)-\chi_p^{-1}(\xi_2)} \right).
\end{align}
Combining \eqref{Eq. Ytilde BA} and \eqref{e:Khat_p(xi1,xi2)}, we obtain 
\begin{equation}\label{Eq. Kernel hat Ytilde}
\begin{split}
\widehat{K}_p (\xi_1, \xi_2) 
&=-\frac{e^{-\frac{i\pi}{3}}}{2 \pi i\,\sqrt{\chi_p'(\chi_p^{-1}(\xi_1))\,\chi_p'(\chi_p^{-1}(\xi_2))}\,\left[\chi_p^{-1}(\xi_1)-\chi_p^{-1}(\xi_2)\right]}\\
& \times \left[(\mathcal{B}\mathcal{A})_{11}(\chi_p^{-1}(\xi_1))(\mathcal{B}\mathcal{A})_{21}(\chi_p^{-1}(\xi_2))-(\mathcal{B}\mathcal{A})_{11}(\chi_p^{-1}(\xi_2))(\mathcal{B}\mathcal{A})_{21}(\chi_p^{-1}(\xi_1)) \right].
\end{split}
\end{equation}
\par
\textbf{Analysis of the second term in the product \eqref{Eq. Kernel hat Ytilde}.} In this step, we analyze the second term in the product \eqref{Eq. Kernel hat Ytilde}, which, for the sake of brevity, is denoted by $\widehat{K}_{p}^{(s)}$. We have:
\begin{align}\label{eq:Kps}
\widehat{K}_{p}^{(s)}(\xi_1, \xi_2)  & = (\mathcal{B}\mathcal{A})_{11}(\chi_p^{-1}(\xi_1))(\mathcal{B}\mathcal{A})_{21}(\chi_p^{-1}(\xi_2))-(\mathcal{B}\mathcal{A})_{11}(\chi_p^{-1}(\xi_2))(\mathcal{B}\mathcal{A})_{21}(\chi_p^{-1}(\xi_1))\\ \nonumber
&= \mathrm{Ai}(\xi_1) \mathrm{Ai}'(\xi_2) \Delta_1 - \mathrm{Ai}(\xi_2) \mathrm{Ai}'(\xi_1) \Delta_2 + \mathrm{Ai}(\xi_1) \mathrm{Ai}(\xi_2) \Delta_3 + \mathrm{Ai}'(\xi_1) \mathrm{Ai}'(\xi_2) \Delta_4\\ \nonumber
&+\det[\mathcal{B}(\xi_1)]\left( \mathrm{Ai}(\xi_1) \mathrm{Ai}'(\xi_2) - \mathrm{Ai}(\xi_2) \mathrm{Ai}'(\xi_1) \right),
\end{align}
where we have introduced the notation:
\begin{align*}
    \Delta_1 &=\begin{vmatrix}
       \mathcal{B}_{11}(\chi_p^{-1}(\xi_1)) & \mathcal{B}_{12}(\chi_p^{-1}(\xi_2))- \mathcal{B}_{12}(\chi_p^{-1}(\xi_1))\\
       \mathcal{B}_{21}(\chi_p^{-1}(\xi_1))&  \mathcal{B}_{22}(\chi_p^{-1}(\xi_2))- \mathcal{B}_{22}(\chi_p^{-1}(\xi_1))
   \end{vmatrix},\\
   \Delta_2&=\begin{vmatrix}
       \mathcal{B}_{11}(\chi_p^{-1}(\xi_2))-\mathcal{B}_{11}(\chi_p^{-1}(\xi_1)) & \mathcal{B}_{12}(\chi_p^{-1}(\xi_1))\\
       \mathcal{B}_{21}(\chi_p^{-1}(\xi_2))- \mathcal{B}_{21}(\chi_p^{-1}(\xi_1))&  \mathcal{B}_{22}(\chi_p^{-1}(\xi_1))
   \end{vmatrix},\\
   \Delta_3&=\begin{vmatrix}
       \mathcal{B}_{11}(\chi_p^{-1}(\xi_1)) & \mathcal{B}_{11}(\chi_p^{-1}(\xi_2))- \mathcal{B}_{11}(\chi_p^{-1}(\xi_1))\\
       \mathcal{B}_{21}(\chi_p^{-1}(\xi_1))&  \mathcal{B}_{21}(\chi_p^{-1}(\xi_2))- \mathcal{B}_{21}(\chi_p^{-1}(\xi_1))
   \end{vmatrix},\\
   \Delta_4&=\begin{vmatrix}
       \mathcal{B}_{12}(\chi_p^{-1}(\xi_1)) & \mathcal{B}_{12}(\chi_p^{-1}(\xi_2))- \mathcal{B}_{12}(\chi_p^{-1}(\xi_1))\\
       \mathcal{B}_{22}(\chi_p^{-1}(\xi_1))&  \mathcal{B}_{22}(\chi_p^{-1}(\xi_2))- \mathcal{B}_{22}(\chi_p^{-1}(\xi_1))
   \end{vmatrix}.
\end{align*}
\par
First, we analyze the term $\det[\mathcal{B}(\xi_1)])$ in \eqref{eq:Kps}. For that, we write the error matrix $E$, defined in
\eqref{eq:def E}, as
\begin{equation}\label{eq:Emat}
E = I_2 + \frac{1}{p} \, \mathcal{E}, \quad \mathcal{E} = \begin{pmatrix}
\mathcal{E}_{11} & \mathcal{E}_{12}\\
\mathcal{E}_{21} & \mathcal{E}_{22}
\end{pmatrix}    
\end{equation}
where
\[
\mathcal{E}_{i,j}(\lambda)=\mathcal{O}(1),
\]
uniformly for $\lambda\in\mathbb{C}$,  $i,j\in \{1,2\}$ (see~\eqref{eq. Error matrix}). Here, recall that the implicit constant in $\mathcal{O}(1)$ depends only on the contour $\Sigma_{E}$, and the radius of the disks $\mathcal{U}_{-1}$ and $\mathcal{U}_{1}$ in Figure~\ref{fig Contour E} (see Remark~\ref{remark error matrix}). Moreover,  since $\mathcal{B}(\lambda):=E(\lambda)Q(\lambda)$ and $\det Q(\lambda)=-2i \pi e^{i \pi/3}$ we have the following exact expression:
\begin{equation}\label{eq:eqdet}
\det[\mathcal{B}(\xi_1)] = -2 \pi i \, e^{i \pi/3} \left( 1 + \frac{1}{p} \left( \cE_{11}(\xi_1) + \cE_{22}(\xi_1)\right) + \frac{1}{p^2} \, \det[\mathcal{E}(\xi_1)]\right).
\end{equation} 
In what follows, we analyze the determinants $\Delta_j$ and derive their
asymptotic estimates. We first group the terms according to their explicit
powers of $p$. The coefficients of these powers may themselves depend on $p$,
so we estimate them separately. We then combine these estimates with the displayed
powers of $p$, which gives the actual asymptotic order of each determinant. We have:
\begin{align}\label{e:dets}
\Delta_1&=\Delta_1^{(0)}+\Delta_1^{(-1)}p^{-1}+\Delta_1^{(-2)}p^{-2}, \qquad
\Delta_2=\Delta_2^{(0)}+\Delta_2^{(-1)}p^{-1}+\Delta_2^{(-2)}p^{-2}\\\\
\Delta_3&=  \Delta_3^{(-2/3)} \, p^{-2/3}+\Delta_3^{(-5/3)} \, p^{-5/3}, \qquad
\Delta_4=  \Delta_4^{(-4/3)} \, p^{-4/3}+\Delta_4^{(-7/3)} \, p^{-7/3}.
\end{align}
\textbf{Analysis of the determinant $\Delta_3$.} The coefficient $\Delta_3^{(-2/3)} $ is given by:
\begin{align}\label{eq:est}
 \Delta_3^{(-2/3)} &= (-1)^{1/3} \pi c_2^2 (1+\lambda_1)^{1/4} (1+\lambda_2)^{1/4} \left( \frac{i\left[ \cE_{11}(\chi_p^{-1}(\xi_2)) - \cE_{11}(\chi_p^{-1}(\xi_1))  \right] + \left[ \cE_{12}(\chi_p^{-1}(\xi_2)) - \cE_{12}(\chi_p^{-1}(\xi_1)) \right] }{\left[q\left(\frac{\xi_1}{c_1p^{2/3}}\right)\right]^{1/4} \left[q\left(\frac{\xi_2}{c_1p^{2/3}}\right)\right]^{1/4}} \right)  \\ \nonumber
& + (-1)^{1/3} \pi c_2^2 (1+\lambda_1)^{1/4} (1+\lambda_2)^{1/4} \left( \frac{\left[ \cE_{21}(\chi_p^{-1}(\xi_2)) - \cE_{21}(\chi_p^{-1}(\xi_1)) \right] + i \left[ \cE_{22}(\chi_p^{-1}(\xi_1)) - \cE_{22}(\chi_p^{-1}(\xi_2)) \right]}{\left[q\left(\frac{\xi_1}{c_1p^{2/3}}\right)\right]^{1/4} \left[q\left(\frac{\xi_2}{c_1p^{2/3}}\right)\right]^{1/4}}\right).
\end{align}
Before we continue with the analysis, we fix for the rest of the paper the following notation:
\[w \coloneqq \chi_p^{-1}(\xi).\]
Since the function $\frac{(-1)^{1/3} \pi c_2^2 (1+\lambda_1)^{1/4} (1+\lambda_2)^{1/4}}{\left[q\left(\frac{\xi_1}{c_1p^{2/3}}\right)\right]^{1/4} \left[q\left(\frac{\xi_2}{c_1p^{2/3}}\right)\right]^{1/4}}$ is bounded, there exists a constant $C$ such that 
\begin{align*}
\left|\Delta_3^{(-2/3)}\right|
&\leq C \Big(
\left|\cE_{11}(w_2)
      -\cE_{11}(w_1)\right|
+
\left|\cE_{12}(w_2)
      -\cE_{12}(w_1)\right|+
\left|\cE_{21}(w_2)
      -\cE_{21}(w_1)\right|
+
\left|\cE_{22}(w_1)
      -\cE_{22}(w_2)\right|
\Big).
\end{align*}
Since, for $i,j=1,2$, $\mathcal{E}_{ij}=\mathcal{O}(1)$ and $\mathcal{E}'_{ij}=\mathcal{O}(1)$, uniformly as $p \to \infty$, we use Lemma~\ref{lem: Bound ElamdaQ2-ElamdaQ1} for $H=\mathcal{E}$, to bound the differences $\cE_{ij}(w_2) - \cE_{ij}(w_1)$. Thus, there exists a constant $c>0$ independent of $p$ such that:
\begin{align}\label{eq:est1}
\boxed{|\Delta_3^{(-2/3)}| \leq c \, \, p^{-2/3} |\xi_2 - \xi_1|}
\end{align}
The coefficient $\Delta_3^{(-5/3)} $ is given by:
\begin{align}
\Delta_3^{(-5/3)} &= \frac{(-1)^{1/3} \pi c_2^2 (1+\lambda_1)^{1/4} (1+\lambda_2)^{1/4}}{\left[q\left(\frac{\xi_1}{c_1p^{2/3}}\right)\right]^{1/4} \left[q\left(\frac{\xi_2}{c_1p^{2/3}}\right)\right]^{1/4}}( -(\cE_{11}(w_2)-i\cE_{12}(w_2))(\cE_{21}(w_1)-i\cE_{22}(w_1))\\&+(\cE_{11}(w_1)-i\cE_{12}(w_1))(\cE_{21}(w_2)-i\cE_{22}(w_2))).
\end{align}
Similar to $\Delta_3^{(-2/3)} $, we use Lemma~\ref{lem: Bound ElamdaQ2-ElamdaQ1} to bound differences of the form $\cE_{ij}(w_2) - \cE_{ij}(w_1)$. To produce such differences, we first add and subtract the following product:
\[
\left(\mathcal{E}_{11}(w_1)-i\mathcal{E}_{12}(w_1)\right)
\left(\mathcal{E}_{21}(w_1)-i\mathcal{E}_{22}(w_1)\right).
\]
Doing so, we obtain:
\begin{align*}
&-\left(\mathcal{E}_{11}(w_2)-i\mathcal{E}_{12}(w_2)\right)
 \left(\mathcal{E}_{21}(w_1)-i\mathcal{E}_{22}(w_1)\right)+\left(\mathcal{E}_{11}(w_1)-i\mathcal{E}_{12}(w_1)\right)
 \left(\mathcal{E}_{21}(w_2)-i\mathcal{E}_{22}(w_2)\right)\\
&=
\Bigl[
\mathcal{E}_{11}(w_1)-\mathcal{E}_{11}(w_2)
-i\bigl(\mathcal{E}_{12}(w_1)-\mathcal{E}_{12}(w_2)\bigr)
\Bigr]
\left(\mathcal{E}_{21}(w_1)-i\mathcal{E}_{22}(w_1)\right)\\
&\qquad
+\left(\mathcal{E}_{11}(w_1)-i\mathcal{E}_{12}(w_1)\right)
\Bigl[
\mathcal{E}_{21}(w_2)-\mathcal{E}_{21}(w_1)
-i\bigl(\mathcal{E}_{22}(w_2)-\mathcal{E}_{22}(w_1)\bigr)
\Bigr].
\end{align*}
Therefore, by the triangle inequality:
\begin{align*}
&\Big|
-\left(\mathcal{E}_{11}(w_2)-i\mathcal{E}_{12}(w_2)\right)
 \left(\mathcal{E}_{21}(w_1)-i\mathcal{E}_{22}(w_1)\right)
+\left(\mathcal{E}_{11}(w_1)-i\mathcal{E}_{12}(w_1)\right)
 \left(\mathcal{E}_{21}(w_2)-i\mathcal{E}_{22}(w_2)\right)
\Big|\\
&\leq
\left(
\left|\mathcal{E}_{11}(w_1)-\mathcal{E}_{11}(w_2)\right|
+
\left|\mathcal{E}_{12}(w_1)-\mathcal{E}_{12}(w_2)\right|
\right)
\left(
|\mathcal{E}_{21}(w_1)|+|\mathcal{E}_{22}(w_1)|
\right)\\
&\qquad
+
\left(
|\mathcal{E}_{11}(w_1)|+|\mathcal{E}_{12}(w_1)|
\right)
\left(
\left|\mathcal{E}_{21}(w_2)-\mathcal{E}_{21}(w_1)\right|
+
\left|\mathcal{E}_{22}(w_2)-\mathcal{E}_{22}(w_1)\right|
\right).
\end{align*}
Directly, by Lemma~\ref{lem: Bound ElamdaQ2-ElamdaQ1}, together with the fact that $\frac{(-1)^{1/3} \pi c_2^2 (1+\lambda_1)^{1/4} (1+\lambda_2)^{1/4}}{\left[q\left(\frac{\xi_1}{c_1p^{2/3}}\right)\right]^{1/4} \left[q\left(\frac{\xi_2}{c_1p^{2/3}}\right)\right]^{1/4}}$ is bounded, there exists a constant $c>0$ independent of $p$, such that:
\begin{align}\label{eq:est2}
\boxed{|\Delta_3^{(-5/3)}| \leq c \, \, p^{-2/3} |\xi_2 - \xi_1|}
\end{align}
Combining the two estimates \eqref{eq:est1}-\eqref{eq:est2} yields that there exists a constant $c>0$ independent of $p$ such that:
\begin{equation}\label{eq. bound delta3}
\boxed{|\Delta_3| \leq c \, \, p^{-4/3} |\xi_2 - \xi_1|}
\end{equation}
\textbf{Analysis of the determinant $\Delta_4$.} The coefficient $\Delta_4^{(-4/3)} $ is given by:
\begin{align*}
\Delta_4^{(-4/3)} &=\frac{(-1)^{1/3} \pi \left[q\left(\frac{\xi_1}{c_1 p^{2/3}}\right)\right]^{1/4} \left[q\left(\frac{\xi_2}{c_1p^{2/3}}\right)\right]^{1/4}}{c_2^2 (1+\lambda_1)^{1/4} (1+\lambda_2)^{1/4}}( i\left( \cE_{11}(w_1) - \cE_{11}(w_2)  \right) + \left( \cE_{12}(w_2) - \cE_{12}(w_1) \right) \\&+\left(\cE_{21}(w_2) - \cE_{21}(w_1) \right) + i \left( \cE_{22}(w_2) - \cE_{22}(w_1) \right) )
\end{align*}
Since the pre-factor function in $\Delta_4^{(-4/3)}$ is bounded, we use a similar argument as before to show that there exists a constant $c>0$ independent of $p$ such that:
\begin{align*}
|\Delta_4^{(-4/3)}| \leq c \, \, p^{-2/3} |\xi_2 - \xi_1|
\end{align*}
The coefficient $\Delta_4^{(-7/3)}$ admits a representation analogous to that of $\Delta_3^{(-5/3)}$, and therefore a similar argument yields
\begin{align*}
|\Delta_4^{(-7/3)}| \leq c \, \, p^{-2/3} |\xi_2 - \xi_1|.
\end{align*}
Combining these two, we show that there exists constant $c>0$ independent of $p$ such that:
\begin{equation}\label{eq. bound delta 4}
\boxed{|\Delta_4| \leq c \, p^{-2} |\xi_2 - \xi_1|}
\end{equation}
\textbf{Analysis of the determinant $\Delta_1$.} Recall that 
\[\Delta_1=\Delta_1^{(0)}+\Delta_1^{(-1)}p^{-1}+\Delta_1^{(-2)}p^{-2}.\]
The coefficient $\Delta_1^{(0)}$ is given by:
\begin{align*}
\Delta_1^{(0)}= 2 \pi i e^{i \pi/3}- 2 \pi i e^{i \pi/3} \frac{(1+\lambda_1)^{1/4} \left[q\left(\frac{\xi_2}{c_1p^{2/3}}\right)\right]^{1/4}}{(1+\lambda_2)^{1/4} \left[q\left(\frac{\xi_1}{c_1p^{2/3}}\right)\right]^{1/4}}.
\end{align*}
Using equation \eqref{eq:eq515}, we can rewrite $\Delta_1^{(0)}$ as follows:
\begin{align*}
\Delta_1^{(0)}&=\frac{2 \pi i e^{i \pi/3}}{\left( 2+\frac{\xi_2}{c_1 p^{2/3}} \right)^{1/4}\left[q\left(\frac{\xi_1}{c_1p^{2/3}}\right)\right]^{1/4}} \left( \left( 2+\frac{\xi_2}{c_1 p^{2/3}} \right)^{1/4}\left[q\left(\frac{\xi_1}{c_1p^{2/3}}\right)\right]^{1/4} - \left( 2+\frac{\xi_1}{c_1 p^{2/3}} \right)^{1/4}\left[q\left(\frac{\xi_2}{c_1p^{2/3}}\right)\right]^{1/4}\right)
\end{align*}
and the expression in the parenthesis can be written as follows:
\begin{align}\label{e:parenth}
2^{1/4} \left( \left[q\left(\frac{\xi_1}{c_1p^{2/3}}\right)\right]^{1/4} - \left[q\left(\frac{\xi_2}{c_1p^{2/3}}\right)\right]^{1/4} \right) +q_0^{1/4} \left( \left( 2+\frac{\xi_2}{c_1 p^{2/3}} \right)^{1/4} - \left( 2+\frac{\xi_1}{c_1 p^{2/3}} \right)^{1/4}\right)\nonumber \\
+ \left( \left( 2+\frac{\xi_2}{c_1 p^{2/3}} \right)^{1/4} -2^{1/4}\right) \left( \left[q\left(\frac{\xi_1}{c_1 p^{2/3}}\right)\right]^{1/4} -q_0^{1/4} \right)-\left( \left( 2+\frac{\xi_1}{c_1 p^{2/3}} \right)^{1/4} -2^{1/4}\right) \left( \left[q\left(\frac{\xi_2}{c_1 p^{2/3}}\right)\right]^{1/4} -q_0^{1/4} \right).
\end{align}
We begin by considering the first two differences. We write:
\begin{align*}
\left[q\left(\frac{\xi}{c_1p^{2/3}}\right)\right]^{1/4}&= q_0^{1/4} + \frac{q_1 \, \xi }{4 c_1 \, q_0^{3/4}} \, p^{-2/3}+\left[q\left(\frac{\xi}{c_1p^{2/3}}\right)\right]^{1/4}-\left(q_0^{1/4} + \frac{q_1 \, \xi }{4 c_1 \, q_0^{3/4}} \, p^{-2/3} \right)
\end{align*}
and thus the first difference in \eqref{e:parenth} is:
\begin{align}\label{e:firstdf}
\frac{q_1}{4 c_1 \, q_0^{3/4}} \,(\xi_1-\xi_2) \, p^{-2/3}+ K(\xi_1)-K(\xi_2)
\end{align}
where 
\begin{align*}
K(\xi)&= \left[q\left(\frac{\xi}{c_1p^{2/3}}\right)\right]^{1/4}-\left(q_0^{1/4} + \frac{q_1 \, \xi }{4 c_1 \, q_0^{3/4}} \, p^{-2/3} \right).
\end{align*}
Since $q$ is an analytic function near $\xi=0$, it admits a convergent series for $|\xi|< c_1\, p^{2/3}$ and thus, we have:
\[K(\xi)=\sum_{j=2}^\infty \hat{q}_j \left( \frac{\xi}{c_1 p^{2/3}}\right)^j.\]
Then:
\[K'(\xi)=\frac{1}{c_1 p^{2/3}} \sum_{j=2}^\infty j \, \hat{q}_j \left( \frac{\xi}{c_1 p^{2/3}}\right)^{j-1}= \frac{\xi}{c_1^2 p^{4/3}} \sum_{j=2}^\infty j \, \hat{q}_j  \left( \frac{\xi}{c_1 p^{2/3}}\right)^{j-2} \] 
where the series is uniformly bounded on compact subsets of $|\xi|<c_1 p^{2/3}$. Therefore, there exists constant $C$ independent of $p$ such that:
\begin{equation}\label{e:est1}
|K'(\xi)| \leq C \,\, |\xi| p^{-4/3}.    
\end{equation}
Similarly, we can write:
\begin{align*}
\left( 2+\frac{\xi}{c_1 p^{2/3}} \right)^{1/4}=2^{1/4} + 2^{1/4} \frac{\xi}{8 c_1} p^{-2/3} + \left( 2+\frac{\xi}{c_1 p^{2/3}} \right)^{1/4} -\left( 2^{1/4} + 2^{1/4} \frac{\xi}{8 c_1} p^{-2/3}\right)
\end{align*}
and thus the second difference in \eqref{e:parenth} becomes
\begin{align}\label{e:secdf}
 \frac{ 2^{1/4}}{8 c_1} \, (\xi_2-\xi_1) \, p^{-2/3} + \Lambda(\xi_2)-\Lambda(\xi_1)
\end{align}
where
\[\Lambda(\xi)= \left( 2+\frac{\xi}{c_1 p^{2/3}} \right)^{1/4} -\left( 2^{1/4} + 2^{1/4} \frac{\xi}{8 c_1} p^{-2/3}\right).\]
Since the function $\left( 2+\frac{\xi}{c_1 p^{2/3}} \right)^{1/4}$ is analytic near $\xi=0$, it admits a convergent infinite series for $|\xi|<c_1 p^{2/3}$:
\begin{align*}
\Lambda(\xi)&= \sum_{j=2}^\infty \hat{\lambda}_j \left( \frac{\xi}{c_1 p^{2/3}}\right)^j.
\end{align*}
Repeating the previous argument for $K(\xi)$, we can find a constant $C$ independent of $p$ such that
\begin{equation}\label{e:est2}
|\Lambda'(\xi)| \leq C \,\, |\xi| p^{-4/3}.  
\end{equation}
For the second part in equation \eqref{e:parenth}, we notice that: 
\begin{align*}
\left( 2+\frac{\xi}{c_1 p^{2/3}} \right)^{1/4} -2^{1/4} &= \frac{2^{1/4}}{8 c_1} \xi \, p^{-2/3} A(\xi), \quad A(\xi) = \sum_{j=0}^\infty \hat{a}_j \left( \frac{\xi}{c_1 p^{2/3}}\right)^j, \quad |\xi|<c_1 p^{2/3}\\\\
\left[q\left(\frac{\xi}{c_1 p^{2/3}}\right)\right]^{1/4} -q_0^{1/4} &= \frac{q_1}{4 c_1 q_0^{3/4}} \xi \, p^{-2/3} B(\xi), \quad B(\xi) =\sum_{j=0}^\infty \hat{b}_j \left( \frac{\xi}{c_1 p^{2/3}}\right)^j, \quad |\xi|<c_1 p^{2/3}
\end{align*}
where $A(0)=1$, $B(0)=1$, and $|A'(\xi)| \leq C_1 p^{-2/3}$, $|B'(\xi)| \leq C_2 p^{-2/3}$ for some constants $C_1$, $C_2$. Thus, we can write the second part of \eqref{e:parenth} as:
 \begin{align}\label{e:est3}
     &\left( \left( 2+\frac{\xi_2}{c_1 p^{2/3}} \right)^{1/4} -2^{1/4}\right) \left( \left[q\left(\frac{\xi_1}{c_1 p^{2/3}}\right)\right]^{1/4} -q_0^{1/4} \right)-\left( \left( 2+\frac{\xi_1}{c_1 p^{2/3}} \right)^{1/4} -2^{1/4}\right) \left( \left[q\left(\frac{\xi_2}{c_1 p^{2/3}}\right)\right]^{1/4} -q_0^{1/4} \right) \nonumber\\ 
    &= \frac{2^{1/4} \, q_1}{32 c_1^2 q_0^{3/4}}\, \xi_1 \xi_2 \,  p^{-4/3} \left[ B(\xi_1) (A(\xi_2)-A(\xi_1))-A(\xi_1) (B(\xi_2)-B(\xi_1))\right]\\ \nonumber
    & \leq \left| \frac{2^{1/4} \, q_1}{32 c_1^2 q_0^{3/4}}\right| \, |\xi_1| |\xi_2| \,  p^{-4/3} \left[ |B(\xi_1)| |A(\xi_2)-A(\xi_1)|+|A(\xi_1)| |B(\xi_2)-B(\xi_1)|\right]
 \end{align}
 where $|B(\xi_1)|$, $|A(\xi_1)|$ are bounded, $|A(\xi_2)-A(\xi_1)| \leq C p^{-2/3}|\xi_2-\xi_1|$, and the same estimate holds for $|B(\xi_2)-B(\xi_1)|$. We now recall equation \eqref{e:parenth}, which becomes:
\begin{align}
& \frac{2^{1/4}}{8 c_1 q_0^{3/4}} (2q_1-q_0)p^{-2/3}(\xi_1-\xi_2)+ 2^{1/4}\left(K(\xi_1)-K(\xi_2) \right)+q_0^{1/4}\left(\Lambda(\xi_2)-\Lambda(\xi_1) \right)\\&
+\left( \left( 2+\frac{\xi_2}{c_1 p^{2/3}} \right)^{1/4} -2^{1/4}\right) \left( \left[q\left(\frac{\xi_1}{c_1 p^{2/3}}\right)\right]^{1/4} -q_0^{1/4} \right)-\left( \left( 2+\frac{\xi_1}{c_1 p^{2/3}} \right)^{1/4} -2^{1/4}\right) \left( \left[q\left(\frac{\xi_2}{c_1 p^{2/3}}\right)\right]^{1/4} -q_0^{1/4} \right)\\
&=\frac{2^{1/4}}{8 c_1 q_0^{3/4}} (2q_1-q_0)p^{-2/3}(\xi_1-\xi_2)+ \mathcal{O}\left(p^{-4/3} \left( |\xi_1|+|\xi_2|\right) |\xi_1-\xi_2| \right) + \mathcal{O} \left( p^{-6/3} |\xi_1| |\xi_2||\xi_1-\xi_2|\right).
\end{align}
Here, we used equations \eqref{e:firstdf}, \eqref{e:est1}, \eqref{e:secdf}, \eqref{e:est2} and \eqref{e:est3}. Moreover:
\begin{align*}
\Delta_1^{(0)}&=\frac{2 \pi i e^{i \pi/3}}{\left( 2+\frac{\xi_2}{c_1 p^{2/3}} \right)^{1/4}\left[q\left(\frac{\xi_1}{c_1p^{2/3}}\right)\right]^{1/4}} \\ & \times \left( \frac{2^{1/4}}{8 c_1 q_0^{3/4}} (2q_1-q_0)p^{-2/3}(\xi_1-\xi_2)+ \mathcal{O}\left(p^{-4/3} \left( |\xi_1|+|\xi_2|\right) |\xi_1-\xi_2| \right) + \mathcal{O} \left( p^{-6/3} |\xi_1| |\xi_2||\xi_1-\xi_2|\right)\right)
\end{align*}
where the pre-factor is bounded for large $p$, and for the present estimate only its leading-order scalar contribution is relevant, namely
\[
\frac{2\pi i e^{i\pi/3}}
{\left(2+\frac{\xi_2}{c_1p^{2/3}}\right)^{1/4}
\left[q\left(\frac{\xi_1}{c_1p^{2/3}}\right)\right]^{1/4}}
=
\frac{(-1)^{5/6}2^{3/4}\pi}{q_0^{1/4}}
+\mathcal{O}\!\left(p^{-2/3}\right).
\]
Thus:
\begin{align}\label{e:box1}
\boxed{\Delta_1^{(0)} = \left( \frac{(-1)^{5/6} \pi}{4 c_1 q_0} (2 q_1-q_0)\right) p^{-2/3}(\xi_1-\xi_2)+\mathcal{O}\left(p^{-4/3} \left( |\xi_1|+|\xi_2|\right) |\xi_1-\xi_2| \right) + \mathcal{O} \left( p^{-6/3} |\xi_1| |\xi_2||\xi_1-\xi_2|\right)}
\end{align}
Next, the coefficient $\Delta_1^{(-1)}$ is given by:
\[
\begin{aligned}
((-1)^{5/6}\pi)^{-1}\Delta_1^{(-1)}
&=
2\bigl(\cE_{11}(w_1)+\cE_{22}(w_1)\bigr)
-\frac{F(\xi_1)}{F(\xi_2)}P(w_1,w_2) \\
\end{aligned}
\]
where we have introduced the notation:
\begin{equation}\label{e:funF}
F(\xi):=
\frac{\left(2+\frac{\xi}{c_1p^{2/3}}\right)^{1/4}}
{\left[q\left(\frac{\xi}{c_1p^{2/3}}\right)\right]^{1/4}}\equiv \frac{2^{1/4}}{q_0^{1/4}} \sum_{j=0}^{\infty} f_j \left( \frac{\xi}{c_1 \, p^{2/3}}\right)^j
\end{equation}
\begin{align*}
P(w_1,w_2)
&=
\cE_{11}(w_1)-\cE_{11}(w_2)
+\cE_{22}(w_1)-\cE_{22}(w_2) \\
&\quad
+i\bigl(\cE_{12}(w_2)-\cE_{12}(w_1)\bigr)
+i\bigl(\cE_{21}(w_1)-\cE_{21}(w_2)\bigr) \\
&\quad
+2\cE_{11}(w_2)+2\cE_{22}(w_2).
\end{align*}
$F(\xi)$ satisfies $F(0)=\frac{2^{1/4}}{q_0^{1/4}}$, and:
\[F'(\xi)=\frac{2^{1/4}}{q_0^{1/4}}  \, \frac{1}{c_1 \,p^{2/3}} \sum_{j=1}^{\infty} j \, f_j \left( \frac{\xi}{c_1 \, p^{2/3}}\right)^{j-1}.\]
$F$ and $F^{-1}$ are uniformly bounded on compact subsets of $|\xi|<c_1 \, p^{2/3}$, and \
\[
|F'(\xi)|\leq C \,\, p^{-2/3},
\]
which implies 
\begin{equation}\label{e:est4}
|F(\xi_1)-F(\xi_2)|
\leq C \,\, p^{-2/3}|\xi_1-\xi_2|.    
\end{equation}
We rewrite $\Delta_1^{(-1)}$ as follows:
\[
\begin{aligned}
((-1)^{5/6}\pi)^{-1}\Delta_1^{(-1)}
&=
\Bigl[
2\bigl(\cE_{11}(w_1)+\cE_{22}(w_1)\bigr)
-P(w_1,w_2)
\Bigr] 
-\left(
\frac{F(\xi_1)}{F(\xi_2)}-1
\right)P(w_1,w_2).
\end{aligned}
\]
The term:
\[
2\bigl(\cE_{11}(w_1)+\cE_{22}(w_1)\bigr)
-P(w_1,w_2)\]
reduces into:
\begin{align}\label{eq:leadterm}
\cE_{11}(w_1)-\cE_{11}(w_2)+i(\cE_{12}(w_1)-\cE_{12}(w_2))-i(\cE_{21}(w_1)-\cE_{21}(w_2))+\cE_{22}(w_1)-\cE_{22}(w_2)
\end{align}
which we can bound by Lemma~\ref{lem: Bound ElamdaQ2-ElamdaQ1}. Therefore, there exists a constant $c>0$ independent of $p$ such that:
\[
\left|
2\bigl(\cE_{11}(w_1)+\cE_{22}(w_1)\bigr)
-P(w_1,w_2)
\right|
\leq
c \,\, p^{-2/3}|\xi_1-\xi_2|.
\]
Notice that $P(w_1,w_2)$ is uniformly bounded since all the differences of the form $\cE_{ij}(w_1)-\cE_{ij}(w_2)$ are bounded by Lemma~\ref{lem: Bound ElamdaQ2-ElamdaQ1}, and $\cE_{ij}(w)$ are uniformly bounded in $p$. Thus, equation \eqref{e:est4} gives:
\[
\left|
\left(
\frac{F(\xi_1)}{F(\xi_2)}-1
\right)
P(w_1,w_2)
\right|
\leq  \left| \frac{F(\xi_1)-F(\xi_2)}{F(\xi_2)}\right| \left|P(w_1,w_2) \right| \leq
C \,\, p^{-2/3}|\xi_1-\xi_2|.
\]
Therefore,
\[
|\Delta_1^{(-1)}|
\leq
C \,\, p^{-2/3}|\xi_1-\xi_2|,
\]
and equivalently,
\[
\boxed{\Delta_1^{(-1)} \, p^{-1}
=\mathcal{O}\!\left(p^{-5/3} (\xi_1-\xi_2)\right)}
\]
The coefficient $\Delta_1^{(-2)}$ is given by:
\begin{align*}
\left((-1)^{1/3} \pi \right)^{-1} \Delta_1^{(-2)} &= 2i \left( -\cE_{12}(w_1) \cE_{21}(w_1) +\cE_{11}(w_1) \cE_{22}(w_1) \right) +\frac{F(\xi_1)}{F(\xi_2)} P(w_1,w_2)
\end{align*}
where $F(\xi)$ as in equation \eqref{e:funF}, and here $P(w_1,w_2)$ is defined as:
\[
\begin{aligned}
P(w_1,w_2)
&=
\left( \cE_{11}(w_2)+i \cE_{12}(w_2)\right)\left(\cE_{21}(w_1)-i \cE_{22}(w_1) \right)-\left(\cE_{11}(w_1)-i \cE_{12}(w_1)\right)\left(\cE_{21}(w_2)+i \cE_{22}(w_2)\right).
\end{aligned}
\]
We rewrite $\Delta_1^{(-2)}$ as:
\begin{align*}
    \left((-1)^{1/3} \pi \right)^{-1} \Delta_1^{(-2)} &=\left[2i \left( -\cE_{12}(w_1) \cE_{21}(w_1) +\cE_{11}(w_1) \cE_{22}(w_1) \right) +P(w_1,w_2)\right]+\left( \frac{F(\xi_1)}{F(\xi_2)} -1\right)P(w_1,w_2).
\end{align*}
The term:
\[ 2i \left( -\cE_{12}(w_1) \cE_{21}(w_1) +\cE_{11}(w_1) \cE_{22}(w_1) \right) +P(w_1,w_2)\]
writes explicitly:
\begin{align*}
&2i \left( -\cE_{12}(w_1) \cE_{21}(w_1) +\cE_{11}(w_1) \cE_{22}(w_1) \right) +
    \left( \cE_{11}(w_2)+i \cE_{12}(w_2)\right)\left(\cE_{21}(w_1)-i \cE_{22}(w_1) \right)\\&-\left(\cE_{11}(w_1)-i \cE_{12}(w_1)\right)\left(\cE_{21}(w_2)+i \cE_{22}(w_2)\right)
\end{align*}
and can be reduced into
\[
\begin{aligned}
&\left[
\mathcal E_{11}(w_2)-\mathcal E_{11}(w_1)
+i\left(\mathcal E_{12}(w_2)-\mathcal E_{12}(w_1)\right)
\right]
\left(\mathcal E_{21}(w_1)-i\mathcal E_{22}(w_1)\right)
\\
&\quad
-\left(\mathcal E_{11}(w_1)-i\mathcal E_{12}(w_1)\right)
\left[
\mathcal E_{21}(w_2)-\mathcal E_{21}(w_1)
+i\left(\mathcal E_{22}(w_2)-\mathcal E_{22}(w_1)\right)
\right].
\end{aligned}
\]
Using Lemma~\ref{lem: Bound ElamdaQ2-ElamdaQ1} and the fact that all entries $\cE_{ij}$ are uniformly bounded in $p$, we can show that the above term is bounded by a constant times $p^{-2/3} |\xi_1-\xi_2|$. Moreover, notice that $P(w_1,w_2)$ is uniformly bounded in $p$ since it consists of terms $\cE_{ij}$ which are all uniformly bounded in $p$. Using equation \eqref{e:est4}, there exists a constant $C$ independent of $p$ such that:
\[
\left|
\left(
\frac{F(\xi_1)}{F(\xi_2)}-1
\right)
P(w_1,w_2)
\right|
\leq
C \,\, p^{-2/3}|\xi_1-\xi_2|.
\]
Therefore,
\[
\boxed{\Delta_1^{(-2)} \, p^{-2}
=\mathcal{O}\!\left(p^{-8/3} (\xi_1-\xi_2)\right)}
\]
\textbf{Analysis of the determinant $\Delta_2$.} Recall that we can write:
\[\Delta_2=\Delta_2^{(0)}+\Delta_2^{(-1)}p^{-1}+\Delta_2^{(-2)}p^{-2}\\\\\]
and the coefficient $\Delta_2^{(0)}$ is given by:
\begin{align*}
\Delta_2^{(0)} = 2 \pi i e^{i \pi/3}- 2 \pi i e^{i \pi/3} \frac{(1+\lambda_2)^{1/4} \left[q\left(\frac{\xi_1}{c_1p^{2/3}}\right)\right]^{1/4}}{(1+\lambda_1)^{1/4} \left[q\left(\frac{\xi_2}{c_1p^{2/3}}\right)\right]^{1/4}}.
\end{align*}
We can rewrite $\Delta_2^{(0)}$ as follows:
\begin{align*}
\Delta_2^{(0)}&=\frac{2 \pi i e^{i \pi/3}}{\left( 2+\frac{\xi_1}{c_1 p^{2/3}} \right)^{1/4}\left[q\left(\frac{\xi_2}{c_1p^{2/3}}\right)\right]^{1/4}} \times \left( \left( 2+\frac{\xi_1}{c_1 p^{2/3}} \right)^{1/4}\left[q\left(\frac{\xi_2}{c_1p^{2/3}}\right)\right]^{1/4} - \left( 2+\frac{\xi_2}{c_1 p^{2/3}} \right)^{1/4}\left[q\left(\frac{\xi_1}{c_1p^{2/3}}\right)\right]^{1/4}\right)
\end{align*}
and we can follow a similar analysis as for the coefficient $\Delta_1^{(0)}$ (simply interchanging the indices $1$ and $2$ in the $\xi$-variable) to show the following for $\Delta_2^{(0)}$:
\begin{align}\label{e:box2}
\boxed{\Delta_2^{(0)} = \left( \frac{(-1)^{5/6} \pi}{4 c_1 q_0} (2 q_1-q_0)\right) p^{-2/3}(\xi_2-\xi_1)+\mathcal{O}\left(p^{-4/3} \left( |\xi_1|+|\xi_2|\right) |\xi_1-\xi_2| \right) + \mathcal{O} \left( p^{-6/3} |\xi_1| |\xi_2||\xi_1-\xi_2|\right)}.
\end{align}
Next, the coefficient $\Delta_2^{(-1)}$ is given by:
\[
\begin{aligned}
((-1)^{5/6}\pi)^{-1}\Delta_2^{(-1)}
&=
\Bigl[
2\bigl(\cE_{11}(w_1)+\cE_{22}(w_1)\bigr)
-P(w_1,w_2)
\Bigr] 
-\left(
\frac{F(\xi_2)}{F(\xi_1)}-1
\right)P(w_1,w_2).
\end{aligned}
\]
where $F(\xi)$ as in equation \eqref{e:funF}, and $P(w_1,w_2)$ here is defined to be:
\[
\begin{aligned}
P(w_1,w_2)
&=
\cE_{11}(w_1)-\cE_{11}(w_2)
+\cE_{22}(w_1)-\cE_{22}(w_2) \\
&\quad
-i\bigl(\cE_{12}(w_2)-\cE_{12}(w_1)\bigr)
-i\bigl(\cE_{21}(w_1)-\cE_{21}(w_2)\bigr) \\
&\quad
+2\cE_{11}(w_2)+2\cE_{22}(w_2).
\end{aligned}
\]
Using a similar analysis as for the coefficient $\Delta_1^{(-1)}$, we can show that:
\[
\boxed{\Delta_2^{(-1)} \, p^{-1}
=\mathcal{O}\!\left(p^{-5/3} (\xi_1-\xi_2)\right)}
\]
Lastly, we estimate $\Delta_2^{(-2)}$, which is given by:
\begin{align*}
\left((-1)^{1/3}\pi\right)^{-1}\Delta_2^{(-2)}
&=\left[2i\left(
-\cE_{12}(w_1)\cE_{21}(w_1)
+\cE_{11}(w_1)\cE_{22}(w_1)
\right)
+P(w_1,w_2)\right]
+
\left(
\frac{F(\xi_2)}{F(\xi_1)}-1
\right)Q(w_1,w_2).
\end{align*}
Here $P(w_1,w_2)$ is:
\begin{align*}
P(w_1,w_2)
={}&
-\left(\cE_{11}(w_2)-i\cE_{12}(w_2)\right)
\left(\cE_{21}(w_1)+i\cE_{22}(w_1)\right)
+
\left(\cE_{11}(w_1)+i\cE_{12}(w_1)\right)
\left(\cE_{21}(w_2)-i\cE_{22}(w_2)\right).
\end{align*}
The term:
\[2i\left(
-\cE_{12}(w_1)\cE_{21}(w_1)
+\cE_{11}(w_1)\cE_{22}(w_1)
\right)
+P(w_1,w_2)\]
reduces into:
\begin{align*}
&-\left[
\cE_{11}(w_2)-\cE_{11}(w_1)
-i\left(\cE_{12}(w_2)-\cE_{12}(w_1)\right)
\right]
\left(\cE_{21}(w_1)+i\cE_{22}(w_1)\right)
\\
&\quad
+\left(\cE_{11}(w_1)+i\cE_{12}(w_1)\right)
\left[
\cE_{21}(w_2)-\cE_{21}(w_1)
-i\left(\cE_{22}(w_2)-\cE_{22}(w_1)\right)
\right].
\end{align*}
At this stage, using similar analysis as for the coefficient $\Delta_1^{(-2)} $ we can show that
\[
\boxed{\Delta_2^{(-2)} \, p^{-2}
=\mathcal{O}\!\left(p^{-8/3} (\xi_1-\xi_2)\right)}
\]
This concludes all the estimates for $\Delta_j$'s. 
\par
\textbf{Analysis of the first term in the product \eqref{Eq. Kernel hat Ytilde}.} Next, we analyze the first term in the product \eqref{Eq. Kernel hat Ytilde} which, for brevity, we denote by $\widehat{K}_{p}^{(f)}$:
\begin{equation}\label{Eq. DFF}
\widehat{K}_{p}^{(f)}(\xi_1, \xi_2) = \frac{1}{\sqrt{\chi_p'(\chi_p^{-1}(\xi_1))\,\chi_p'(\chi_p^{-1}(\xi_2))}\,\left[\chi_p^{-1}(\xi_1)-\chi_p^{-1}(\xi_2)\right]}.
\end{equation}
From equation \eqref{eq:eq515}, we have:
\begin{align*}
   \frac{d}{d\xi} \lambda(\xi)&=\frac{a_1}{c_1\,p^{2/3}}+\frac{1}{\xi}\sum_{j=2}^{\infty}j\,a_j\left(\frac{\xi}{c_1\,p^{2/3}}\right)^j =\frac{a_1}{c_1\,p^{2/3}}+\mathcal{O}\left(p^{-4/3}\right), 
\end{align*}
for $|\xi|<c_1\, p^{2/3}$. Then, the first factor in the denominator of \eqref{Eq. DFF} becomes:
\begin{align*}
  \sqrt{\chi_p'(\chi_p^{-1}(\xi_1))\,\chi_p'(\chi_p^{-1}(\xi_2))} & = \left(\sqrt{\left(\frac{a_1}{c_1\,p^{2/3}}+\frac{1}{\xi_1}\sum_{j=2}^{\infty}j\,a_j\left(\frac{\xi_1}{c_1\,p^{2/3}}\right)^j\right)\left(\frac{a_1}{c_1\,p^{2/3}}+\frac{1}{\xi_2}\sum_{j=2}^{\infty}j\,a_j\left(\frac{\xi_2}{c_1\,p^{2/3}}\right)^j\right)}\right)^{-1}\\
  &=\frac{c_1\,p^{2/3}}{a_1}\left(\sqrt{\left(1+\frac{c_1\,p^{2/3}}{a_1\xi_1}\sum_{j=2}^{\infty}j\,a_j\left(\frac{\xi_1}{c_1\,p^{2/3}}\right)^j\right)\left(1+\frac{c_1\,p^{2/3}}{a_1\xi_2}\sum_{j=2}^{\infty}j\,a_j\left(\frac{\xi_2}{c_1\,p^{2/3}}\right)^j\right)}\right)^{-1}
\end{align*}
while the second factor in the denominator of \eqref{Eq. DFF} becomes:
\begin{align*}
\chi_p^{-1}(\xi_1)-\chi_p^{-1}(\xi_2)&=\sum_{j=0}^{\infty}a_j\left(\frac{\xi_1}{c_1\,p^{2/3}}\right)^j-\sum_{j=0}^{\infty}a_j\left(\frac{\xi_2}{c_1\,p^{2/3}}\right)^j\\
      &=\frac{a_1}{c_1\,p^{2/3}}(\xi_1-\xi_2)+\sum_{j=2}^{\infty} \frac{a_j}{(c_1\,p^{2/3})^j}(\xi_1^j-\xi_2^j)\\
      &=\frac{a_1}{c_1\,p^{2/3}}(\xi_1-\xi_2)\left[ 1+\sum_{j=2}^{\infty} \frac{a_j\,c_1\,p^{2/3}}{a_1(\xi_1-\xi_2)(c_1\,p^{2/3})^j}(\xi_1^j-\xi_2^j)\right].
\end{align*}
Then, equation \eqref{Eq. DFF} becomes:
\begin{align}\label{eq First and second factor}
\widehat{K}_{p}^{(f)}(\xi_1, \xi_2) & = \frac{1}{\xi_1-\xi_2}\frac{\sqrt{\left(1+\frac{c_1\,p^{2/3}}{a_1\xi_1}\sum_{j=2}^{\infty}j\,a_j\left(\frac{\xi_1}{c_1\,p^{2/3}}\right)^j\right)\left(1+\frac{c_1\,p^{2/3}}{a_1\xi_2}\sum_{j=2}^{\infty}j\,a_j\left(\frac{\xi_2}{c_1\,p^{2/3}}\right)^j\right)}}{ 1+\sum_{j=2}^{\infty} \frac{a_j\,c_1\,p^{2/3}}{a_1(\xi_1-\xi_2)(c_1\,p^{2/3})^j}(\xi_1^j-\xi_2^j)} \nonumber \\
         &=\frac{1}{\xi_1-\xi_2}[1+R(\xi_1,\xi_2)],
\end{align}
where we have introduced the notation
\begin{align}\label{e:estR}
R(\xi_1,\xi_2)&=-1+\frac{\sqrt{\left(1+\frac{c_1\,p^{2/3}}{a_1\xi_1}\sum_{j=2}^{\infty}j\,a_j\left(\frac{\xi_1}{c_1\,p^{2/3}}\right)^j\right)\left(1+\frac{c_1\,p^{2/3}}{a_1\xi_2}\sum_{j=2}^{\infty}j\,a_j\left(\frac{\xi_2}{c_1\,p^{2/3}}\right)^j\right)}}{ 1+\sum_{j=2}^{\infty} \frac{a_j\,c_1\,p^{2/3}}{a_1(\xi_1-\xi_2)(c_1\,p^{2/3})^j}(\xi_1^j-\xi_2^j)} \nonumber \\&=\mathcal{O}\left(p^{-4/3}\right).
\end{align}
The implicit constant in $\mathcal{O}\left(p^{-4/3}\right)$ is uniform for all $\xi_1, \xi_2\in\mathcal{U}_{1}$ (it depends on the radius of the disk).

Before we put everything together, we simplify the notation in \eqref{eq:Kps}. To this end we introduce the following notation:
\begin{equation*}
K_{Airy}^{(U)}:=\mathrm{Ai}(\xi_1)\mathrm{Ai}'(\xi_2)-\mathrm{Ai}(\xi_2)\mathrm{Ai}'(\xi_1), 
\end{equation*}

\begin{align*}
\det[\mathcal{B}(\xi_1)] &= -2 \pi i \, e^{i \pi/3} \left( 1 + \frac{1}{p} \left( \cE_{11}(\xi_1) + \cE_{22}(\xi_1)\right) + \frac{1}{p^2} \, \det[\mathcal{E}(\xi_1)]\right)= -2 \pi i \, e^{i \pi/3} \left( 1 + \frac{1}{p} \mathbb{D}\right),
\end{align*} 
where 
\begin{equation}\label{eq. D O(1)}
\mathbb{D}=\mathbb{D}(\xi_1):= \cE_{11}(\xi_1) + \cE_{22}(\xi_1) + \frac{1}{p} \, \det[\mathcal{E}(\xi_1)]=\mathcal{O}(1),
\end{equation}
and 
\begin{equation}\label{eq. def boxes}
\begin{aligned}
\Box_1 &:= \mathrm{Ai}(\xi_1)\mathrm{Ai}'(\xi_2),
&\qquad
\Box_2 &:= \mathrm{Ai}(\xi_2)\mathrm{Ai}'(\xi_1),\\
\Box_3 &:= \mathrm{Ai}(\xi_1)\mathrm{Ai}(\xi_2),
&\qquad
\Box_4 &:= \mathrm{Ai}'(\xi_2)\mathrm{Ai}'(\xi_1).
\end{aligned}
\end{equation}
Using this we can write:
\begin{equation}\label{eq:Kps boxes}
  \widehat{K}_{p}^{(s)}(\xi_1, \xi_2)=-2 \pi i \, e^{i \pi/3}\left(1+\frac{1}{p}\mathbb{D}\right)K_{Airy}^{(U)}+
\Delta_1\Box_{1}- \Delta_2 \Box_{2} +\Delta_3\Box_{3} + \Delta_4 \Box_{4}.
\end{equation}
Hence, replacing \eqref{eq:Kps boxes} and \eqref{eq First and second factor} in \eqref{Eq. Kernel hat Ytilde} we obtain:
\begin{equation}\label{Eq. Kernel hat Ytilde box}
\begin{split}
\widehat{K}_p (\xi_1, \xi_2) 
&=-\frac{e^{-\frac{i\pi}{3}}}{2 \pi i}\frac{1+R(\xi_1,\xi_2)}{\xi_1-\xi_2}\left\{-2 \pi i \, e^{i \pi/3}\left(1+\frac{1}{p}\mathbb{D}\right)K_{Airy}^{(U)}+
\Delta_1\Box_{1}- \Delta_2 \Box_{2} +\Delta_3\Box_{3} + \Delta_4 \Box_{4} \right\}\\
&= K_{Airy}(\xi_1,\xi_2)-\frac{e^{-\frac{i\pi}{3}}}{2 \pi i}\frac{1}{\xi_1-\xi_2}\left(\Delta_1^{(00)}\,\Box_1-\Delta_2^{(00)}\,\Box_2\right)+\widetilde{R}(\xi_1,\xi_2),
\end{split}
\end{equation}
where
\begin{equation}\label{eq. R tilde full}
\begin{split}
\widetilde{R}(\xi_1,\xi_2)&=\frac{\mathbb{D}}{p}K_{Airy}(\xi_1,\xi_2)-\frac{e^{-\frac{i\pi}{3}}}{2 \pi i}\frac{1}{\xi_1-\xi_2}\left\{\left(\Delta_1-\Delta_1^{(00)}\right)\,\Box_1-\left(\Delta_2-\Delta_2^{(00)}\right)\,\Box_2+\Delta_3\Box_{3} + \Delta_4 \Box_{4}\right\}\\
&-\frac{e^{-\frac{i\pi}{3}}}{2 \pi i}\frac{R(\xi_1,\xi_2)}{\xi_1-\xi_2}\left\{-2 \pi i \, e^{i \pi/3}\left(1+\frac{1}{p}\mathbb{D}\right)K_{Airy}^{(U)}+
\Delta_1\,\Box_{1}- \Delta_2\,\Box_{2} +\Delta_3\,\Box_{3} + \Delta_4\, \Box_{4} \right\},
\end{split}
\end{equation}
and
\begin{equation*}
    \Delta_1^{(00)} = \left( \frac{(-1)^{5/6} \pi}{4 c_1 q_0} (2 q_1-q_0)\right) p^{-2/3}(\xi_1-\xi_2),\qquad  \Delta_2^{(00)} = \left( \frac{(-1)^{5/6} \pi}{4 c_1 q_0} (2 q_1-q_0)\right) p^{-2/3}(\xi_2-\xi_1).
\end{equation*}
Then, note that
\begin{align*}
\frac{e^{-\frac{i\pi}{3}}}{2 \pi i}\frac{1}{\xi_1-\xi_2}\left(\Delta_1^{(00)}\,\Box_1-\Delta_2^{(00)}\,\Box_2\right) &=-\frac{1}{2}   \frac{(2 q_1-q_0) }{4 c_1 q_0}\left( \Box_1+\Box_2 \right) \, p^{-2/3} \\
&=\frac{1}{4 c_1}\left( \frac{1}{2} +\frac{2}{5}\frac{h_p'(1)}{h_p(1)} \right)\,\left( \Box_1+\Box_2 \right) \, p^{-2/3} .
\end{align*}
Therefore, we can write 
\begin{equation}\label{asym exp}
\widehat{K}_p (\xi_1, \xi_2) = K_{\Airy}(\xi_1,\xi_2) + \widehat{C}_p \left( \mathrm{Ai}(\xi_1) \mathrm{Ai}'(\xi_2) + \mathrm{Ai}(\xi_2) \mathrm{Ai}'(\xi_1) \right) \, p^{-2/3} + \widetilde{R}(\xi_1,\xi_2),
\end{equation}
with
\begin{equation*}
    \widehat{C}_p :=\frac{1}{4 c_1}\left( \frac{1}{2} +\frac{2}{5}\frac{h_p'(1)}{h_p(1)} \right).
\end{equation*}
Since the error matrix $E$ is analytic in $\mathcal{U}_{1}$, and the domain $\chi_{p}(\mathcal{U}_{1})$ expands to cover the entire complex plane as $p \to \infty$, relation~\eqref{asym exp} holds for any fixed $\xi_{1},\xi_{2} \in \mathbb{C}$, for all $p$ sufficiently large.
\vspace{0.2cm}

\textbf{Uniform bound if the error term $\widetilde{R}$.} Consider a fixed constant $L>0$. We will use the standard bounds on the Airy function and its derivative to prove the existence of constants $C,k>0$ and an integer $M>0$, depending only on $L$, such that, for all $p\ge M$ and all $\xi_1,\xi_2>-L$:
\begin{equation}\label{eq. bound Rtilde 4 cases}
\left|\widetilde{R}(\xi_1,\xi_2)\right|
\le
\begin{cases}
Cp^{-1}e^{-k(\xi_1+\xi_2)},
& (\xi_1,\xi_2)\in(0,\infty)^2,\\[1mm]
Cp^{-1}e^{-k\xi_1},
& (\xi_1,\xi_2)\in(0,\infty)\times(-L,0],\\[1mm]
Cp^{-1}e^{-k\xi_2},
& (\xi_1,\xi_2)\in(-L,0]\times(0,\infty),\\[1mm]
Cp^{-1},
& (\xi_1,\xi_2)\in(-L,0]^2.
\end{cases}
\end{equation}
To prove this, the triangle inequality in~\eqref{eq. R tilde full} yields
\begin{equation}\label{eq. Abs R tilde terms}
\begin{split}
\left|\widetilde{R}(\xi_1,\xi_2)\right|&\le \frac{|\mathbb{D}|}{p}\Bigl|K_{Airy}(\xi_1,\xi_2)\Bigr|\\
&+\frac{1}{2 \pi}\left\{\left|\frac{\Delta_1-\Delta_1^{(00)}}{\xi_1-\xi_2}\right|\,|\,\Box_1|+\left|\frac{\Delta_2-\Delta_2^{(00)}}{\xi_1-\xi_2}\right|\,|\,\Box_2|+\left|\frac{\Delta_3}{\xi_1-\xi_2}\right| |\,\Box_{3}| + \left|\frac{\Delta_4}{\xi_1-\xi_2}\right| |\,\Box_{4}|\right\}\\
&+\frac{1}{2 \pi }\Bigl|R(\xi_1,\xi_2)\Bigr|\left\{2 \pi \left|1+\frac{\mathbb{D}}{p}\right| \Bigl|K_{Airy}\Bigr|+
\left|\frac{\Delta_1}{\xi_1-\xi_2}\right|\,|\,\Box_{1}|+ \left|\frac{\Delta_2}{\xi_1-\xi_2}\right|\,|\,\Box_{2}| +\left|\frac{\Delta_3}{\xi_1-\xi_2}\right| |\,\Box_{3}| + \left|\frac{\Delta_4}{\xi_1-\xi_2}\right| |\,\Box_{4}| \right\}.
\end{split}
\end{equation}
Here, recall the relations~\eqref{eq. D O(1)} and~\eqref{e:estR}:
\begin{equation}\label{eq. D and R order}
\mathbb{D}=\mathcal{O}(1)\quad \text{and} \quad R=\mathcal{O}\left(p^{-4/3}\right).
\end{equation}
Also, the estimates on the terms $\Delta_i$, $i=1,2,3,4$, yield
\begin{equation}\label{eq. order delta12 lower}
\left|\frac{\Delta_1-\Delta_1^{(00)}}{\xi_1-\xi_2}\right|=\mathcal{O}\left(p^{-4/3}\right),\quad \left|\frac{\Delta_2-\Delta_2^{(00)}}{\xi_1-\xi_2}\right|=\mathcal{O}\left(p^{-4/3}\right), \qquad
\end{equation}
and relations \eqref{eq. bound delta3}, \eqref{eq. bound delta 4} give:
\begin{equation}\label{eq. order deltas final}
    \left|\frac{\Delta_i}{\xi_1-\xi_2}\right|=
    \begin{cases}
        \mathcal{O}\left(1\right),& i=1\\
        \mathcal{O}\left(1\right),& i=2\\
        \mathcal{O}\left(p^{-4/3}\right),& i=3\\
        \mathcal{O}\left(p^{-2}\right),& i=4.\\
        \end{cases}
\end{equation}
It is important to notice that the implicit constants in notation $\mathcal{O}\left( \,\cdot\,\right)$, in \eqref{eq. D and R order} and \eqref{eq. order deltas final} depend only on the radius of the disks $\mathcal{U}_{-1}$, $\mathcal{U}_{1}$, , and on the remaining part of the contour $\Sigma_E$ shown in Figure~\ref{fig Contour E}.

On the other hand, the standard Airy bounds imply that there exists a constant $C>0$, such that for every $\xi\in\mathbb{R}$,
\begin{equation}\label{eq. Airy sharp bound}
|\mathrm{Ai}(\xi)| \leq \frac{C}{(1+|\xi|)^{1/4}}
\begin{cases}
e^{-\frac23\xi^{3/2}}, & \xi\ge 0,\\
1, & \xi<0.
\end{cases}
\end{equation}
and 
\begin{equation}\label{eq. Airy derivative sharp bound}
|\mathrm{Ai}'(\xi)| \leq C(1+|\xi|)^{1/4}
\begin{cases}
e^{-\frac23\xi^{3/2}}, & \xi\ge 0,\\
1, & \xi<0,
\end{cases}
\end{equation}
were $C$ is an absolute constant independent of $\xi$. Moreover, the bounds~\eqref{eq. Airy sharp bound} and~\eqref{eq. Airy derivative sharp bound} imply that there exist constants $C,k>0$ such that, for $\xi>0$,
\begin{equation}\label{eq. bounds Ai Ai' positive}
|\mathrm{Ai}(\xi)|
\le Ce^{-k\xi},
\qquad
|\mathrm{Ai}'(\xi)|
\le Ce^{-k\xi}.
\end{equation}
Here, note that the constant $k$ can be taken as $k=2/3$ and the constant $C$ depends only on $k$.

We will use the integral representation of the Airy kernel:
\[
K_{\mathrm{Airy}}(x,y)= \int_0^\infty
\mathrm{Ai}(x+t)\,\mathrm{Ai}(y+t)\,dt,
\]
to prove that the term 
\[\frac{|\mathbb{D}|}{p}\Bigl|K_{Airy}(\xi_1,\xi_2)\Bigr|
\]
in \eqref{eq. bound Rtilde 4 cases}, satisfies the same estimates as $\tilde{R}$ in~\eqref{eq. bound Rtilde 4 cases}, with the constant $C$ depending only on the values $L$ and $k$ ($k=2/3$). Indeed, we must consider the following cases:
\begin{itemize}
\item For $(\xi_1,\xi_2)\in(0,\infty)^2$, the bound~\eqref{eq. bounds Ai Ai' positive} yields
\begin{align}
|K_{\Airy}(\xi_1,\xi_2)|&= \int_0^\infty
\mathrm{Ai}(\xi_1+t)\,\mathrm{Ai}(\xi_2+t)\,dt\\
&\le C^2\int_0^\infty e^{-k(\xi_1+t)} e^{-k(\xi_2+t)}\,dt\\
&= C^2 e^{-k(\xi_1+\xi_2)} \int_0^\infty e^{-2kt}\,dt\\
&= \frac{C^2}{2k}e^{-k(\xi_1+\xi_2)}\\
&\le C e^{-k(\xi_1+\xi_2)},
\end{align}
where in the last inequality the constant $C$ has been redefined. 

\item For $(\xi_1,\xi_2)\in(0,\infty)\times(-L,0]$, since the airy function is bounded on \([-L,\infty)\), ~\eqref{eq. bounds Ai Ai' positive} gives
\begin{align}
|K_{\Airy}(\xi_1,\xi_2)| &\le C e^{-k\xi_1}
\int_0^\infty e^{-kt}|\mathrm{Ai}(\xi_2+t)|\,dt\\
&\le Ce^{-k\xi_1},
\end{align}
where again, $C$ has been redefined in the last inequality.

\item The same argument from the previous case implies that for $(\xi_1,\xi_2)\in(-L,0]\times(0,\infty)$, we have
\[ |K_{\Airy}(\xi_1,\xi_2)| \le Ce^{-k\xi_2}. \]

\item Finally, for $(\xi_1,\xi_2)\in [-L,0]^2$, since this is a compact set we have 
\[ |K_{\Airy}(\xi_1,\xi_2)|\le C_L.\]
\end{itemize}
Therefore, after changing the constants if necessary, and since $\mathbb{D}=\mathcal{O}(1)$, we conclude that $\frac{|\mathbb{D}|}{p}\Bigl|K_{Airy}(\xi_1,\xi_2)\Bigr|$ satisfies the estimates in~\eqref{eq. bound Rtilde 4 cases}. Note that, in the last case the constant $C_L$ depends on size of the rectangle $[-L,0]^2$, so to establish~\eqref{eq. bound Rtilde 4 cases} using the same constant $C$ in all regions, the dependence of $C$ will be only on the value $L$ and the decay rate $k$ (which can be taken as $k=2/3$).

Now, we claim that all the terms $\Box_i$, $i=1,2,3,4$,  defined in ~\eqref{eq. def boxes}, have an upper bound as $\tilde{R}$ in~\eqref{eq. bound Rtilde 4 cases}, but without the factor $1/p$. The other factors in~\eqref{eq. Abs R tilde terms} will compensate with a faster decay. We only prove this bound for the term $\Box_1$, since the argument is similar in the other cases.

Recall that 
\[
\Box_1=\mathrm{Ai}(\xi_1)\mathrm{Ai}'(\xi_2).
\]
Hence, we have the following cases:
\begin{itemize}
    \item For $(\xi_1,\xi_2)\in(0,\infty)^2$, 
     \[ |\mathrm{Ai}(\xi_1)\mathrm{Ai}'(\xi_2)| \le C e^{-k(\xi_1+\xi_2)}. \]
    \item For $(\xi_1,\xi_2)\in(0,\infty)\times(-L,0]$, \[|\mathrm{Ai}(\xi_1)\mathrm{Ai}'(\xi_2)|\le Ce^{-k\xi_1}C_L \le C_L e^{-k\xi_1}. \]
    \item For  $(\xi_1,\xi_2)\in(-L,0]\times(0,\infty)$, \[ |\mathrm{Ai}(\xi_1)\mathrm{Ai}'(\xi_2)| \le C e^{-k\xi_2}. \]
    \item Finally, for $(\xi_1,\xi_2)\in [-L,0]^2$, \[ |\mathrm{Ai}(\xi_1)\mathrm{Ai}'(\xi_2)| \le C_L. \]
\end{itemize}
Again, after changing the constants if necessary, this proves our claim holds for $\Box_1$. As above, to establish~\eqref{eq. bound Rtilde 4 cases} using the same constant $C$ in all regions, $C$ will depend only on $L$ and $k$.

Finally, note that the estimates \eqref{eq. D and R order}, \eqref{eq. order delta12 lower} and~\eqref{eq. order deltas final}, together with the estimates on
the terms $\Box_i$, $i=1,2,3,4$, imply that the expressions
\begin{equation}
    \frac{1}{2 \pi}\left\{\left|\frac{\Delta_1-\Delta_1^{(00)}}{\xi_1-\xi_2}\right|\,|\,\Box_1|+\left|\frac{\Delta_2-\Delta_2^{(00)}}{\xi_1-\xi_2}\right|\,|\,\Box_2|+\left|\frac{\Delta_3}{\xi_1-\xi_2}\right| |\,\Box_{3}| + \left|\frac{\Delta_4}{\xi_1-\xi_2}\right| |\,\Box_{4}|\right\}
\end{equation}
and
\begin{equation}
\frac{1}{2 \pi }\Bigl|R(\xi_1,\xi_2)\Bigr|\left\{2 \pi \left|1+\frac{\mathbb{D}}{p}\right| \Bigl|K_{Airy}\Bigr|+
\left|\frac{\Delta_1}{\xi_1-\xi_2}\right|\,|\,\Box_{1}|+ \left|\frac{\Delta_2}{\xi_1-\xi_2}\right|\,|\,\Box_{2}| +\left|\frac{\Delta_3}{\xi_1-\xi_2}\right| |\,\Box_{3}| + \left|\frac{\Delta_4}{\xi_1-\xi_2}\right| |\,\Box_{4}| \right\},
\end{equation}
in \eqref{eq. bound Rtilde 4 cases}, satisfy the same estimates as $\tilde{R}$ in~\eqref{eq. bound Rtilde 4 cases}, but with the smaller prefactor $p^{-4/3}$ instead of $p^{-1}$.

Therefore $\widetilde{R}$ satisfies \eqref{eq. bound Rtilde 4 cases}, which completes the proof of Proposition \ref{proposition Khat to K_airy}.

\subsection{Proof of Proposition~\ref{prop prop53}}\label{proof Prop trace D0}
Let us denote 
\[
\mathcal{P} := (\mathbb{I}-\mathcal{K}_{\Airy})^{-1}.
\]
First, we consider the action of the operator $\mathcal{D}_0$, whose kernel is defined as
\begin{align*}
D_0(x,y)
= \widehat{C}_p\,p^{-2/3}
\bigl(\mathrm{Ai}(x)\mathrm{Ai}'(y)
+\mathrm{Ai}'(x)\mathrm{Ai}(y)\bigr).
\end{align*}
Let $f\in L^2([s,\infty))$, we have:
\begin{align*}
\left( \mathcal{D}_0 f \right) (x) &= \int_s^{\infty}  \widehat{C}_p\, p^{-2/3}  \left( \mathrm{Ai}(x) \mathrm{Ai}'(y) + \mathrm{Ai}(y) \mathrm{Ai}'(x) \right) f(y) \, dy\\
&= p^{-2/3} \, \widehat{C}_p  \left( \mathrm{Ai}(x) \int_s^{\infty} \mathrm{Ai}'(y) f(y) \, dy + \mathrm{Ai}'(x)\int_s^{\infty} \mathrm{Ai}(y) f(y) \, dy \right).
\end{align*}
Applying $\mathcal{P}$ to both sides and using its linearity gives
\begin{align*}
\left( \left( \mathbb{I}-\mathcal{K}_{Airy} \right)^{-1} \mathcal{D}_0 f \right)& (x) = p^{-2/3} \,\widehat{C}_p \left [ \mathcal{P} \left(  \mathrm{Ai}(x) \int_s^{\infty} \mathrm{Ai}'(y) f(y) \, dy\right) + \mathcal{P} \left( \mathrm{Ai}'(x)\int_s^{\infty} \mathrm{Ai}(y) f(y) \, dy \right)\right]\\
&= p^{-2/3} \, \widehat{C}_p  \left[ \left( \int_s^{\infty} \mathrm{Ai}'(y) f(y) \, dy \right) \left( \mathcal{P}  \mathrm{Ai}\right)(x) + \left( \int_s^{\infty} \mathrm{Ai}(y) f(y) \, dy \right) \left( \mathcal{P}  \mathrm{Ai}'\right)(x) \right]\\
&= p^{-2/3} \, \widehat{C}_p  \left[ \left( \int_s^{\infty} \mathrm{Ai}'(y) f(y) \, dy \right) q(x,s) + \left( \int_s^{\infty} \mathrm{Ai}(y) f(y) \, dy \right) p(x,s) \right]\\
&= p^{-2/3} \, \widehat{C}_p  \int_s^{\infty} \left( \mathrm{Ai}'(y) q(x,s) + \mathrm{Ai}(y) p(x,s) \right) f(y) \, dy,
\end{align*}
where recall that:
\begin{equation}
\begin{aligned}
q(x,s) =
\left(\mathcal{P}
\mathrm{Ai}\right)(x), \quad p(x,s) =
\left(\mathcal{P}\mathrm{Ai}'\right)(x).
\end{aligned}
\end{equation}
Therefore, the kernel of $(\mathbb I-\mathcal K_{\Airy})^{-1}\mathcal D_0$ equals:
\[ p^{-2/3} \, \widehat{C}_p\left( \mathrm{Ai}'(y) q(x,s) + \mathrm{Ai}(y) p(x,s) \right).\]
Hence, the trace of the operator $\left( \mathbb{I}-  \mathcal K_{\Airy}\right)^{-1} \mathcal{D}_o$ is given by
\begin{align*}
  \mathrm{Tr} \left( \left( \mathbb{I}-\mathcal{K}_{\Airy} \right)^{-1} \mathcal{D}_0 \right) &=  p^{-2/3} \, \widehat{C}_p \int_s^{\infty} \mathrm{Ai}'(y) \, q(y,s) + \mathrm{Ai}(y) \, p(y,s) \, dy =2\, p^{-2/3} \, \widehat{C}_p\, v(s).
\end{align*}

\subsection{Proof of Proposition~\ref{lem:bound-D-D0}}\label{proof Prop Fredholm D to D0}
Fix $L>0$ and let $L_0>-L$. Our goal is to prove the relation
\begin{equation}\label{eq. trunc}
     \det \left( \mathbb{I}- \mathcal{P}\mathcal{D}\right)_{L^2([L_0,\infty))}=\det \left( \mathbb{I}-\mathcal{P} \mathcal{D} _0\right)_{L^2([L_0,\infty))}+\mathcal{O}\left(p^{-1}\right),
\end{equation}
with the implicit constant in $\mathcal O(p^{-1})$ independent of $L_0$.

We will use the Hilbert--Carleman determinant of the trace class operators $\mathcal{P} \mathcal{D}$ and $\mathcal{P} \mathcal{D}_0$, defined by:
\begin{align*}
{\det}_2 (\mathbb{I}-\mathcal{P} \mathcal{D}) := \det (\mathbb{I}-\mathcal{P} \mathcal{D})_{L^2([L_0,\infty))} e^{\mathrm{Tr(\mathcal{P} \mathcal{D})}}, \quad {\det}_2 (\mathbb{I}-\mathcal{P} \mathcal{D}_0) := \det (\mathbb{I}-\mathcal{P} \mathcal{D}_0)_{L^2([L_0,\infty))} e^{\mathrm{Tr(\mathcal{P} \mathcal{D}_0)}}.
\end{align*}
To prove \eqref{eq. trunc}, it suffices to establish the following two estimates:
\begin{equation}\label{eq. key 1 Carleman}
    {\det}_2 (\mathbb{I}-\mathcal{P} \mathcal{D}) ={\det}_2 (\mathbb{I}-\mathcal{P} \mathcal{D}_0) +\mathcal{O}\left(p^{-1}\right)
\end{equation}
\begin{equation}\label{eq. key 2 Carleman}
    \Tr(\mathcal{P} \mathcal{D}) = \Tr(\mathcal{P} \mathcal{D}_0)+\mathcal{O}\left(p^{-1}\right).
\end{equation}
Indeed,  assuming that these estimates hold uniformly for $L_0>-L$. Then
\begin{align*}
\det \left( \mathbb{I}- \mathcal{P}\mathcal{D}\right)_{L^2([L_0,\infty))}&= {\det}_2 (\mathbb{I}-\mathcal{P} \mathcal{D})\, e^{\mathrm{-Tr(\mathcal{P} \mathcal{D})}}\\
&=\left({\det}_2 (\mathbb{I}-\mathcal{P} \mathcal{D}_0) +\mathcal{O}\left(p^{-1}\right)\right) \,e^{-\left( \Tr(\mathcal{P} \mathcal{D}_0)+\mathcal{O}\left(p^{-1}\right) \right)}\\
&={\det}_2 (\mathbb{I}-\mathcal{P} \mathcal{D}_0)\,e^{-\Tr(\mathcal{P} \mathcal{D}_0)}+\mathcal{O}\left(p^{-1}\right)\\
&=\det \left( \mathbb{I}- \mathcal{P}\mathcal{D}_0\right)_{L^2([L_0,\infty))}+\mathcal{O}\left(p^{-1}\right).
\end{align*}
Hence, it remains to prove the estimates \eqref{eq. key 1 Carleman} and \eqref{eq. key 2 Carleman}. We first note that
\begin{equation}\label{eq. def D0 in proof}
D_0(x,y) = \widehat{C}_p\,p^{-2/3}
\bigl( \mathrm{Ai}(x)\mathrm{Ai}'(y) +
\mathrm{Ai}'(x)\mathrm{Ai}(y) \bigr),
\end{equation}
where $\widehat C_p$ is uniformly bounded (see~\eqref{eq:errt}). Moreover,  recall that the operators $\mathcal{P}$ and $\mathcal{D}_0$ act on $L^2([L_0,\infty))$, with $L_0>-L$. Then 
\begin{equation}\label{eq. uniform bound P, D0}
\|\mathcal P\|_{op} =\mathcal O(1) \qquad \text{and}\qquad \lVert\mathcal{D}_0\rVert_2=\mathcal{O}\left(p^{-2/3}\right),
\end{equation}
uniformly for $L_0>-L$, with the implicit constant in $\mathcal O(\,\cdot\,)$ depending on the fixed value $L$.
\vspace{0.2cm}

\textbf{Proof of the estimate~\eqref{eq. key 1 Carleman}}. By Theorem~6.5 in~\cite{simon1977notes}:
\begin{equation}\label{eq. Thm Gohberg}
\left|{\det}_2 (\mathbb{I}-\mathcal{P} \mathcal{D})-{\det}_2 (\mathbb{I}-\mathcal{P} \mathcal{D}_0)\right|
\leq 
\lVert \mathcal{P}\mathcal{D}-\mathcal{P}\mathcal{D}_0\rVert_2 \exp\left( \frac{1}{2}(\lVert \mathcal{P}\mathcal{D}\rVert_2+\lVert \mathcal{P}\mathcal{D}_0\rVert_2 +1)^2\right).
\end{equation}
For the second term in the exponential factor, Eq.~\eqref{eq. uniform bound P, D0} and the ideal property of the Hilbert--Schmidt norm yield
\begin{equation}\label{eq. order PD_0}
   \lVert\mathcal{P}\mathcal{D}_0\rVert_2 \le \lVert\mathcal{P}\rVert_{op} \lVert\mathcal{D}_0\rVert_2 =\mathcal{O}\left(p^{-2/3}\right),
\end{equation}
with the implicit constant in $\mathcal O(p^{-2/3})$ independent of $L_0$, for all $L_0>-L$.

For the first term in the exponential factor, we have
\begin{align}\label{eq. bound HS PD}
   \lVert\mathcal{P}\mathcal{D}\rVert_2 &\le \lVert\mathcal{P}\rVert_{op}\, \lVert\mathcal{D}\rVert_2\nonumber \\
   & = \lVert\mathcal{P}\rVert_{op} \,\lVert \widehat{\mathcal{K}}_p^{(E)} - \mathcal{K}_{Airy}\rVert_2 \nonumber\\
   & =\lVert\mathcal{P}\rVert_{op} \,\lVert \Pi_p\, \left(\mathcal{K}_{Airy}+ \mathcal{D}_0\right) \,\Pi_p +\widetilde{\mathcal R}^{(E)}- \mathcal{K}_{Airy}\rVert_2\nonumber\\
   & \le\lVert\mathcal{P}\rVert_{op} \left(
   \|\Pi_p\mathcal K_{\mathrm{Airy}}\Pi_p -\mathcal K_{\mathrm{Airy}}\|_2 + \|\Pi_p\,\mathcal D_0\,\Pi_p \|_2 + \| \widetilde{\mathcal R}^{(E)}\|_2\right).
\end{align}
Here, \(\widetilde{\mathcal R}^{(E)}\) denotes the integral operator, on \(L^2([L_0,\infty))\), whose kernel is obtained by extending \(\widetilde R\) from \([L_0,L_p]^2\) to \([L_0,\infty)^2\), setting it equal to zero outside its original domain. Also, recall the definition of $L_p$ in~\eqref{eq. def L_p} and the expression for $\chi_p$ in~\eqref{eq. chi_p asymtotics}: 
\begin{equation}\label{eq. def asymp L_p}
L_p:=\chi_p(\lambda_D)=p^{\frac{2}{3}}c_1\Omega(\lambda_D)>0.
\end{equation}
Here,  $\Omega(\lambda_D)$ and $c_1=c_1(p)=\left( -h_p(1)\right)^{2/3}$ are positive quantities, uniformly bounded away from zero (see~\eqref{eq. def Omega}, \eqref{eq. def constant c1} and~\eqref{eq:consh1}).

Now, consider the projection operator  $\Pi_p$, defined by
\[
(\Pi_p f)(x)=\mathds{1}_{[L_0,L_p]}(x)f(x),
\]
and recall the integral representation of the Airy kernel:
\[
K_{\mathrm{Airy}}(x,y)
=
\int_0^\infty
\mathrm{Ai}(x+t)\,\mathrm{Ai}(y+t)\,dt.
\]
Moreover, consider the integral operator 
\[
\mathcal A:L^2([0,\infty))\longrightarrow L^2([L_0,\infty)),
\]
with kernel $A(x,t)=\mathrm{Ai}(x+t)$. Then, the Airy operator can be factored as 
\[
\mathcal K_{\mathrm{Airy}}=\mathcal A\,\mathcal A^*.
\]
Hence, we can write
\begin{align*}
\mathcal K_{\mathrm{Airy}} - \Pi_p\mathcal K_{\mathrm{Airy}}\Pi_p
&= \mathcal A\mathcal A^*- \Pi_p \mathcal A\mathcal A^* \Pi_p\\
&=(\mathbb{I}-\Pi_p)\mathcal A\mathcal A^*+
\Pi_p\mathcal A\mathcal A^*(\mathbb{I}-\Pi_p).
\end{align*}
By the triangle inequality and the ideal property of the Hilbert--Schmidt norm, we have
\begin{equation}\label{eq. trace norm bd1}
\begin{split}
\left\|\mathcal K_{\mathrm{Airy}}-\Pi_p\mathcal K_{\mathrm{Airy}}\Pi_p\right\|_2
&\leq \lVert (\mathbb{I}-\Pi_p)\mathcal A\rVert_2 \lVert\mathcal A^*\rVert_{op}+\lVert \Pi_p\mathcal A\rVert_{op} \lVert\mathcal A^*(\mathbb{I}-\Pi_p)\rVert_2\\
&\leq
2\lVert\mathcal A\rVert_{op} \, \lVert (\mathbb{I}-\Pi_p)\mathcal A\rVert_2.
\end{split} 
\end{equation}
Here,
\begin{align*}
\lVert (\mathbb{I}-\Pi_p)\mathcal A\rVert_2^2
&=
\int_{L_0}^{\infty}\int_0^\infty
\left|
\mathbf{1}_{(L_p,\infty)}(x)
\mathrm{Ai}(x+t)
\right|^2\,dt\,dx\\
&=
\int_{L_p}^{\infty}\int_0^\infty
\left|\mathrm{Ai}(x+t)\right|^2\,dt\,dx.
\end{align*}
Then, using the standard Airy bound~\eqref{eq. Airy sharp bound} for $p$ sufficiently large so that $L_p\geq 1$, we obtain
\begin{align}\label{eq. estimate Airy tail}
    \lVert  (\mathbb{I}-\Pi_p)\mathcal A\rVert_2^2
    &\leq
    C\int_{L_p}^{\infty}\int_0^\infty
    e^{-\frac{2}{3}(x+t)^{3/2}}\,dt\,dx\\
    &\le C e^{-\frac{2}{3}L_p^{3/2}},
\end{align}
where $C$ has been adjusted in the last step. Recall that $L_p=p^{2/3}c_1\Omega(\lambda_D)$, with $c_1\Omega(\lambda_D)>\left(\frac{3}{2}k\right)^{2/3}$, for some positive constant $k$ independent of $p$ (see~\eqref{eq. def asymp L_p}). Then~\eqref{eq. trace norm bd1} and \eqref{eq. estimate Airy tail} imply that
\begin{equation}\label{eq. bd HS dif Ai}
   \left\|\mathcal K_{\mathrm{Airy}} -\Pi_p\mathcal K_{\mathrm{Airy}}\Pi_p\right\|_2 = \mathcal{O}\left( e^{-k\,p}\right).   
\end{equation}
Note that, in this case, the estimate relies only on the standard Airy bound on the tail, so the implicit constant in $\mathcal O\left( e^{-k\,p}\right)$ is independent even on $L$.

Next, note that the kernel of $\Pi_p\mathcal D_0\Pi_p$ is:
\[
\mathds{1}_{[L_0,L_p]}(x)\,
\widehat{C}_p\,p^{-2/3}
\bigl(
\mathrm{Ai}(x)\mathrm{Ai}'(y)
+
\mathrm{Ai}'(x)\mathrm{Ai}(y)
\bigr)\,
\mathds{1}_{[L_0,L_p]}(y).
\]
Then, Eq.~\eqref{eq. uniform bound P, D0} and the ideal property of the Hilbert--Schmidt norm give
\begin{equation}\label{eq. bd HS dif D0}
\left\|\Pi_p\mathcal D_0\Pi_p\right\|_2 \le \|\Pi_p\|_{op}  \left\|\mathcal D_0\right\|_2 \|\Pi_p\|_{op}=\mathcal O\left(p^{-2/3}\right),
\end{equation}
with the implicit constant in $\mathcal O\left( p^{-2/3}\right)$ independent of $L_0$, for $L_0>-L$.

Finally, we write
\begin{align*}
\left\|\widetilde{\mathcal R}^{(E)}\right\|_2^2&=
\int_{L_0}^{\infty}\int_{L_0}^{\infty}
\left|\widetilde R(x,y)\right|^2\,dx\,dy\\
& \le\int_{-L}^{\infty}\int_{-L}^{\infty}
\left|\widetilde R(x,y)\right|^2\,dx\,dy\\
& =\int_{-L}^{0}\int_{-L}^{0}
\left|\widetilde R(x,y)\right|^2\,dx\,dy+
\int_{-L}^{0}\int_{-L}^{\infty}\left|\widetilde R(x,y)\right|^2\,dx\,dy+
\int_{-L}^{\infty} \int_{-L}^{0}\left|\widetilde R(x,y)\right|^2\,dx\,dy+
\int_{0}^{\infty}\int_{0}^{\infty}
\left|\widetilde R(x,y)\right|^2\,dx\,dy,
\end{align*}
where $\widetilde R$ is defined as zero on $[L_0,\infty]^2\setminus[L_0,L_p)^2$. Then the bounds for $\widetilde R(x,y)$ in~\eqref{eq. bounds R tilde prop} yield
\begin{equation}\label{eq. bd HS norm R^E}
\left\|\widetilde{\mathcal R}^{(E)}\right\|_2
=\mathcal{O}\left(p^{-1}\right),
\end{equation}
with the implicit constant in $\mathcal O(p^{-1})$ independent of $L_0$, for all $L_0>-L$.

Therefore, Eq.~\eqref{eq. bound HS PD}, together with the estimates \eqref{eq. bd HS dif Ai}, \eqref{eq. bd HS dif D0} and \eqref{eq. bd HS norm R^E}, show that
\begin{equation}\label{eq. order PD}
   \lVert\mathcal{P}\mathcal{D}\rVert_2 =\mathcal{O}\left(p^{-2/3}\right),
\end{equation}
where the implicit constant in $\mathcal O(p^{-2/3})$ depends on the fixed value $L$ and is uniform for all $L_0>-L$.

Finally, we need to estimate the term $ \lVert \mathcal{P}\mathcal{D}-\mathcal{P}\mathcal{D}_0\rVert_2 $. To this end, first note that
\begin{align}\label{eq. bound HS PD-PD0}
    \lVert \mathcal{P}\mathcal{D}-\mathcal{P}\mathcal{D}_0\rVert_2 &\le \lVert \mathcal{P} \lVert_{op} \lVert \mathcal{D}-\mathcal{D}_0\rVert_2\nonumber\\
    &
    =\lVert \mathcal{P} \lVert_{op} \,\lVert \Pi_p\, \left(\mathcal{K}_{Airy}+ \mathcal{D}_0\right) \,\Pi_p +\widetilde{\mathcal R}^{(E)}- \mathcal{K}_{Airy}-\mathcal{D}_0\rVert_2\nonumber\\
   & \le\lVert\mathcal{P}\rVert_{op} \left(
   \|\Pi_p\mathcal K_{\mathrm{Airy}}\Pi_p -\mathcal K_{\mathrm{Airy}}\|_2 + \|\Pi_p\,\mathcal D_0\,\Pi_p -\mathcal{D}_0\|_2 + \| \widetilde{\mathcal R}^{(E)}\|_2\right).
\end{align}
We have already estimates for the first factor (see~\eqref{eq. uniform bound P, D0}), and for the first and third terms from the parentheses in~\eqref{eq. bound HS PD-PD0}. It therefore remains to estimate the second term:
\begin{align}\label{eq. D0 - Do truncated}
    \|\Pi_p\,\mathcal D_0\,\Pi_p-\mathcal D_0\|_2
    &=
    \|\mathcal D_0-\Pi_p\,\mathcal D_0\,\Pi_p\|_2\nonumber\\
    &=
    \|\mathcal D_0-\mathcal D_0\Pi_p
    +\mathcal D_0\Pi_p-\Pi_p\mathcal D_0\Pi_p\|_2\nonumber\\
    &=
    \|\mathcal D_0(\mathbb I-\Pi_p)
    +(\mathbb I-\Pi_p)\mathcal D_0\Pi_p\|_2\nonumber\\
    &\leq
    \|\mathcal D_0(\mathbb I-\Pi_p)\|_2
    +\|(\mathbb I-\Pi_p)\mathcal D_0\Pi_p\|_2.
\end{align}
Now, note that the kernel of
\(\mathcal D_0(\mathbb I-\Pi_p)\) is:
\[
\widehat{C}_p\,p^{-2/3}
\Bigl(
\mathrm{Ai}(x)\mathrm{Ai}'(y)
+
\mathrm{Ai}'(x)\mathrm{Ai}(y)
\Bigr)
\mathds{1}_{(L_p,\infty)}(y).
\]
Then
\begin{equation}\label{eq. D0 I-P}
\lVert\mathcal D_0(\mathbb I-\Pi_p)\rVert_2
\le
|\widehat{C}_p|p^{-2/3}
\Bigl(
\lVert\mathrm{Ai}\rVert_2
\lVert(\mathbb I-\Pi_p)\mathrm{Ai}'\rVert_2+
\lVert\mathrm{Ai}'\rVert_2
\lVert(\mathbb I-\Pi_p)\mathrm{Ai}\rVert_2
\Bigr).
\end{equation}
Similarly, the kernel of
\((\mathbb I-\Pi_p)\mathcal D_0\Pi_p\) is:
\[
\widehat{C}_p\,p^{-2/3}
\mathds{1}_{(L_p,\infty)}(x)
\Bigl(
\mathrm{Ai}(x)\mathrm{Ai}'(y)
+
\mathrm{Ai}'(x)\mathrm{Ai}(y)
\Bigr)
\mathds{1}_{[L_0,L_p]}(y).
\]
Then
\begin{align}\label{eq. I-P D0 P_p}
\lVert(\mathbb I-\Pi_p)\mathcal D_0\Pi_p\rVert_2
&\le
|\widehat{C}_p|p^{-2/3}
\Bigl(
\lVert(\mathbb I-\Pi_p)\mathrm{Ai}\rVert_2
\lVert \Pi_p\mathrm{Ai}'\rVert_2 +
\lVert(\mathbb I-\Pi_p)\mathrm{Ai}'\rVert_2
\lVert \Pi_p\mathrm{Ai}\rVert_2
\Bigr)\nonumber\\
&\le |\widehat{C}_p|p^{-2/3}
\Bigl(
\lVert(\mathbb I-\Pi_p)\mathrm{Ai}\rVert_2
\lVert \mathrm{Ai}'\rVert_2 +
\lVert(\mathbb I-\Pi_p)\mathrm{Ai}'\rVert_2
\lVert \mathrm{Ai}\rVert_2
\Bigr).
\end{align}
Hence, \eqref{eq. D0 - Do truncated}, together with the estimates \eqref{eq. D0 I-P} and \eqref{eq. I-P D0 P_p}, yield
\begin{equation}\label{eq:D0-HS-trunc-bound}
\begin{split}
\lVert \Pi_p\mathcal D_0\Pi_p-\mathcal D_0\rVert_2
\le
2|\widehat{C}_p|p^{-2/3}
\Bigl(
\lVert\mathrm{Ai}\rVert_2
\lVert(\mathbb I-\Pi_p)\mathrm{Ai}'\rVert_2+
\lVert\mathrm{Ai}'\rVert_2
\lVert(\mathbb I-\Pi_p)\mathrm{Ai}\rVert_2
\Bigr).
\end{split}
\end{equation}
Here, note that
\[
\|(\mathbb I-\Pi_p)\mathrm{Ai}\|_2^2
=
\int_{L_p}^{\infty}|\mathrm{Ai}(x)|^2\,dx,
\qquad
\|(\mathbb I-\Pi_p)\mathrm{Ai}'\|_2^2
=
\int_{L_p}^{\infty}|\mathrm{Ai}'(x)|^2\,dx.
\]
Moreover, the standard Airy bounds in \eqref{eq. Airy sharp bound} imply, in particular, that for all $x$ sufficiently large:
\[
|\mathrm{Ai}(x)|
\le Cx^{-1/4}e^{-\frac23x^{3/2}},
\qquad
|\mathrm{Ai}'(x)|
\le Cx^{1/4}e^{-\frac23x^{3/2}}.
\]
Consequently,
\begin{equation}\label{eq. (I-P_p)A A' bound}
\|(\mathbb I-\Pi_p)\mathrm{Ai}\|_2
\le
C L_p^{-1/2}e^{-\frac23L_p^{3/2}},
\qquad
\|(\mathbb I-\Pi_p)\mathrm{Ai}'\|_2
\le
C e^{-\frac23L_p^{3/2}}.
\end{equation}
Thus \eqref{eq:D0-HS-trunc-bound} together with \eqref{eq. (I-P_p)A A' bound} imply that there exist a constant $C$ such that
\begin{align*}
\|\Pi_p\mathcal D_0\Pi_p-\mathcal D_0\|_2 & \le
 C\, |\widehat C_p|p^{-2/3}
e^{-\frac23 L_p^{3/2}}.
\end{align*}
Again, since $L_p=p^{2/3}c_1\Omega(\lambda_D)>0$, with $c_1\Omega(\lambda_D)>\left(\frac{3}{2}k\right)^{2/3}$, for some positive constant $k$ independent of $p$ (see~\eqref{eq. def asymp L_p}), and since $\widehat{C}_p$ is uniformly bounded, then after adjusting constants we have
\begin{equation}\label{eq. bd trace class dif D0}
\|\Pi_p\mathcal D_0\Pi_p-\mathcal D_0\|_2=  \mathcal{O}\left(
e^{-k\,p}\right).
\end{equation}
Here the implicit constant in $\mathcal O\left( e^{-k\,p}\right)$ is uniform for all $L_0>-L$. Moreover, in this case $C$ is independent of $L$ as we are using only bounds on the tail of the Airy function.

Therefore, Eq.~\eqref{eq. bound HS PD-PD0}, together with the estimates \eqref{eq. bd HS dif Ai}, \eqref{eq. bd HS norm R^E} and \eqref{eq. bd trace class dif D0} give 
\begin{equation}\label{eq. Gohberg key 1}
    \lVert \mathcal{P}\mathcal{D}-\mathcal{P}\mathcal{D}_0\rVert_2 
    = \mathcal{O}\left(p^{-1}\right),
\end{equation}
 with the implicit constant in $\mathcal O\left(p^{-1}\right)$ depending the fixed value $L$ and uniform for all $L_0>-L$.

Finally, our claimed estimate~\eqref{eq. key 1 Carleman} follows from \eqref{eq. Thm Gohberg} and the estimates \eqref{eq. order PD_0}, \eqref{eq. order PD} and \eqref{eq. Gohberg key 1}.
\vspace{0.2cm}

\textbf{Proof of the estimate~\eqref{eq. key 2 Carleman}}. We begin by defining the operator $\mathcal Q$ by means of the identity:
\[
\mathcal P
=
\left(\mathbb I-\mathcal K_{\mathrm{Airy}}\right)^{-1}
=
\mathbb I+\mathcal Q.
\]
Note that
\[
\mathcal Q =\mathcal P-\mathbb I=
\mathcal P\mathcal K_{\mathrm{Airy}}.
\]
Since \(\mathcal P\) is bounded and
$\mathcal K_{\mathrm{Airy}}$ is trace class,  the ideal property implies that the operator $\mathcal Q$ is trace class. In particular $\mathcal Q$ is Hilbert--Schmidt. Therefore,
\begin{align}\label{eq:trace-PD-D0-bound}
\left|
\operatorname{Tr}\!\left(
\mathcal P(\mathcal D-\mathcal D_0)
\right)
\right|
&=
\left|
\operatorname{Tr}(\mathcal D-\mathcal D_0)
+
\operatorname{Tr}\!\left(
\mathcal Q(\mathcal D-\mathcal D_0)
\right)
\right|
\nonumber\\
&\leq
\left|\operatorname{Tr}(\mathcal D-\mathcal D_0)\right|
+\left|\operatorname{Tr}\!\left(
\mathcal Q(\mathcal D-\mathcal D_0)
\right)
\right| \nonumber\\
&\leq
\left|\operatorname{Tr}(\mathcal D-\mathcal D_0)\right|
+\| Q(\mathcal D-\mathcal D_0)
\|_1 \nonumber\\
&\leq
\left|\operatorname{Tr}(\mathcal D-\mathcal D_0)\right|
+
\|\mathcal Q\|_2
\|\mathcal D-\mathcal D_0\|_2.
\end{align}
Since $\mathcal Q$ acts on $L^2([L_0,\infty))$, with $L_0>-L$, then $\|\mathcal Q\|_2$ can be uniformly bounded for all  $L_0>-L$, with the bound depending only on $L$. Also, recall that (see~\eqref{eq. bound HS PD-PD0}, together with \eqref{eq. bd HS dif Ai}, \eqref{eq. bd HS norm R^E}, and  \eqref{eq. bd trace class dif D0}):
\begin{equation}\label{eq. order HS D-D0}
\|\mathcal D-\mathcal D_0\|_2
\le \|  \Pi_p\mathcal K_{\mathrm{Airy}}\Pi_p -\mathcal K_{\mathrm{Airy}}  \|_2+\|\Pi_p\mathcal D_0\Pi_p
-\mathcal D_0 \|_2+\| \widetilde{\mathcal R}^{(E)} \|_2 = \mathcal O\left(p^{-1}\right).
\end{equation}

On the other hand, since $\mathcal D-\mathcal D_0$ is trace class, the diagonal trace formula gives
\begin{align}
\operatorname{Tr}(\mathcal D-\mathcal D_0)
&=
-\int_{L_p}^{\infty}
K_{\mathrm{Airy}}(x,x)\,dx
-\int_{L_p}^{\infty}
D_0(x,x)\,dx
+
\int_{L_0}^{L_p}
\widetilde R(x,x)\,dx.
\label{eq:trace-D-D0-decomposition}
\end{align}
Here, the standard Airy bounds \eqref{eq. Airy sharp bound}  imply, in particular, that there exist a constant $C>0$ such that, for $x>0$,
\[
\left|K_{\mathrm{Airy}}(x,x)\right|
\le C e^{-\frac{4}{3}x^{3/2}}.
\]
Since $L_p=p^{2/3}c_1\Omega(\lambda_D)>0$, with $c_1\Omega(\lambda_D)>\left(\frac{3}{4}k\right)^{2/3}$, for some positive constant $k$ independent of $p$, then 
\begin{equation}\label{eq. trace K_Airi}
\int_{L_p}^{\infty}
K_{\mathrm{Airy}}(x,x)\,dx
= \mathcal{O}\left(e^{-kp}\right).
\end{equation}
Here, the implicit constant in $\mathcal{O}\left(e^{-kp}\right)$ is independent on $L$, as we are using the Airy bound on the tail.

Moreover,
\[
D_0(x,x)
= 2\,\widehat C_p\,p^{-2/3}
\mathrm{Ai}(x)\mathrm{Ai}'(x).
\]
Then
\begin{equation}\label{eq. trace D0}
\left|
\int_{L_p}^{\infty}D_0(x,x)\,dx
\right| \le
2\left|\widehat C_p\right|\,p^{-2/3}
\int_{L_p}^{\infty} \left|
\mathrm{Ai}(x)\mathrm{Ai}'(x)\right|\,dx
= \mathcal{O}\left(e^{-kp}\right),
\end{equation}
with $k$ as in \eqref{eq. trace K_Airi}, and the implicit constant $C$ in $\mathcal O\left(e^{-kp}\right)$ independent of $L_0$ (even of $L$). 

Finally, the bounds for $\widetilde R$ in~Proposition~\ref{proposition Khat to K_airy} imply that
\begin{align}\label{eq. order trace R tilde}
\left| \int_{L_0}^{L_p}\widetilde R(x,x)\,dx
\right| &\leq \int_{L}^{\infty} \left|\widetilde R(x,x)\right|\,dx\\
&= \int_{L}^{0} \left|\widetilde R(x,x)\right|\,dx+\int_{0}^{\infty} \left|\widetilde R(x,x)\right|\,dx  \\
&=\mathcal O\left(p^{-1}\right),
\end{align}
 with the implicit constant in $\mathcal O\left(p^{-1}\right)$ depending the fixed value $L$ and uniform for all $L_0>-L$.

Hence, using \eqref{eq. trace K_Airi} \eqref{eq. trace D0} and~\eqref{eq. order trace R tilde} in Eq.~\eqref{eq:trace-D-D0-decomposition}, we obtain
\begin{equation}\label{eq. trace D-D0}
\left|
\operatorname{Tr}(\mathcal D-\mathcal D_0)
\right| =
\mathcal O\left(p^{-1}\right).
\end{equation} 

Therefore, Eq.~\eqref{eq:trace-PD-D0-bound}, together with the estimates \eqref{eq. order HS D-D0} and \eqref{eq. trace D-D0} imply that
\[
\operatorname{Tr}\!\left(
\mathcal P(\mathcal D-\mathcal D_0)
\right)
=
\mathcal O\left(p^{-1}\right),
\]
with the implicit constant in $\mathcal O\left(p^{-1}\right)$ uniform for all $L_0>-L$. This proves the claimed estimate \eqref{eq. key 2 Carleman} and completes the proof of Proposition~\ref{lem:bound-D-D0}. 


\section{Painlev\'{e} II representation of the correction term}\label{Proof Prop 5.4}

In this section, we establish the connection between the correction term in the asymptotic expansion \eqref{eq:mres2} and the Painlev\'{e} II equation. To this end, let $s\in\mathbb{R}$ be fixed, and throughout this section, consider operators acting on $L^2([s,\infty))$. We recall the following quantities:
\begin{equation}\label{eq. def p q u v}
\begin{aligned}
q(x,s) &= \left((\mathbb{I}-\mathcal K_{\mathrm{Airy}})^{-1}\mathrm{Ai}\right)(x),
& p(x,s) &= \left((\mathbb{I}-\mathcal K_{\mathrm{Airy}})^{-1}\mathrm{Ai}'\right)(x),\\
u(s) &= \left\langle \mathrm{Ai},(\mathbb{I}-\mathcal K_{\mathrm{Airy}})^{-1}\mathrm{Ai}\right\rangle,
& v(s) &= \left\langle \mathrm{Ai},(\mathbb{I}-\mathcal K_{\mathrm{Airy}})^{-1}\mathrm{Ai}'\right\rangle,
\end{aligned}
\end{equation}
where $\langle\cdot,\cdot\rangle$ denotes the inner product on $L^2([s,\infty))$. The following proposition provides an explicit representation of $v(s)$ in terms of the Painlev\'{e} II function.

\begin{proposition}
Let $v(s)$ be defined in \eqref{eq. def p q u v}. Then,
\begin{equation}\label{eq: defn v}
v(s) =
\frac{1}{2}\left(\int_s^{\infty}q^2(t)\,dt\right)^2
-\frac{q^2(s)}{2},
\end{equation}
where $q(s)$ is the same solution of the Painlev\'{e} II equation specified in \eqref{eq. PII}.
\end{proposition}

\begin{remark}
Note that the Airy operator $\mathcal{K}_{\Airy}$ acts on $L^2([s,\infty))$, hence the $s$-dependence in the definitions of $p$, $q$, $u$ and $v$, enters through the domain $[s,\infty)$. Moreover, the inverse $(\mathbb I-\mathcal K_{\mathrm{Airy}})^{-1}$ exists, since the Fredholm determinant $\det(\mathbb I-\mathcal K_{\mathrm{Airy}}) = F_2(s)$ coincides with the Tracy--Widom c.d.f., which is strictly positive for every $s \in \mathbb{R}$ (see \cite{tracy1994level}).
\end{remark}

\begin{proof}
We split the proof of this proposition into several steps. First we obtain some key identities that will be used to establish differential equations for the functions $q$, $p$, $u$, and $v$ defined in \eqref{eq. def p q u v}, with respect to $s$. Then we prove that $q$ solves the Painlev\'{e} II equation.

\textbf{Key relations.} A straightforward differentiation of the Airy kernel yields:
\begin{equation}\label{eq:key}
\partial_x K_{\textnormal{Airy}}(x,y) + \partial_y K_{\textnormal{Airy}}(x,y) = -\mathrm{Ai}(x) \mathrm{Ai}(y).
\end{equation}
Now, let $R(x,y)$ denote the resolvent kernel, that is, the kernel of the operator $\mathcal{R}=(\mathbb I-\mathcal K_{\Airy})^{-1}\mathcal K_{\Airy}$. Then, a useful identity follows:
\begin{equation}\label{eq:resol}
(\mathbb I-\mathcal K_{\Airy})^{-1}
=(\mathbb I-\mathcal K_{\Airy})^{-1}(\mathbb I-\mathcal K_{\Airy}+\mathcal K_{\Airy})
=\mathbb I+(\mathbb I-\mathcal K_{\Airy})^{-1}\mathcal K_{\Airy}
=\mathbb I+\mathcal{R}.
\end{equation}
The definition of the resolvent operator is equivalently rewritten as:
\begin{align}\label{e:resolvent}
    \mathcal{R}= \mathcal K_{\Airy}+ \mathcal K_{\Airy} \, \mathcal{R}
\end{align}
where the operator $\mathcal K_{\Airy} \, \mathcal{R}$ is defined as follows
\begin{align*}
\left( \mathcal K_{\Airy} \, \mathcal{R} \, f\right)(x)&=\int_s^{\infty} K_{\Airy}(x,y) \left( \mathcal{R} f\right)(y) \, dy= \int_s^\infty \left[ \int_s^\infty K_{\Airy}(x,y) R(y,z) dy\right] f(z) dz.
\end{align*}
Therefore, equation \eqref{e:resolvent} at the level of the kernels, implies:
\begin{align}\label{e:rst}
R(x,s)=K_{\Airy}(x,s)+\int_s^{\infty} K_{\Airy}(x,y) R(y,s) \, dy.
\end{align}
This is the identity that will be used below. Equivalently:
\begin{equation}\label{eq:key1}
R(x,s)=(\mathbb I-\mathcal K_{\Airy})^{-1}\left(K_{\Airy}(x,s)\right).
\end{equation}
Next, since $p$ and $q$ both depend on $x$ and $s$, we derive expressions for their partial derivatives with respect to both variables. To do this, we rewrite the equations for $p$ and $q$ as follows:
\begin{align*}
(\mathbb I-\mathcal K_{\Airy})\,q(x,s) = \mathrm{Ai}(x), \quad (\mathbb I-\mathcal K_{\Airy})\,p(x,s) = \mathrm{Ai}'(x),
\end{align*}
which, for all $x>s$, can be written in the integral form:
\begin{align}\label{eq:pandq}
q(x,s)=\mathrm{Ai}(x)+ \int_s^{\infty} K_{\Airy}(x,y) \, q(y,s) \, dy, \quad p(x,s)=\mathrm{Ai}'(x)+ \int_s^{\infty} K_{\Airy}(x,y) \, p(y,s) \, dy.
\end{align}
Hence,
\begin{align}\label{e:eqs}
\frac{\partial q(x,s)}{\partial x} &= \mathrm{Ai}'(x) + \int_s^{\infty} \partial_x K_{\Airy}(x,y) \, q(y,s) \, dy, \\
\frac{\partial q(x,s)}{\partial s} &= -K_{\Airy}(x,s) q(s,s) + \int_s^{\infty} K_{\Airy}(x,y) \, \frac{\partial q(y,s)}{\partial s} \, dy,\\
\frac{\partial p(x,s)}{\partial x} &= x \mathrm{Ai}(x) + \int_s^{\infty} \partial_x K_{\Airy}(x,y) \, p(y,s) \, dy, \\
\frac{\partial p(x,s)}{\partial s} &= -K_{\Airy}(x,s) p(s,s) + \int_s^{\infty} K_{\Airy}(x,y) \, \frac{\partial p(y,s)}{\partial s} \, dy
\end{align}
where in the third equation we used the fact that the Airy function satisfies $\mathrm{Ai}''(x)= x \mathrm{Ai}(x)$. Using relation \eqref{eq:key1}, we can rewrite the second equation in \eqref{e:eqs} as follows:
\begin{align*}
(\mathbb I-\mathcal K_{\Airy}) \left( \frac{\partial q(x,s)}{\partial s} \right) &= -K_{\Airy}(x,s) q(s,s) \\
\frac{\partial q(x,s)}{\partial s}&= (\mathbb I-\mathcal K_{\Airy})^{-1} \left(  -K_{\Airy}(x,s) q(s,s) \right)\\
\frac{\partial q(x,s)}{\partial s}&=- q(s,s) (\mathbb I-\mathcal K_{\Airy})^{-1} \left( K_{\Airy}(x,s) \right)\\
\frac{\partial q(x,s)}{\partial s}&=- q(s,s) R(x,s).
\end{align*}
Similarly, we can write the last equation in \eqref{e:eqs} as:
\[\frac{\partial p(x,s)}{\partial s}=- p(s,s) R(x,s).\]
We recall these two expressions below, which we will use later:
\begin{equation}\label{eq. part p,q, wrt s}
\frac{\partial q(x,s)}{\partial s}=- q(s,s) R(x,s),\quad
\frac{\partial p(x,s)}{\partial s}=- p(s,s) R(x,s). 
\end{equation}
\textbf{Differential equation  for $q$.} Subtracting the equation for $p(x,s)$ in \eqref{eq:pandq}, from the equation for $\frac{\partial q}{\partial x}$  in \eqref{e:eqs}, and using the relation \eqref{eq:key}, we obtain
\begin{align*}
\frac{\partial q(x,s)}{\partial x} - p(x,s) &= -\int_s^{\infty} \left( \partial_y K_{\Airy}(x,y) + \mathrm{Ai}(x)\, \mathrm{Ai}(y) \right) q(y,s) \, dy - \int_s^{\infty} K_{\Airy}(x,y) p(y,s) \, dy.
\end{align*}
Integration by parts on the first integral of the right hand side yields:
\begin{align*}
\frac{\partial q(x,s)}{\partial x} - p(x,s) &= K_{\Airy}(x,s) q(s,s) + \int_s^{\infty} K_{\Airy}(x,y) \left( \frac{\partial q(y,s)}{\partial y} - p(y,s) \right) \, dy - \mathrm{Ai}(x) u(s).
\end{align*}
We can rewrite the last relation into the following Fredholm integral form:
\begin{align*}
(\mathbb I-\mathcal K_{\Airy}) \left( \frac{\partial q(x,s)}{\partial x} - p(x,s) \right) &= K_{\Airy}(x,s) q(s,s) - \mathrm{Ai}(x) u(s).
\end{align*}
Applying $(\mathbb I-\mathcal K_{\Airy})^{-1}$ to both sides and using the definition of $q$ together with relation \eqref{eq:key1}, we obtain:
\begin{align}\label{eq:eq15}
\frac{\partial q(x,s)}{\partial x} - p(x,s) &= q(s,s) R(x,s) - q(x,s) u(s).
\end{align}
Next, we add equation \eqref{eq:eq15} to the equation for $\frac{\partial q}{\partial s}$ in \eqref{e:rst}, and combine relations \eqref{e:rst}-- \eqref{eq. part p,q, wrt s} to obtain:
\begin{align}\label{eq:eq16}
\frac{\partial q}{\partial  s}(x,s) + \frac{\partial q(x,s)}{\partial x} - p(x,s)
&=-K_{\Airy}(x,s) q(s,s) + \int_s^{\infty} K_{\Airy}(x,y) \, \frac{\partial q(y,s)}{\partial s} \, dy+q(s,s) R(x,s) - q(x,s) u(s) \nonumber \\
&=- q(x,s) u(s).
\end{align}
Since $q$ is a function of both $x$ and $s$, if we define $q(s):=q(s,s)$, then 
\begin{equation}\label{eq:eq17}
q'(s) = \left[\frac{\partial q}{\partial s}(x,s) + \frac{\partial q}{\partial x}(x,s)\right]_{x=s}.    
\end{equation}
Therefore, evaluating equation \eqref{eq:eq16} at $x=s$ and using \eqref{eq:eq17}, we obtain
\[\boxed{q'(s)=p(s)-q(s) u(s)}\]
where $p(s):=p(s,s)$.

\textbf{Differential equation  for $p$.} We proceed similarly to the derivation of the differential equation for $q$. Subtracting $x$ times the equation for $q(x,s)$ in \eqref{eq:pandq} from the equation for $\frac{\partial p}{\partial x}$  in \eqref{e:eqs}, and using the relation \eqref{eq:key}, we obtain
\begin{align*}
\frac{\partial p(x,s)}{\partial x} - x q(x,s) &= - \int_s^{\infty} \left( \partial_y K_{\Airy}(x,y) + \mathrm{Ai}(x) \, \mathrm{Ai}(y)\right) \, p(y,s) \, dy - x \int_s^{\infty} K_{\Airy}(x,y) \, q(y,s) \, dy.
\end{align*}
Integration by parts on the first integral of the right hand side yields
\begin{align*}
\frac{\partial p(x,s)}{\partial x} - x q(x,s) &= K_{\Airy}(x,s) \, p(s,s) - \mathrm{Ai}(x) v(s) + \int_s^{\infty} K_{\Airy}(x,y) \frac{\partial p (y,s)}{\partial y} \, dy - x \int_s^{\infty} K_{\Airy}(x,y) \, q(y,s) \, dy.
\end{align*}
To rewrite the previous relation in a Fredholm integral form, we first add and subtract the quantity $\int_s^{\infty} K_{\Airy}(x,y) y \, q(y,s) \, dy$ on the right-hand side. This yields
\begin{align*}
\frac{\partial p(x,s)}{\partial x} - x q(x,s) &= K_{\Airy}(x,s) \, p(s,s) - \mathrm{Ai}(x) v(s) + \int_s^{\infty} K_{\Airy}(x,y) \left( \frac{\partial p(y,s)}{\partial y} - y \, q(y,s) \right) \, dy \\ & + \int_s^{\infty} K_{\Airy}(x,y) (y-x) \, q(y,s) \, dy.
\end{align*}
Moving the first integral of the right-hand side on the other side, and using the definition of the Airy kernel on the last term, the previous equation can be expressed in the following Fredholm integral form:
\begin{equation*}
(\mathbb I-\mathcal K_{\Airy}) \left( \frac{\partial p(x,s)}{\partial x} - x q(x,s) \right) = K_{\Airy}(x,s) \, p(s,s) - \mathrm{Ai}(x) v(s) + \int_s^{\infty} \left( -\mathrm{Ai}(x) \mathrm{Ai}'(y) + \mathrm{Ai}'(x) \mathrm{Ai}(y) \right) \, q(y,s) \, dy. 
\end{equation*}
Using the definitions of $u$ and $v$, together with the self-adjointness of $(\mathbb I-\mathcal K_{\Airy})^{-1}$, we obtain
\begin{equation*}
(\mathbb I-\mathcal K_{\Airy}) \left( \frac{\partial p(x,s)}{\partial x} - x q(x,s) \right) = K_{\Airy}(x,s) \, p(s,s) - 2 \mathrm{Ai}(x) v(s) + \mathrm{Ai}'(x) u(s). 
\end{equation*}
Applying $(\mathbb I-\mathcal K_{\Airy})^{-1}$ to both sides and using relation \eqref{eq:key1}, this gives
\begin{align}\label{eq: dif p0}
\frac{\partial p(x,s)}{\partial x} - x q(x,s) &= (\mathbb I-\mathcal K_{\Airy})^{-1} \left( K_{\Airy}(x,s) \, p(s,s) - 2 \mathrm{Ai}(x) v(s) + \mathrm{Ai}'(x) u(s) \right)\nonumber \\
&= -2v(s) q(x,s) + u(s) p(x,s) + p(s,s) R(x,s). 
\end{align}
Again, since $p$ is a function of both $x$ and $s$, if we define $p(s):=p(s,s)$, then 
\begin{align}\label{eq:eq18}
p'(s) = \left(\frac{\partial p}{\partial s}(x,s) + \frac{\partial p}{\partial x}(x,s) \right)\Bigg|_{x=s}.
\end{align}
Therefore, evaluating \eqref{eq: dif p0} and  \eqref{eq. part p,q, wrt s} at $x=s$, and using \eqref{eq:eq18}, we obtain:
\[\boxed{p'(s) = s\, q(s) - 2 v(s) q(s) + u(s) p(s)}\]
\textbf{Differential equation  for $u$.} Recall that $u$ is defined as
\[u(s) = \left\langle \mathrm{Ai}, (I-\mathcal K_{\mathrm{Airy}})^{-1} \mathrm{Ai} \right\rangle= \int_s^{\infty} \mathrm{Ai}(x) \, q(x,s) \, dx.\]
Taking the derivative with respect to $s$ and using relation \eqref{eq. part p,q, wrt s}, we obtain
\begin{align*}
u'(s) &= -\mathrm{Ai}(s) q(s,s) + \int_s^{\infty} \mathrm{Ai}(x) \, \frac{\partial q(x,s)}{\partial s} \, dx \\
 &= -\mathrm{Ai}(s) q(s,s) - \int_s^{\infty} \mathrm{Ai}(x) q(s,s) \, R(x,s) \, dx \\
 &= -q(s,s) \left( \mathrm{Ai}(s)+ \int_s^{\infty} \mathrm{Ai}(x) \, R(x,s) \, dx \right).
\end{align*}
Now,  using the definition of $q$ and relation \eqref{eq:resol}, we obtain
\begin{align}\label{eq:eq19}
q(x,s) = \big((\mathbb I-\mathcal K_{\mathrm{Airy}})^{-1} \mathrm{Ai}\big)(x) = (\mathbb I+\mathcal R) \mathrm{Ai}(x).
\end{align}
Since $\mathcal K_{\Airy}$ is self-adjoint, so is
$(\mathbb I-\mathcal K_{\Airy})^{-1}$, and consequently the resolvent
operator $\mathcal R$ is self-adjoint. Hence its kernel is symmetric. Therefore, we have:
\begin{align*}
\boxed{u'(s) = -q^2(s)}
\end{align*}
\textbf{Differential equation  for $v$.} Recall that $v$ is defined as
\[v(s) =\left\langle \mathrm{Ai}, (\mathbb I-\mathcal K_{\mathrm{Airy}})^{-1} \mathrm{Ai}' \right\rangle= \int_s^{\infty} \mathrm{Ai}(x) \, p(x,s) \, dx.\]
Differentiating with respect to $s$ and using \eqref{eq. part p,q, wrt s}, we get
\begin{align*}
v'(s) &= -\mathrm{Ai}(s) p(s,s) + \int_s^{\infty} \mathrm{Ai}(x) \, \frac{\partial p(x,s)}{\partial s} \, dx \\
 &= -\mathrm{Ai}(s) p(s,s) - \int_s^{\infty} \mathrm{Ai}(x) p(s,s) \, R(x,s) \, dx \\
 &= -p(s,s) \left( \mathrm{Ai}(s)+ \int_s^{\infty} \mathrm{Ai}(x) \, R(x,s) \, dx \right)\\
 &= -p(s) \, q(s).
\end{align*}
 Lastly, using equation \eqref{eq:eq19}, this gives the differential equation  for $v$:
 \[\boxed{v'(s) = -p(s) \, q(s)}\]
\textbf{The Painlev\'{e} II equation.} We have obtained the following differential equations:
\begin{align*}
q'(s) &= p(s)-q(s) u(s),\\
p'(s) &= s q(s) - 2 v(s) q(s) + u(s) p(s),\\
u'(s) &= -q^2(s) ,\\
v'(s) &= -p(s) \, q(s).
\end{align*}
Differentiating the first one once again with respect to $s$ gives:
\begin{align}\label{eq:eqPVII}
q'' &= p' -q' u - q u'\\
&= sq + pu -2 q v - q' u + q^3\\
&= sq + q^3 + q(u^2-2v).
\end{align}
Here, in the second equality we used the differential equations for $p'$ and $u'$, while in the third equality we used the differential equation for $q'$. Next, we claim that 
\[u^2-2v = q^2.\]
Indeed, taking the derivative of the expression $u^2-2v$, and using the differential equations for $u'$, $v'$ and $q'$, we have
\begin{align*}
\left( u^2-2v \right)' &= 2u u'-2v'\\
&= 2q(p-u q)\\
&= 2q \, q'.
\end{align*}
Integrating gives
\begin{equation}\label{eq:conC}
u^2-2v = q^2 + C.
\end{equation}
Since the Airy function and its derivative decay exponentially as $s\to\infty$, the definitions of $q$, $u$, and $v$ imply that
\[
q(s)\to0,\qquad u(s)\to0,\qquad v(s)\to0.
\]
Therefore, letting $s\to\infty$ in \eqref{eq:conC} gives $C=0$. Furthermore, recalling that $u'(s)=-q^2(s)$ and
$u(s)\to 0$ as $s\to\infty$, we obtain
\[
u(s)=\int_s^\infty q^2(t)\,dt.
\]
Substituting this expression into \eqref{eq:conC}, with
$C=0$, yields
\[
v(s)=\frac12\left(\int_s^\infty q^2(t)\,dt\right)^2
-\frac{q^2(s)}{2},
\]
which establishes the claimed representation of $v(s)$. Finally, substituting \eqref{eq:conC} into equation \eqref{eq:eqPVII}, we get PV-II:
\[q''=sq + 2 q^3.\]
 \end{proof}


\section{Conclusion}\label{s:conclusion}

In this work, we investigate the asymptotic distribution of the largest eigenvalue $\theta_{\max}$ of the random matrix ${\mathbf J} = ({\mathbf A}+{\mathbf B})^{-1} {\mathbf B}$, where ${\mathbf A}$ and ${\mathbf B}$ are standard complex Wishart matrices with parameter $m$ and $n$, respectively. The high-dimensional regime is potentially unbalanced in the sense that, as $p \to \infty$, $m$ may grow faster than $p$, whereas $n$ grows proportionally to $p$. The asymptotic behavior of $\theta_{\max}$ is characterized by means of a Riemann-Hilbert (RH) problem for Jacobi polynomials with modified Jacobi weights. We establish the type 2 Tracy-Widom asymptotic distribution of $\theta_{\max}$. This provides an alternative, RH-based, proof of the same claim in Johnstone \cite{johnstone:2008} for the balanced regime, and establishes a conjecture in Jiang \cite{jiang:2009} for the unbalanced regime. In addition, we obtain an explicit correction of order \(p^{-2/3}\) to the Tracy–Widom approximation, with a remainder of order \(\mathcal O\left(p^{-1}\right)\), and show the exponential decay, as a function of $p$, of the upper tail of the distribution of $\theta_{\max}$.

This paper leads to several open problems, of which we mention three. First, it naturally points to the investigation of the asymptotic distribution of the smallest eigenvalue of the unbalanced complex Jacobi ensemble (cf.\ Lyu and Chen \cite{lyu:chen:2021}). Second, there is the ``extreme opposite" unbalanced regime where $n \gg p$ and $m \approx p$. Third, the case $\beta = 2$ is a natural first step in the study of the large Jacobi ensembles for general $\beta$. Hence, with a view towards a broad set of applications, our work points to the study of the extreme eigenvalues of the unbalanced \textit{real} Jacobi ensemble.

\section*{Acknowledgments}
G.D. is partially supported by the award NSF DMS 2515732. J.L. and A.G. are supported by a PIMS postdoctoral research fellowship. K.M. is partially supported by the award NSF DMS 2554750.

\appendix 

\section{Technical results}\label{s:technical_results}


\subsection[Auxiliary results for a away from the edge of the spectrum]{Auxiliary results for $a$ away from the endpoint $\beta_p$}

In this subsection, we present useful lemmas that will be used to study the eigenvalue statistics when $\varphi_p^{-1}(a) \in (\lambda_D, \varphi_p^{-1}(1))$.

\begin{lemma}\label{lem:lemma1}
    Let $h$ be a real-valued function such that $h'$ is continuous on the interval $[a, b]$, and let $x, y$ be two points on $(a,b)$ with $x<y$. Consider 
    \begin{equation*}        A(x,y)=\frac{h(x)-h(y)}{x-y}.
    \end{equation*}
Then, $A$ is uniformly bounded when $x \rightarrow y$.
\end{lemma}
\begin{proof}
Since $h'$ is continuous on $[a, b]$, then there exists $M >0$ such that $h'(\lambda)\le M$, for all $\lambda \in [a, b]$. Set $\lambda = t y + (1-t) x$, then $| h'(\lambda) | \le M$, for all $t \in [0,1]$. Hence, we have:
\begin{align*}
  |A(x,y)| &= \left|\frac{h(x)-h(y)}{x - y}\right|\\
           &= \left|\frac{1}{y-x} \int_0^1 \frac{d}{dt} h(ty + (1-t)x) \, dt\right|\\
           &= \left|\frac{1}{y-x} \int_0^1 h'(t y + (1-t)x) (y-x) \, dt\right|\\
           &\le M.
\end{align*}
\end{proof}

\begin{lemma}\label{lem:lemma2}
For $\lambda_1, \lambda_2 \in (\lambda_D, \varphi_p^{-1}(1))$, there exists a constant $k > 0$ such that
\[ \widetilde{K}_p(\lambda_1,\lambda_2)=  \mathcal{O}(e^{-k p}).\]
\end{lemma}

\begin{proof}
Recall from \eqref{eq: Ytilde outside discs} that for $\lambda \in (\lambda_D, \varphi_p^{-1}(1))$, we have:
\begin{subequations}\label{eq:eq67}
\begin{gather}
\widetilde{Y}_{11}(\lambda) =\frac{i}{2}e^{p\,g(\varphi_p(\lambda))}\left(-i f(\lambda) E_{11}(\lambda) + \widetilde{f}(\lambda) E_{12}(\lambda)\right),
\end{gather}
\begin{gather}
\widetilde{Y}_{21}(\lambda) =\frac{i}{2}e^{p(g(\varphi_p(\lambda))-\ell)}\left(-i f(\lambda) E_{21}(\lambda) + \widetilde{f}(\lambda) E_{22}(\lambda)\right),
\end{gather}
\end{subequations}
where $f(\lambda)$, $\widetilde{f}(\lambda)$ are given in \eqref{eq:eqfftil}. Hence, from relation \eqref{eq:eqlastrelation} we have:
{\small
\begin{align*}
 \widetilde{K}_p (\lambda_1, \lambda_2)  & = -\frac{\sqrt{w_p(\varphi_p(\lambda_1)) w_p(\varphi_p(\lambda_2))}\,\,e^{p(g(\varphi_p(\lambda_1))+g(\varphi_p(\lambda_2))-\ell)}}{4 \pi i(\beta_p - \alpha_p)} \left(\frac{\widetilde{Y}_{11}(\lambda_1)\widetilde{Y}_{21}(\lambda_2)-\widetilde{Y}_{11}(\lambda_2)\widetilde{Y}_{21}(\lambda_1)}{\lambda_1-\lambda_2}\right)\\
 &=-\frac{\sqrt{w_p(\varphi_p(\lambda_1))w_p(\varphi_p(\lambda_2))}\,\,e^{p(g(\varphi_p(\lambda_1))+g(\varphi_p(\lambda_2))-\ell)}}{8 \pi i} \left( \frac{A(\lambda_1,\lambda_2) + B(\lambda_1,\lambda_2) + C(\lambda_1,\lambda_2) + D(\lambda_1,\lambda_2)}{\lambda_1-\lambda_2} \right)
\end{align*}
}
where the functions $A$, $B$, $C$ and $D$ are given explicitly by the following expressions:
{\small
\begin{align*}
A(\lambda_1,\lambda_2) &= f(\lambda_1) E_{11}(\lambda_1) \Big( f(\lambda_1) E_{21}(\lambda_1) - f(\lambda_2) E_{21}(\lambda_2)\Big) - f(\lambda_1) E_{21}(\lambda_1) \Big( f(\lambda_1) E_{11}(\lambda_1) - f(\lambda_2) E_{11}(\lambda_2)\Big)\\
B(\lambda_1,\lambda_2) &= i f(\lambda_1) E_{11}(\lambda_1) \left( \widetilde{f}(\lambda_1) E_{22}(\lambda_1) - \widetilde{f}(\lambda_2) E_{22}(\lambda_2) \right) -i \widetilde{f}(\lambda_1) E_{22}(\lambda_1) \Big( f(\lambda_1) E_{11}(\lambda_1) - f(\lambda_2) E_{11}(\lambda_2) \Big)\\
C(\lambda_1,\lambda_2) &= -i \widetilde{f}(\lambda_1) E_{12}(\lambda_1) \Big( f(\lambda_1) E_{21}(\lambda_1) - f(\lambda_2) E_{21}(\lambda_2) \Big) - i f(\lambda_1) E_{21}(\lambda_1) \left( \widetilde{f}(\lambda_1) E_{12}(\lambda_1) - \widetilde{f}(\lambda_2) E_{12}(\lambda_2) \right)\\
D(\lambda_1,\lambda_2) &= \widetilde{f}(\lambda_1) E_{22}(\lambda_1) \left( \widetilde{f}(\lambda_1) E_{12}(\lambda_1) - \widetilde{f}(\lambda_2) E_{12}(\lambda_2) \right) - \widetilde{f}(\lambda_1) E_{12}(\lambda_1) \left( \widetilde{f}(\lambda_1) E_{22}(\lambda_1) - \widetilde{f}(\lambda_2) E_{22}(\lambda_2) \right).
\end{align*}
}
Rearranging terms in the expressions for $A$, $B$, $C$, and $D$, we obtain:
{\small
\begin{align*}
A(\lambda_1,\lambda_2) + B(\lambda_1,\lambda_2) + C(\lambda_1,\lambda_2) + D(\lambda_1,\lambda_2) &= \left( i f(\lambda_1) E_{11}(\lambda_1) - \widetilde{f}(\lambda_1) E_{12}(\lambda_1) \right) \left( \psi(\lambda_1) - \psi(\lambda_2) \right) \nonumber \\ &+ \left( \widetilde{f}(\lambda_1) E_{22}(\lambda_1) -i f(\lambda_1) E_{21}(\lambda_1) \right) \left( \widetilde{\psi}(\lambda_1) - \widetilde{\psi}(\lambda_2) \right),
\end{align*}
}
where we have used the notation
\begin{gather*}
\psi(\lambda) = -i f(\lambda) E_{21}(\lambda) + \widetilde{f}(\lambda) E_{22}(\lambda), \quad \widetilde{\psi}(\lambda) = -i f(\lambda) E_{11}(\lambda) + \widetilde{f}(\lambda) E_{12}(\lambda).
\end{gather*}
This implies that
{\small
\begin{align}\label{eq:eq610}
\widetilde{K}_p (\lambda_1, \lambda_2) &=- \frac{\sqrt{w_p(\varphi_p(\lambda_1))w_p(\varphi_p(\lambda_2))}\,\,e^{p(g(\varphi_p(\lambda_1)) + g(\varphi_p(\lambda_2)) - \ell)}}{8 \pi i} \nonumber \\
& \times \left( \left( i f(\lambda_1) E_{11}(\lambda_1) - \widetilde{f}(\lambda_1) E_{12}(\lambda_1) \right) \frac{ \psi(\lambda_1) - \psi(\lambda_2)}{\lambda_1 - \lambda_2} +  \left( \widetilde{f}(\lambda_1) E_{22}(\lambda_1) -i f(\lambda_1) E_{21}(\lambda_1) \right) \frac{ \widetilde{\psi}(\lambda_1) - \widetilde{\psi}(\lambda_2)}{\lambda_1 - \lambda_2} \right).
\end{align}
}
Since $f(\lambda)$ and $\widetilde{f}(\lambda)$ are continuous on $(\lambda_D, \varphi_p^{-1}(1))$, and all entries of the error matrix $E$ are analytic, it follows that the functions 
\[ i f(\lambda) E_{11}(\lambda) - \widetilde{f}(\lambda) E_{12}(\lambda), \quad \widetilde{f}(\lambda) E_{22}(\lambda) - i f(\lambda) E_{21}(\lambda)\]
are continuous and consequently bounded on $(\lambda_D, \varphi_p^{-1}(1))$. Moreover, since $\psi$ and $\widetilde{\psi}$ satisfy the requirements of Lemma~\ref{lem:lemma1}, then the functions
\[\frac{\psi(\lambda_1) - \psi(\lambda_2)}{\lambda_1 - \lambda_2}, \quad \frac{\widetilde{\psi}(\lambda_1) - \widetilde{\psi}(\lambda_2)}{\lambda_1 - \lambda_2}\]
are bounded when $\lambda_1 \rightarrow \lambda_2$, for $\lambda_1, \lambda_2 \in (\lambda_D,\varphi_p^{-1}(1))$. Therefore, the entire expression inside the parenthesis in \eqref{eq:eq610} is bounded. Furthermore, note that:

\begin{align*}
\sqrt{w_p(\varphi_p(\lambda_1))w_p(\varphi_p(\lambda_2))}&\,\,e^{p(g(\varphi_p(\lambda_1)) + g(\varphi_p(\lambda_2))-\ell)} \\ 
&= \sqrt{\varphi_p(\lambda_1)^{p \, b_{\infty}}(1-\varphi_p(\lambda_1))^{a_{\infty} \, p^{1+\nu}}\,\varphi_p(\lambda_2)^{p \, b_{\infty}}(1-\varphi_p(\lambda_2))^{a_{\infty} \, p^{1+\nu}}}e^{p(g(\varphi_p(\lambda_1)) + g(\varphi_p(\lambda_2)) -\ell)}\\
&= e^{p\Big(g(\varphi_p(\lambda_1)) + g(\varphi_p(\lambda_2)) -\ell+\frac{b_{\infty}}{2}\ln \varphi_p(\lambda_1) +\frac{b_{\infty}}{2}\ln \varphi_p(\lambda_2)+\frac{a_{\infty} \, p^{\nu}}{2}\ln(1-\varphi_p(\lambda_1))+\frac{a_{\infty} \, p^{\nu}}{2}\ln(1-\varphi_p(\lambda_1))\Big)}\\
&=e^{p\Big(g(\varphi_p(\lambda_1)) + \frac{b_{\infty}}{2}\ln \varphi_p(\lambda_1)+\frac{a_{\infty} \, p^{\nu}}{2}\ln(1-\varphi_p(\lambda_1))-\frac{\ell}{2}+g(\varphi_p(\lambda_2))+\frac{b_{\infty}}{2}\ln \varphi_p(\lambda_2)+\frac{a_{\infty} \, p^{\nu}}{2}\ln(1-\varphi_p(\lambda_2))-\frac{\ell}{2}\Big)}.
\end{align*}
Since $g$ is analytic on the interval $ (1 + \epsilon, \varphi_p^{-1}(1)) $, equation \eqref{e:eq9} implies that, for any $\lambda$ in this interval we have
\begin{equation*}
2g(\varphi_p(\lambda)) \leq - b_{\infty} \log \varphi_p(\lambda) - a_{\infty} \, p^{\nu} \log (1 - \varphi_p(\lambda)) + \ell.
\end{equation*}
This shows that the coefficient of $p$ in the previous exponential term is negative. Therefore, the kernel $\widetilde{K}_p(\lambda_1,\lambda_2)$ decays exponentially as $p \to \infty$.
\end{proof}


\subsection[Auxiliary results for a near the edge of the spectrum]{
Auxiliary results for $a$ near the endpoint $\beta_p$}

In this subsection, we present useful lemmas that will be used to study the eigenvalue statistics when $a$ is near $\beta_p$.

\begin{lemma}\label{lem:lemma3}
If $\lambda_1 \in (\lambda_D, \varphi_p^{-1}(1))$, and $\lambda_2,\lambda_3 \in (\varphi_p^{-1}(a), \lambda_D-\delta)$, where $\delta>0$, then there exists a constant $k>0$ such that
\[\widetilde{K}_{p}(\lambda_1,\lambda_2) \widetilde{K}_p(\lambda_3,\lambda_1) = \mathcal{O}(e^{-k p}).\]
\end{lemma}

\begin{proof}

Recall from \eqref{eq:eqYtilde11} that for $\lambda \in (\varphi_p^{-1}(a), \lambda_D)$, we have:
\begin{align}
    \widetilde{Y}_{11}(\lambda) = e^{\frac{p \ell}{2}} \Big(E_{11}(\lambda)\left(Q \mathcal{A}\right)(\lambda)_{11} + E_{12}(\lambda)  \left(Q \mathcal{A}\right)(\lambda)_{21}  \Big)\,e^{-\frac{i\pi}{6}}w_p^{-1/2}(\varphi_p(\lambda)) \nonumber \\
    \widetilde{Y}_{21}(\lambda) = e^{-\frac{p \ell}{2}} \Big(E_{21}(\lambda) \left(Q \mathcal{A}\right)(\lambda)_{11} + E_{22}(\lambda)  \left(Q \mathcal{A}\right)(\lambda)_{21}  \Big)\,e^{-\frac{i\pi}{6}}w_p^{-1/2}(\varphi_p(\lambda)) 
\end{align} 
and the functions $\left(Q \mathcal{A}\right)(\lambda)_{11}$ and $\left(Q \mathcal{A}\right)(\lambda)_{21} $ are given in equation \eqref{eq QA}, and are bounded away from $\lambda=1$. Also, from \eqref{eq: Ytilde outside discs}, for $\lambda \in (\lambda_D, \varphi_p^{-1}(1))$, we have
\begin{align*}
\widetilde{Y}_{11}(\lambda) = e^{p g(\varphi_p(\lambda))} \left(-i E_{11}(\lambda) f(\lambda) + E_{12}(\lambda) \widetilde{f}(\lambda) \right),\\
\widetilde{Y}_{21}(\lambda) = e^{p (g(\varphi_p(\lambda))- \ell)} \left( -i E_{21}(\lambda) f(\lambda) + E_{22}(\lambda) \widetilde{f}(\lambda) \right),
\end{align*}
where the functions $f$ and $\widetilde{f}$ are given in \eqref{eq:eqfftil}, and are bounded away from $\lambda  =1$. Using the above equations, we have:
\begin{align}\label{eq:eqs612}
\widetilde{Y}_{11}^2(\lambda) = \mathcal{O}(e^{2p g(\varphi_p(\lambda))}), & \quad \lambda \in (\lambda_D,\varphi_p^{-1}(1)), \nonumber \\
\widetilde{Y}_{21}^2(\lambda) =  \mathcal{O}(e^{2p (g(\varphi_p(\lambda)) - \ell)}), & \quad \lambda \in (\lambda_D,\varphi_p^{-1}(1)), \nonumber \\
\widetilde{Y}_{21}^2(\lambda) = \mathcal{O}(e^{- p \ell}), & \quad \lambda \in (\varphi_p^{-1}(a),\lambda_D), \nonumber \\ 
\widetilde{Y}_{11}(\lambda) \widetilde{Y}_{21}(\lambda) = \mathcal{O}(e^{p(2 g(\varphi_p(\lambda)) - \ell)}), & \quad \lambda \in (\lambda_D,\varphi_p^{-1}(1)). 
\end{align}
Moreover, since $g$ is analytic on $(\lambda_D, \varphi_p^{-1}(1))$, inequality \eqref{e:eq9} implies that
\begin{equation}\label{eq:eqineq}
2 g(\varphi_p(\lambda)) + b_{\infty} \log(\varphi_p(\lambda)) + a_{\infty} p^{\nu} \log(1-\varphi_p(\lambda)) - \ell < 0, \quad \lambda \in (\lambda_D, \varphi_p^{-1}(1)).
\end{equation}
Let $\lambda_1 \in (\lambda_D, \varphi_p^{-1}(1))$ and $\lambda_2,\lambda_3 \in (\varphi_p^{-1}(a), \lambda_D-\delta)$, then
\begin{equation}\label{eq:eq614}
\begin{split}
\widetilde{K}_{p}(\lambda_1,\lambda_2) &\widetilde{K}_p(\lambda_3,\lambda_1) = \frac{w_p(\varphi_p(\lambda_1)) \sqrt{w_p(\varphi_p(\lambda_2)) w_p(z(\lambda_3))}}{\pi^2 (\lambda_1-\lambda_2)(\lambda_1-\lambda_3)}\\
& \times \left[\widetilde{Y}_{11}(\lambda_1)\widetilde{Y}_{21}(\lambda_2)-\widetilde{Y}_{11}(\lambda_2)\widetilde{Y}_{21}(\lambda_1)\right]\left[\widetilde{Y}_{11}(\lambda_3)\widetilde{Y}_{21}(\lambda_1)-\widetilde{Y}_{11}(\lambda_1)\widetilde{Y}_{21}(\lambda_3)\right].
\end{split}
\end{equation}

Using the definitions of $\widetilde{Y}_{11},\,\widetilde{Y}_{21}$ on the corresponding interval, and equations \eqref{eq:eqs612}, we obtain that the following holds for each of the terms in the right-hand side of \eqref{eq:eq614}:
\begin{align*}
w_{p}(\varphi_p(\lambda_1)) \widetilde{Y}_{11}(\lambda_1) \widetilde{Y}_{21}(\lambda_1) \widetilde{Y}_{21}(\lambda_2)\widetilde{Y}_{11}(\lambda_3) & = q_1(\lambda_1,\lambda_2,\lambda_3) \, e^{p[2 g(\varphi_p(\lambda_1)) + b_{\infty} \log(\varphi_p(\lambda_1)) + a_{\infty} p^{\nu} \log(1-\varphi_p(\lambda_1)) - \ell]}  \\
w_{p}(\varphi_p(\lambda_1))  \widetilde{Y}_{11}(\lambda_1) \widetilde{Y}_{21}(\lambda_1) \widetilde{Y}_{11}(\lambda_2) \widetilde{Y}_{21}(\lambda_3) & =  q_2(\lambda_1,\lambda_2,\lambda_3) \, e^{p[2 g(\varphi_p(\lambda_1)) + b_{\infty} \log(\varphi_p(\lambda_1)) + a_{\infty} p^{\nu} \log(1-\varphi_p(\lambda_1)) - \ell]}\\
w_{p}(\varphi_p(\lambda_1))  \widetilde{Y}_{11}^2(\lambda_1) \widetilde{Y}_{21}(\lambda_2) \widetilde{Y}_{21}(\lambda_3) & = q_3(\lambda_1,\lambda_2,\lambda_3) \, e^{p[2 g(\varphi_p(\lambda_1)) + b_{\infty} \log(\varphi_p(\lambda_1)) + a_{\infty} p^{\nu} \log(1-\varphi_p(\lambda_1)) - \ell]}\\
w_{p}(\varphi_p(\lambda_1))  \widetilde{Y}_{21}^2(\lambda_1) \widetilde{Y}_{11}(\lambda_2) \widetilde{Y}_{11}(\lambda_3) & = q_4(\lambda_1,\lambda_2,\lambda_3) \, e^{p[2 g(\varphi_p(\lambda_1)) + b_{\infty} \log(\varphi_p(\lambda_1)) + a_{\infty} p^{\nu} \log(1-\varphi_p(\lambda_1)) - \ell]}
\end{align*}
where $q_j(\lambda_1,\lambda_2,\lambda_3)$, for $j=1, \cdots 4$ are bounded functions. Since  $\lambda_1$ is bounded away from $\lambda_2,\lambda_3$, the above identities imply that we can rewrite \eqref{eq:eq614} as
\[\widetilde{K}_{p}(\lambda_1,\lambda_2) \widetilde{K}_p(\lambda_3,\lambda_1) = F(\lambda_1,\lambda_2,\lambda_3) \,e^{p[2 g(\varphi_p(\lambda_1)) + b_{\infty} \log(\varphi_p(\lambda_1)) + a_{\infty} p^{\nu} \log(1-\varphi_p(\lambda_1)) - \ell]} ,\]
where $F(\lambda_1,\lambda_2,\lambda_3)$ is bounded. The desired result follows from the last expression combined with inequality \eqref{eq:eqineq}.
\end{proof}

\begin{lemma}\label{Prop truncating kernel}
For $\varphi_p^{-1}(a)\in (1, \lambda_D)$, there exists a constant $k>0$ such that
\[\det\left( \mathbb{I} - \widetilde{\mathcal{K}}_{p} \right)_{L^2([\varphi_p^{-1}(a),\varphi_p^{-1}(1)])} = \det\left( \mathbb{I} - \widetilde{\mathcal{K}}_{p} \right)_{L^2([\varphi_p^{-1}(a),\lambda_D])}+ \mathcal{O}(e^{-k p}).\]
\end{lemma}

\begin{proof}
By definition
\begin{align*}
\det& \left( 1 - \widetilde{\mathcal{K}}_{p} \right)_{L^2([\varphi_p^{-1}(a),\varphi_p^{-1}(1)])}= 1 + \sum_{j=1}^{p} \frac{(-1)^{j}}{j!} \int_{\varphi_p^{-1}(a)}^{{\varphi_p^{-1}(1)}} \cdots \int_{\varphi_p^{-1}(a)}^{{\varphi_p^{-1}(1)}} \det_{i,k=1,2,\ldots,j} \left[ \widetilde{K}_{p} \left( \lambda_i, \lambda_k \right) \right] d^j\lambda \\
& = 1 + \sum_{j=1}^{p} \frac{(-1)^{j}}{j!} \int_{\varphi_p^{-1}(a)}^{{\lambda_D}} \cdots \int_{\varphi_p^{-1}(a)}^{{\lambda_D}} \det_{i,k=1,2,\ldots,j} \left[ \widetilde{K}_{p} \left( \lambda_i, \lambda_k \right) \right] d^j\lambda \\
& \hspace{2cm}+ \sum_{j=1}^{p} \frac{(-1)^{j}}{j!}\left(\sum_{l=1}^{2^j-1}  \int_{\gamma_1}^{{\rho_1}} \cdots \int_{\gamma_j}^{{\rho_j}} \det_{i,k=1,2,\ldots, j} \left[ \widetilde{K}_{p} \left(\lambda_i, \lambda_k \right) \right] d^j\lambda\right).
\end{align*}
In the second equality, for each $j$ we have split the $j$ integrals over the interval $(\varphi_p^{-1}(a), \varphi_p^{-1}(1))$ into one integral over $(\varphi_p^{-1}(a), \lambda_D)$ and one integral over $(\lambda_D, \varphi_p^{-1}(1))$. Thus, the summation over $l$ corresponds to the $2^j-1$ remaining $j$-folded integrals, where limits of integration are $\gamma_i\in\{\varphi_p^{-1}(a),\lambda_D\}$, and $\rho_i\in\{\lambda_D,\varphi_p^{-1}(1)\}$, with the condition that at least one of the $\gamma_i$'s equals $\lambda_D$. Hence, we can write 
\begin{equation}\label{eq: det(1-k) truncated}
\begin{split}
    \det\left( \mathbb{I} - \widetilde{\mathcal{K}}_{p} \right)_{L^2([\varphi_p^{-1}(a),\varphi_p^{-1}(1)])}&= \det\left( \mathbb{I} - \widetilde{\mathcal{K}}_{p} \right)_{L^2([\varphi_p^{-1}(a),\lambda_D])} \\&+\sum_{j=1}^{p} \frac{(-1)^{j}}{j!}\left(\sum_{l=1}^{2^j-1}  \int_{\gamma_1}^{{\rho_1}} \cdots \int_{\gamma_j}^{{\rho_j}} \det_{i,k=1,2,\ldots, j} \left[ \widetilde{K}_{p} \left(\lambda_i, \lambda_k \right) \right] d^j\lambda\right).
\end{split}
\end{equation}
It remains to show that the sum in \eqref{eq: det(1-k) truncated} decays exponentially as $p$ grows. Indeed, Lemma \ref{lem:lemma2} implies that if all the variables $\lambda_k \in (\lambda_D, \varphi_p^{-1}(1))$, then the resulting $j$-folding integral over $(\lambda_D, \varphi_p^{-1}(1))$ decays exponentially with $p$. We claim that exponential decay also occurs when there is at least one of the $j$-integrals over $(\lambda_D, \varphi_p^{-1}(1))$. Without loss of generality, suppose that the integral in $\lambda_1$-variable is over $(\lambda_D, \varphi_p^{-1}(1))$. First we focus on the integral
\begin{align*}
&\int_{\lambda_D}^{{\varphi_p^{-1}(1)}}  \det_{i,k=1,2,\ldots, j} \left[ \widetilde{K}_{p} \left( \lambda_i, \lambda_k \right) \right]\, d \lambda_1 = \int_{\lambda_D}^{\varphi_p^{-1}(1)}   \det \begin{vmatrix}
\widetilde{K}_p(\lambda_1,\lambda_1) & \cdots & \widetilde{K}_p(\lambda_1,\lambda_j)\\
\vdots & \cdots & \vdots \\
\widetilde{K}_p(\lambda_j,\lambda_1) & \cdots & \widetilde{K}_p(\lambda_j,\lambda_j)
\end{vmatrix} \, d \lambda_1\\ \\[2pt]
&= \int_{\lambda_D}^{\varphi_p^{-1}(1)}\,  \sum_{i=1}^{j} (-1)^{1+i} \widetilde{K}_p(\lambda_1,\lambda_i) \det \begin{vmatrix}
\widetilde{K}_p(\lambda_1,\lambda_1) & \cdots & \widetilde{K}_p(\lambda_1,\lambda_j)\\
\vdots & \cdots & \vdots \\
\widetilde{K}_p(\lambda_j,\lambda_1) & \cdots & \widetilde{K}_p(\lambda_j,\lambda_j)
\end{vmatrix}^{(1,i)} \, d \lambda_1\\
& = \int_{\lambda_D}^{\varphi_p^{-1}(1)} \left[\widetilde{K}_p(\lambda_1,\lambda_1)\det_{l,k=1,2,\ldots, j} [\widetilde{K}(\lambda_l,\lambda_k)]^{(1,1)}\right.\\
&\hspace{2cm} \left.+\left(\sum_{i=2}^{j} (-1)^{i+1} \widetilde{K}_p(\lambda_1,\lambda_i)\sum_{k=2}^{j}(-1)^{k+1}\widetilde{K}_p(\lambda_k,\lambda_1)\det_{l,k=1,2,\ldots, j} [\widetilde{K}(\lambda_l,\lambda_k)]^{(k,i)}\right)\right]\, d \lambda_1\\
& =\int_{\lambda_D}^{\varphi_p^{-1}(1)} \left[\widetilde{K}_p(\lambda_1,\lambda_1)\det_{l,k=1,2,\ldots, j} [\widetilde{K}_p(\lambda_l,\lambda_k)]^{(1,1)}\right.\\
&\hspace{2cm} \left.+\sum_{i=2}^{j} \sum_{k=2}^{j}(-1)^{k+i}\widetilde{K}_p(\lambda_1,\lambda_i)\widetilde{K}_p(\lambda_k,\lambda_1)\det_{l,k=1,2,\ldots, j} [\widetilde{K}_p(\lambda_l,\lambda_k)]^{(k,i)}\right] \,d \lambda_1\,.
\end{align*}
Since $\lambda_1 \in (\lambda_D,\varphi_p^{-1}(1))$, Lemma \ref{lem:lemma2} gives exponential decay for the first term. To prove the exponential decay for the double sum note that if either of the integrals in $\lambda_i$ or in $\lambda_k$ are over $(\lambda_D, \varphi_p^{-1}(1))$, then Lemma \ref{lem:lemma2} implies exponential decay of the term $\widetilde{K}_p(\lambda_1,\lambda_i)\widetilde{K}_p(\lambda_k,\lambda_1)$. On the other hand, if both integrals are over $(\varphi_p^{-1}(a),\lambda_D)$, then we can write: 
\begin{align*}
    \int_{\lambda_D}^{\varphi_p^{-1}(1)} \int_{\varphi_p^{-1}(a)}^{\lambda_D} \int_{\varphi_p^{-1}(a)}^{\lambda_D}&\widetilde{K}_p(\lambda_1,\lambda_i)\widetilde{K}_p(\lambda_k,\lambda_1) \, d \lambda_i\, d \lambda_k\,d \lambda_1\,= \int_{\lambda_D}^{\varphi_p^{-1}(1)} \int_{\varphi_p^{-1}(a)}^{\lambda_D-\delta} \int_{\varphi_p^{-1}(a)}^{\lambda_D-\delta}\widetilde{K}_p(\lambda_1,\lambda_i)\widetilde{K}_p(\lambda_k,\lambda_1) \, d \lambda_i\, d \lambda_k\,d \lambda_1\,\\
    &+ \int_{\lambda_D}^{\varphi_p^{-1}(1)} \int_{\lambda_D-\delta}^{\lambda_D} \int_{\varphi_p^{-1}(a)}^{\lambda_D-\delta}\widetilde{K}_p(\lambda_1,\lambda_i)\widetilde{K}_p(\lambda_k,\lambda_1) \, d \lambda_i\, d \lambda_k\,d \lambda_1\,\\
    &+ \int_{\lambda_D}^{\varphi_p^{-1}(1)} \int_{\varphi_p^{-1}(a)}^{\lambda_D-\delta} \int_{\lambda_D-\delta}^{\lambda_D}\widetilde{K}_p(\lambda_1,\lambda_i)\widetilde{K}_p(\lambda_k,\lambda_1) \, d \lambda_i\, d \lambda_k\,d \lambda_1\,\\
    &+ \int_{\lambda_D}^{\varphi_p^{-1}(1)} \int_{\lambda_D-\delta}^{\lambda_D} \int_{\lambda_D-\delta}^{\lambda_D}\widetilde{K}_p(\lambda_1,\lambda_i)\widetilde{K}_p(\lambda_k,\lambda_1) \, d \lambda_i\, d \lambda_k\,d \lambda_1\,.
\end{align*}
Here, Lemma \ref{lem:lemma3} implies the exponential decay of the first term on the right-hand side. For the remaining three terms, we utilize the kernel representation within a disk centered at $\lambda=1$ of smaller radius $\lambda_D-\delta$. In this setting, at least two of the integrals lie outside the disc, and thus Lemma \ref{lem:lemma2} applies to establish the exponential decay of these terms. With the claim proven, Lemma~\ref{Prop truncating kernel} follows immediately from \eqref{eq: det(1-k) truncated}.
\end{proof}

\begin{lemma}\label{lem: Bound ElamdaQ2-ElamdaQ1}
Let $H$ be a $2 \times 2$ matrix-valued function such that, for $i,j=1,2$, $H_{ij}=\mathcal{O}(1)$ and $H'_{ij}=\mathcal{O}(1)$ uniformly as $p \to \infty$. Then, there exists a constant $c$ such that:
\[\left| H_{ij}(\chi_p^{-1}(\xi_2))-H_{ij}(\chi_p^{-1}(\xi_1)) \right| \leq c \, p^{-2/3} \, \left| \xi_2-\xi_1\right|, \quad i,j=1,2.\]
\end{lemma}
\begin{proof}
For all $i,j$, we have:
\begin{align}\label{eq:estimate}
\left|H_{ij}(\chi_p^{-1}(\xi_2))-H_{ij}(\chi_p^{-1}(\xi_1)) \right| \leq \int_{\chi_p^{-1}(\xi_1)}^{\chi_p^{-1}(\xi_2)} \left| \frac{d}{ds} H_{ij}(s) \right| \, |ds|\leq c\left| \chi_p^{-1}(\xi_2)-\chi_p^{-1}(\xi_1)\right|.
\end{align}
It remains to estimate the difference $\left| \chi_p^{-1}(\xi_2)-\chi_p^{-1}(\xi_1)\right|$. Recall that for $|\xi|< c_1 p^{2/3}$:
\[\chi_p^{-1}(\xi) \equiv \lambda = \sum_{j=0}^\infty a_j \left( \frac{\xi}{c_1 p^{2/3}} \right)^j\] 
and therefore:
\[\left| \chi_p^{-1}(\xi_2)-\chi_p^{-1}(\xi_1)\right| \leq \left| \frac{a_1}{c_1}\right| p^{-2/3} \, |\xi_2-\xi_1|.\]
Replace this into \eqref{eq:estimate}, we get the desired estimate. Note that this estimate is independent of $i,j$.
\end{proof}


\section{Basic properties of complex random variables}\label{s:complex_RVs}

\begin{definition}
Let $\mathbf{X}$ and $\mathbf{Y}$ be two $\mathbb{C}^p$-valued random variables with finite second moments. The \textit{covariance} and the \textit{pseudo-covariance} between $\mathbf{X}$ and $\mathbf{Y}$ are given by, respectively,
$$
\mathbf{Cov}(\mathbf{X},\mathbf{Y}) =\mathbb{E} \big[(\mathbf{X} -\mathbb{E} \mathbf{X})\cdot (\mathbf{Y} -\mathbb{E} \mathbf{Y})^* \big], \quad 
\mathbf{PCov}(\mathbf{X},\mathbf{Y}) =\mathbb{E} \big[(\mathbf{X} -\mathbb{E} \mathbf{X}) \cdot (\mathbf{Y} -\mathbb{E} \mathbf{Y})^{\top} \big].
$$
\end{definition}

\begin{definition}
We say $\mathbf{X} \in \mathbb{C}^p$ is a \textit{complex Gaussian random vector} if $\Re \mathbf{X} \in \mathbb{R}^p$ and $\Im \mathbf{X} \in \mathbb{R}^p$ are Gaussian random vectors, where
$$
\mu :=\mathbb{E} \mathbf{X}, \quad \Sigma_{\textnormal{C}} := \mathbf{Cov}(\mathbf{X},\mathbf{X}), \quad \Sigma_{\textnormal{PC}} := \mathbf{PCov}(\mathbf{X},\mathbf{X}).
$$
 \textbf{Notation}: $\mathbf{X} \sim \mathbb{C} {\mathcal N}(\mu, \Sigma_{\textnormal{C}},\Sigma_{\textnormal{PC}} )$.  
 In particular, for $p = 1$, we say $Z$ is a \textit{standard complex Gaussian random variable} if $\Re Z, \Im Z \stackrel{\textbf{i.i.d.}}\sim {\mathcal N}(0,1/2)$. \textbf{Notation}: $Z \sim \mathbb{C} {\mathcal N}(0,1)$.
\end{definition}
In the following simple proposition, presented without proof, we recap the connection between independence, covariance and pseudo-covariance in the context of complex Gaussian random variables.
\begin{proposition}\label{p:Cov=0=PCov<=>indep}
Let $X$ and $Y$ be complex Gaussian random variables. Then, $\mathbf{Cov}(X,Y) = 0 = \mathbf{PCov} (X,Y)$ if and only if $X$ and $Y$ are independent.
\end{proposition}
The next proposition describes the fundamental properties of the discrete Fourier transform of a Gaussian white noise sequence.
\begin{proposition}\label{p:DFT_is_complex_Gaussian}
Let $\{X_t\}_{t = 1,\hdots,T} \stackrel{\mathbf{i.i.d.}}\sim {\mathcal N}(0,\sigma^2)$, where $\sigma^2 > 0$. Recall that the Fourier frequencies and the discrete Fourier transform are given by relations \eqref{e:omega-j} and \eqref{e:DFT}, respectively. Consider the indices
\begin{equation}\label{e:j,ell=1,...,n/2}
j, \ell = 1,\hdots,\lfloor T /2 \rfloor-1, \quad j \neq \ell.
\end{equation}
Then, 
\begin{itemize}
\item [$(i)$] 
\begin{equation}\label{e:ReF,ImF_i.i.d._N(0,sigma^2/2)}
\Re {\mathcal F}_{T}(\omega_j), \Im {\mathcal F}_{T}(\omega_j) \stackrel{\mathbf{i.i.d.}}\sim {\mathcal N}\Big(0,\frac{\sigma^2}{2}\Big).
\end{equation}
In particular, ${\mathcal F}_{T}(\omega_j)$ is a zero-mean, complex Gaussian random variable and 
\begin{equation}\label{e:|F|^2_sim_sigma^2Exp(1)}
|{\mathcal F}_{T}(\omega_j)|^2 \sim \sigma^2 \hspace{0.5mm}\textnormal{Exp}(1);
\end{equation}
\item [$(ii)$] the random variables ${\mathcal F}_{T}(\omega_j)$ and ${\mathcal F}_{T}(\omega_\ell)$ are independent.
\end{itemize} 

\end{proposition}
\begin{proof}
To show $(i)$, the linearity of the discrete Fourier transform implies that $\Re {\mathcal F}_{T}(\omega_j)$ and $\Im {\mathcal F}_{T}(\omega_j)$ are zero-mean, Gaussian random variables. Now note that
$$
\mathbb{E} \big[ \Re {\mathcal F}_{T}(\omega_j) \Im {\mathcal F}_{T}(\omega_j) \big]
= - \frac{\sigma^2}{T} \sum^{T}_{t=1}  \cos(t \omega_j ) \sin ( t \omega_j  )
= - \frac{\sigma^2}{2T} \sum^{T}_{t=1}  \sin ( 2 t \omega_j  )
$$ 
\begin{equation}\label{e:E(ReIm)=0}
= - \frac{\sigma^2}{2T} \cdot  \frac{\sin\big(\frac{2 \omega_j T}{2}\big) \sin\big(\frac{2\omega_j (T+1)}{2} \big)}{\sin\big(\frac{2\omega_j}{2}\big)} = 0,
\end{equation}
In \eqref{e:E(ReIm)=0}, the second to last equality is a consequence of formula \textbf{1.342.1} in Gradshteyn and Ryzhik \cite{gradshteyn:ryzhik:2014}. Also, the last equality holds because $\sin\big(\frac{2 \omega_j T}{2}\big) = 0$. Therefore, $\Re {\mathcal F}_{T}(\omega_j)$ and $\Im {\mathcal F}_{T}(\omega_j)$ are independent. In addition, for the indices \eqref{e:j,ell=1,...,n/2}, 
$$
\mathbf{Var} \Re {\mathcal F}_{T}(\omega_j) = \frac{\sigma^2}{T} \sum^{T}_{t=1} \cos^2(t \omega_j) = \frac{\sigma^2}{2} + \frac{\sigma^2}{T} \frac{\cos[(T+1) \omega_j] \sin(T \omega_j)}{2 \sin \omega_j}  = \frac{\sigma^2}{2}
$$ 
and
$$
\mathbf{Var} \Im {\mathcal F}_{T}(\omega_j) = \frac{\sigma^2}{T} \sum^{T}_{t=1} \sin^2(t \omega_j) = \frac{\sigma^2}{2} - \frac{\sigma^2}{T} \frac{\cos[(T+1) \omega_j] \sin(T \omega_j)}{2 \sin \omega_j}  = \frac{\sigma^2}{2}
$$ 
(see relations \textbf{1.351.1} and \textbf{1.351.2} in Gradshteyn and Ryzhik \cite{gradshteyn:ryzhik:2014}). This establishes \eqref{e:ReF,ImF_i.i.d._N(0,sigma^2/2)} and, hence, $(i)$.  

We now turn to $(ii)$. Note that
$$
\mathbf{Cov} \big({\mathcal F}_{T}(\omega_j), {\mathcal F}_{T}(\omega_\ell)\big)  =\mathbb{E} \big({\mathcal F}_{T}(\omega_j) \overline{{\mathcal F}_{T}(\omega_\ell)}\big) = \sigma^2 \sum^{T}_{t=1} e^{\Im t (\omega_\ell - \omega_j)} 
$$
\begin{equation}\label{e:Cov(F(omega-j,omega-ell)}
=: \sigma^2 \sum^{T}_{t=1} e^{\Im t \widetilde{\omega}}
= \sigma^2 e^{\Im  \widetilde{\omega}} \frac{1 - e^{\Im  \widetilde{\omega}T}}{1 - e^{\Im  \widetilde{\omega}}} = 0.
\end{equation}
In \eqref{e:Cov(F(omega-j,omega-ell)}, the last equality holds because, under condition \eqref{e:j,ell=1,...,n/2},  $ \widetilde{\omega}T = \frac{2 \pi}{T}(\ell - j) T \in 2 \pi \mathbb{Z}  $ and $e^{\Im  \widetilde{\omega}} \neq 1$. Likewise, 
$$
 \mathbf{PCov} \big({\mathcal F}_{T}(\omega_j), {\mathcal F}_{T}(\omega_\ell)\big)  =\mathbb{E} \big({\mathcal F}_{T}(\omega_j) {\mathcal F}_{T}(\omega_\ell)\big) = \sigma^2 \sum^{T}_{t=1} e^{-\Im t (\omega_\ell + \omega_j)} 
 $$
 \begin{equation}\label{e:PCov(F(omega-j,omega-ell)}
 =: \sigma^2 \sum^{T}_{t=1} e^{\Im t \breve{\omega}}
= \sigma^2 e^{\Im  \breve{\omega}} \frac{1 - e^{\Im  \breve{\omega}T}}{1 - e^{\Im  \breve{\omega}}} = 0.
\end{equation}
In \eqref{e:PCov(F(omega-j,omega-ell)}, the last equality holds because, under condition \eqref{e:j,ell=1,...,n/2},  $ \breve{\omega}T = - \frac{2 \pi}{T}(\ell + j) T \in 2 \pi \mathbb{Z}  $ and $e^{\Im  \breve{\omega}} \neq 1$ under \eqref{e:j,ell=1,...,n/2}.  By Proposition \ref{p:Cov=0=PCov<=>indep} and by \eqref{e:ReF,ImF_i.i.d._N(0,sigma^2/2)},  ${\mathcal F}_{T}(\omega_j)$ and ${\mathcal F}_{T}(\omega_\ell)$ are independent. In other words, $(ii)$ holds. $\Box$\\
\end{proof}

\begin{remark}
Note that, combined, relation \eqref{e:|F|^2_sim_sigma^2Exp(1)} and Proposition \ref{p:DFT_is_complex_Gaussian}, part $(ii)$, non-asymptotically retrieve the claim of Proposition 10.3.2, part $(i)$, in Brockwell and Davis \cite{brockwell:davis:1991}.
\end{remark}


\bibliographystyle{plain}
\bibliography{bio.bib}

\end{document}